\documentclass[12pt]{article}

\usepackage{amssymb,amsmath,bm}
\usepackage{mathrsfs}
\usepackage{stackengine}
\usepackage{constants,color}
\usepackage{enumerate}
\usepackage{appendix}
\usepackage{tikz}
\def\qed{\hfill \mbox{\rule{0.5em}{0.5em}}}
\usepackage{graphicx}
\usepackage{amssymb}
\usepackage{bookmark}
\usepackage{soul,xcolor}

\newcommand{\be}{\begin{equation}}
\newcommand{\ee}{\end{equation}}
\newcommand{\bes}{\begin{equation*}}
\newcommand{\ees}{\end{equation*}}
\newcommand{\ba}{\begin{aligned}}
\newcommand{\ea}{\end{aligned}}
\newcommand{\bi}{\begin{itemize}}
\newcommand{\ei}{\end{itemize}}

\newcommand{\EE}{\mathbb{E}}
\newcommand{\PP}{\mathbb{P}}

\newcommand{\RR}{\mathbb{R}}

\allowdisplaybreaks

\numberwithin{equation}{section}

\newtheorem{theorem}{Theorem}

\newtheorem{proposition}{Proposition}
\newtheorem{lemma}{Lemma}
\newtheorem{definition}{Definition}[section]

\title{Clonal Diversity at Cancer Recurrence}
\author{Kevin Leder, Zicheng Wang}

\begin{document}
\setstcolor{red}
\maketitle
\begin{abstract}
Despite initial success, cancer therapies often fail due to the emergence of drug-resistant cells. In this study, we use a mathematical model to investigate how cancer evolves over time, specifically focusing on the state of the tumor when it recurs after treatment. We use a two-type branching process to capture the dynamics of both drug-sensitive and drug-resistant cells. We analyze the clonal diversity of drug-resistant cells at the time of cancer recurrence, which is defined as the first time the population size of drug-resistant cells exceeds a specified proportion of the initial population size of drug-sensitive cells. We examine two clonal diversity indices: the number of clones and the Simpson's Index. We calculate the expected values of these indices and utilize them to develop statistical methods for estimating model parameters. Additionally, we examine these two indices conditioned on early recurrence in the special case of a deterministically decaying sensitive population, with the aim of addressing the question of whether early recurrence is driven by a single mutation that generates an unusually large family of drug-resistant cells (corresponding to a low clonal diversity), or if it is due to the presence of an unusually large number of mutations causing drug resistance (corresponding to a high clonal diversity). Our findings, based on both indices, support the latter possibility. Furthermore, we demonstrate that the time of cancer recurrence can serve as a valuable indicator of clonal diversity, offering new insights for the treatment of recurrent cancers.
\end{abstract}

\section{Introduction}\label{introduction}

Despite advancements in cancer therapies, tumor cells exhibit a remarkable ability for developing drug resistance, thereby limiting the effectiveness of treatment and leading to cancer recurrence \cite{Lackner2012}. Recurrent human cancers are well-documented to exhibit substantial intratumor heterogeneity. For example, analysis of genomic DNA from recurrent human malignant gliomas reveals a large number of somatic mutations following alkylating agent treatment \cite{Hunter2006}.
 
The level of clonal diversity, or the number of genetically distinct populations, at the time of cancer recurrence is crucial for clinical decision-making and treatment efficacy. For example, in the case of chronic myeloid leukemia, imatinib, a common treatment, can be rendered ineffective by a variety of distinct point mutations that confer drug resistance. Second-line agents, such as dasatinib and nilotinib, are effective against some mutations, but partial resistance can still occur \cite{Gambacorti2011, Burgess2006}. Therefore, understanding the clonal diversity of recurrent tumors has significant implications for guiding treatment strategies\footnote{While this research is motivated by cancer recurrence, the findings can be applied to other areas such as pest control, parasitic infection treatment, and treatment of other diseases caused by viruses and bacteria.}.

To investigate clonal diversity, we utilize a stylized two-type branching process model. In this model, we examine a population of drug-sensitive cancer cells that undergo continuous reduction in population size during therapy. These drug-sensitive cells accumulate driver mutations at a small, individual constant rate. Cells with driver mutations become resistant to the therapy, enabling them to (potentially) escape extinction and lead to cancer recurrence. This work builds upon previous research \cite{JK2013, JKJ2014, BP2020}. Foo and Leder (\cite{JK2013}) examine the pathway by which cancer cells escape treatment. They obtain a uniform in time approximation for the sample paths as the initial size of the tumor approaches infinity. Additionally, they study two important times in the course of cancer recurrence: (1) the time at which the total population size (drug-sensitive cells plus drug-resistant cells) begins to rebound, and (2) the first time at which the size of drug-resistant cells exceeds that of drug-sensitive cells (crossover time). In \cite{JKJ2014}, the authors extend the result of \cite{JK2013} by including random mutational fitness advantage. They consider a more general setting in which each mutation results in a mutant with a birth rate sampled randomly from a distribution. The authors obtain a functional central limit result for the mutant cell process and then establish a central limit theorem for the crossover time. In \cite{BP2020}, the authors investigate the large deviations of cancer recurrence timing. The authors obtain a convergence in probability result for the recurrence time and then apply a large deviations analysis to the event of early recurrence (a similar result is obtained for the crossover time). The authors also obtain the most likely number of mutant clones at cancer recurrence through optimizing the large deviations rate.

Our work is also related to a body of literature that explores the path to extinction of a biological group and its escape from extinction. Iwasa and coauthors (\cite{Iwasa2003}) use a multi-type branching process model to study the escape dynamics of a biological group from biomedical intervention. The authors obtain the probability of a successful escape under different scenarios, considering factors such as the number of point mutations required to confer resistance and the distribution of mutants before intervention. Jagers and co-authors (\cite{JKS2007a}) investigate the path to extinction of a subcritical Markov branching process. They obtain a convergence of finite dimensional distributions result for the path when the time is scaled to $\left[0,1\right]$, where $0$ denotes the starting time, and $1$ denotes the extinction time. They (\cite{JKS2007b}) further extend the result of \cite{JKS2007a} by considering a more general branching process. Sagitov and co-authors (\cite{Sagitov2009}) study the escape from the extinction of a Bienayme-Galton-Watson process. They obtain the limit process conditioned on successful escape as the mutation rate goes to zero. In \cite{Sehl2011}, the author employs a birth-death process to investigate the extinction times of cancer cells and normal cells in response to a therapy. Their findings can be used to evaluate a therapy's safety and efficacy. Lastly, Avanzini and Antal (\cite{avanzini2019cancer}) use a branching process model to study cancer recurrence resulting from latent metastases.

Another related stream of literature focuses on the intratumor heterogeneity induced by mutations in the tumor cell population. Previous studies (\cite{CG1986}, \cite{Michelson1989}, \cite{Komarova2006}, \cite{Haeno2007}, \cite{Durrett2010}, \cite{Durrett2011}, \cite{Keller2015}, \cite{Cheek2018}) have extensively examined the mutation of cancer cells and their resulting drug resistance. Of particular relevance to our work is the study by \cite{Durrett2011}, where the authors examine the intratumor heterogeneity of a tumor during its expansion. They investigate a multi-type branching process model in which each mutation results in a random, additive change in the cell's birth rate. They study both between-generation heterogeneity (where cells with the same number of mutations are grouped together as a generation) and within-generation heterogeneity in the first generation of cells (comprised of cells with only one mutation). For the latter part, they investigate two metrics to measure heterogeneity: (1) Simpson's Index and (2) the fraction of cells that belong to the largest clone. 

In this article, we investigate two indices of clonal diversity: the number of mutant clones and the Simpson's Index of mutant clones at cancer recurrence. In Section \ref{result_diversity_indices}, we obtain the limit of the expectation of the scaled number of mutant clones and the scaled Simpson's Index (at a deterministic time close to the cancer recurrence time) per Lemma \ref{expectation_without_conditioning} and Proposition \ref{simpson_without_conditioning} respectively. In Section \ref{Sec:estimator}, we use our asymptotic results for the number of clones and Simpson's Index to develop statistical methods for estimating model parameters. In Section \ref{sec:results_with_conditioning}, we study the number of mutant clones and the Simpson's Index of mutant clones conditioned on early recurrence in the special case of deterministically decaying sensitive population. In Proposition \ref{proposition_diversity}, we show that the distribution of the number of mutant clones at cancer recurrence conditioned on early recurrence stochastically dominates that without conditioning in the large population limit. In Lemma \ref{expectation_with_conditioning}, we obtain the limit of the expectation of the scaled number of mutant clones conditioned on early recurrence. We find that the expected number of clones conditioned on early recurrence is higher than that without conditioning. In Proposition \ref{period_clones_result}, we show that, conditioned on early recurrence, the number of clones generated in any given sufficiently small time period is concentrated at a larger number than would be expected without conditioning. In Theorem \ref{Simpson_conditioning_I}, we obtain the limit of the expectation of the scaled Simpson's Index of mutant clones conditioned on early recurrence. Simulation results indicate that the Simpson's Index conditioned on early recurrence is smaller than that without conditioning. All of our findings suggest that early recurrence is associated with higher clonal diversity. In particular, our results indicate that early recurrence is driven by a larger than expected number of mutations which leads to a more diverse resistant population.

The remainder of this paper is organized as follows. In Section \ref{model}, we describe our model and present important results from previous works. In Section \ref{result_diversity_indices}, we present results on diversity indices without conditioning. In Section \ref{Sec:estimator}, we provide a set of estimators for the model parameters. In Section \ref{sec:results_with_conditioning}, we present results on diversity indices conditioned on early recurrence. In Sections \ref{Sec:proof_result_diversity_indices} to \ref{Sec:proof_simpson_with_conditioning}, we present the proofs of our main results.


\section{Models and Previous Results}\label{model}

In this section, we describe a two-type branching process model, which has been employed in previous studies to examine the progression of cancer under treatment (see, for example, \cite{JK2013} and \cite{JKJ2014}).

Consider a subcritical birth-death process $\left(Z_0^n\left(t\right)\right)_{t\ge 0}$ with birth rate $r_0$, death rate $d_0$ and net growth rate $\lambda_0 = r_0-d_0 < 0$. $Z_0^n$ represents the population size of \textit{drug-sensitive cells} under a certain treatment. At time $0$, the initial population size is $Z_0^n\left(0\right)=n$. We assume that at time $t$, drug-sensitive cells give birth to a \textit{drug-resistant mutant} and a drug-sensitive cell at rate $Z_0^n\left(t\right) \mu n^{-\alpha}$ for $\alpha \in \left(0,1\right)$.
Each of these mutations results in the creation of a distinct clone (each mutation is distinct under the infinite sites approximation) which is modeled as a Yule process \footnote{If cell death is considered, our results can be applied to the skeleton subpopulation of cells whose descendants do not go extinct.} with birth rate $r_1=\lambda_1 > 0$, and death rate $d_1=0$. We denote this population by $\left(Z_1^n\left(t\right)\right)_{t\ge 0}$. Then $Z_1^n$ is a supercritical branching process with immigration. For each $n\geq 1$, the processes $(Z_0^n,Z_1^n)$ are defined on a common probability space $(\Omega_n, \mathcal{F}^n,\PP_n)$. For ease of notation we will write $\PP$ instead of $\PP_n$.

Define $z_1^n\left(t\right)=\EE Z_1^n\left(t\right)$, then (cf. \cite{bp}):
\begin{align*}
& z_1^n\left(t\right)=\frac{\mu}{\lambda_1-\lambda_0}n^{1-\alpha}e^{\lambda_1 t}\left(1-e^{\left(\lambda_0-\lambda_1\right)t}\right).
\end{align*}

We also present a few useful results for a birth-death process starting from a single cell. Let $Z=\{Z\left(t\right), t\ge 0\}$ denote a birth-death process where $Z\left(0\right)=1$ and each individual cell has birth rate $r_1$, death rate $d_1$, and net growth rate $\lambda_1=r_1-d_1$. The moment generating function of $Z\left(t\right)$ is given by 

\begin{align}\label{generating_function}
\phi_t\left(\theta\right) = \EE \exp\left(\theta Z\left(t\right)\right)=\left\{ \begin{array}{cc} 
\frac{d_1\left(e^{\theta}-1\right)-e^{-\lambda_1t}\left(r_1e^{\theta}-d_1\right)}{r_1\left(e^{\theta}-1\right)-e^{-\lambda_1t}\left(r_1e^{\theta}-d_1\right)}, & \hspace{5mm} \theta<\bar{\theta}_t \\
\infty & \hspace{5mm} \theta\ge\bar{\theta}_t \\
\end{array} \right.
\end{align}
where
\begin{equation}\label{bar_theta}
\bar{\theta}_t\doteq \log\left(\frac{r_1e^{\lambda_1t}-d_1}{r_1e^{\lambda_1 t}-r_1}\right)
\end{equation}
(see page 109 of \cite{bp}). Throughout this paper, we will repeatedly use $\phi_t$ to denote the moment generating function of $Z(t)$. If we let $\Psi_1=\frac{d_1e^{\lambda_1 t}-d_1}{r_1e^{\lambda_1 t}-d_1}$ and $\Psi_2=\frac{r_1e^{\lambda_1 t}-r_1}{r_1e^{\lambda_1 t}-d_1}$, then (see page 6 of \cite{Durrett})
$$
\PP\left(Z\left(t\right)=0\right) = \Psi_1, \quad \PP\left(Z\left(t\right)=n\right) = \left(1-\Psi_1\right)\left(1-\Psi_2\right)\Psi_2^{n-1} \quad \text{for } n\ge 1.
$$
We define the recurrence time as
\begin{equation*}
\gamma_n\left(a\right)=\inf\{t\ge 0: Z_1^n\left(t\right)>an \}
\end{equation*}
for $a>0$. The recurrence time represents the first time that the mutant cell population exceeds a proportion $a$ of the initial population size of drug-sensitive cells. We will often be interested in $\gamma_n(1)$, and use the notation $\gamma_n\equiv \gamma_n(1)$.
We denote by $\zeta_n\left(a\right)$ the unique value of $t$ such that $z_1^n\left(t\right)=an$.
It has been established in \cite{BP2020} that 
$$
\zeta_n\left(a\right)-\frac{1}{\lambda_1}\log\left(\frac{an^{\alpha}\left(\lambda_1-\lambda_0\right)}{\mu}\right)\rightarrow 0
$$ 
as $n\rightarrow \infty$, and $\gamma_n\left(a\right)-\zeta_n\left(a\right)\rightarrow 0$ in probability. Note that we will use the notation $\zeta_n\equiv \zeta_n(1)$. Consider the event of early recurrence such that recurrence happens $y$ units of time earlier than the deterministic limit, i.e., $\{\gamma_n\left(a\right)\le \zeta_n\left(a\right)-y\}$. We have the following large deviations result from \cite{BP2020}.
\begin{theorem}\label{recurrence_LD}
	Assume that $\alpha\in \left(0,1\right)$, then for $y>0$,
	\begin{equation*}
	\lim\limits_{n\rightarrow \infty}\frac{1}{n^{1-\alpha}}\log \PP \left(\gamma_n\left(a\right)\le \zeta_n\left(a\right)-y\right)= -\sup\limits_{\theta\in\left(0,1\right)}\left[\frac{ \mu  \theta e^{y \lambda_1}}{\lambda_1-\lambda_0}-\mu\theta\int_{0}^{\infty}\frac{e^{\lambda_0 s}}{e^{\lambda_1s}-\theta}ds\right].
	\end{equation*}
\end{theorem}
Note that the supremum on the right-hand side of the equation in Theorem \ref{recurrence_LD} is positive, and the solution to this optimization problem will be used frequently throughout the rest of the paper. 
\begin{definition}\label{def_theta_2}
We define $\theta^*_y$ to be the unique solution to the following equation 
\begin{align}\label{theta_2}
\frac{ e^{\lambda_1 y}}{\lambda_1-\lambda_0}= \int_{0}^{\infty}\frac{e^{\lambda_1 s}}{\left(e^{\lambda_1 s}- \theta\right)^2}e^{\lambda_0 s}ds.
\end{align}
\end{definition}
Note that $\theta^*_y$ satisfies the first order optimality condition of the optimization problem that appears in Theorem \ref{recurrence_LD}. Moreover, $\theta^*_y$ is a positive number because when $\theta=0$, the left-hand side of \eqref{theta_2} is larger than the right-hand side of \eqref{theta_2}, and the right-hand side of \eqref{theta_2} is an increasing function of $\theta$.

Throughout this work we will use the following notation for the asymptotic behavior of positive functions:
\begin{align*}
f\left(t\right) \sim g\left(t\right) & \quad\hbox{if $f\left(t\right)/g\left(t\right) \to 1$ as $t \to \infty$}, \\
f\left(t\right) = o\left(g\left(t\right)\right) & \quad\hbox{if $f\left(t\right)/g\left(t\right) \to 0$ as $t \to \infty$}, \\
f\left(t\right) = O\left(g\left(t\right)\right) & \quad \hbox{if $f\left(t\right) \leq C g\left(t\right)$ for all $t$}, \\
f\left(t\right) = \Theta\left(g\left(t\right)\right) & \quad \hbox{if $c g\left(t\right) \le f\left(t\right) \leq C g\left(t\right)$ for all $t$},
\end{align*}
where $C$ and $c$ are positive constants.

\section{Results on diversity indices without conditioning}\label{result_diversity_indices}

In this section, we examine the number of mutant clones and the Simpson's Index of mutant clones at the deterministic time $\zeta_n\left(1\right)$ which is a very good approximation of the cancer recurrence time $\gamma_n\left(1\right)$. We denote by $I_n\left(t\right)$ the number of mutant clones generated in the time period $\left(0,t\right)$. We first obtain the limit of the expected scaled number of mutant clones.
\begin{lemma}\label{expectation_without_conditioning}
\begin{align*}
\lim\limits_{n\rightarrow \infty}\frac{1}{n^{1-\alpha}}\EE\left[I_n\left(\zeta_n\right)\right]=-\frac{\mu}{\lambda_0}.
\end{align*}
\end{lemma}
This result tells us that the number of mutant clones at the deterministic limit of recurrence time is of order $\Theta\left(n^{1-\alpha}\right)$. The limit increases in $\mu$ and decreases in $\left|\lambda_0\right|$, which is expected as a higher mutation rate or a lower decaying rate of sensitive cells leads to more mutant clones.

We then investigate the Simpson's Index of mutant clones. Simpson's Index represents the probability that two randomly chosen cells from the mutant cell population come from the same clone. The Simpson's Index is close to $1$ if a few mutant clones dominate the mutant population. If a large number of mutant clones are similar in size, the Simpson's Index is near zero. 

Let $X_{i,n}$ denote the number of mutants at time $\zeta_n$ which belong to the $i$-th clone. Note that the mutant clones are ordered at random, not in chronological order by when the mutation occurred. Then the Simpson's index is computed by
\begin{align}
R_n\left(\zeta_n\right)=\sum_{i=1}^{I_n\left(\zeta_n\right)}\left(\frac{X_{i,n}}{Z_1^n\left(\zeta_n\right)}\right)^2,
\end{align}
and we define $R_n\left(\zeta_n\right)=0$ when $Z_1^n\left(\zeta_n\right)=0$. We have the following proposition for the large $n$ behavior of the expected Simpson's Index.
\begin{proposition}\label{simpson_without_conditioning}
\begin{align*}
\lim_{n\to\infty}n^{1-\alpha}\EE\left[R_n\left(\zeta_n\right)\right]= \frac{2\left(\lambda_1-\lambda_0\right)^2}{\mu \left(2\lambda_1-\lambda_0\right)}.
\end{align*}
\end{proposition}

This result tells us that the Simpson's Index of mutant clones at the deterministic limit of recurrence time is of order $\Theta\left(n^{\alpha-1}\right)$. We notice that the limit decreases in $\mu$, as a higher mutation rate results in a greater number of mutant clones, which leads to a lower Simpson's Index. We also notice that the limit increases in $\lambda_1$, which is owing to the fact that a higher growth rate of mutants results in larger clone sizes and a lesser number of clones at cancer recurrence, resulting in a higher Simpson's Index.

\section{Estimators for model parameters}\label{Sec:estimator}

In practice, the mutation rate $\mu n^{-\alpha}$, the net growth rate for drug-sensitive cells $\lambda_0$ and drug-resistant cells $\lambda_1$ under a therapy are important parameters in deciding patient treatment plans. Our theoretical results can help generate estimates for these three parameters. Note that in this section for ease of notation, we will use the notation $I_n \equiv I_n(\zeta_n)$ and $R_n \equiv R_n(\zeta_n)$.

From \cite{BP2020}, we know that $\gamma_n-\zeta_n$ converges in probability to zero, which implies that for sufficiently large $n$, with a very high probability,
\begin{align}
\gamma_n\approx \frac{1}{\lambda_1}\log\left(\frac{\lambda_1-\lambda_0}{\mu n^{-\alpha}}\right). \label{test_equation_1}
\end{align}
From Proposition \ref{simpson_without_conditioning}, we know that for sufficiently large $n$,
\begin{align}
\EE\left[R_n\right] \approx \frac{2\left(\lambda_1-\lambda_0\right)^2}{\mu n^{1-\alpha} \left(2\lambda_1-\lambda_0\right)}. \label{test_equation_2}
\end{align}
From Lemma \ref{expectation_without_conditioning}, we know that for sufficiently large $n$,
\begin{align}
\EE\left[I_n\right]\approx -\frac{\mu n^{1-\alpha}}{\lambda_0}. \label{test_equation_3}
\end{align}

First assume that for a given parameter set $(\mu n^{-\alpha}, \lambda_0,\lambda_1)$, we have $M$ independent observations, 
$$
\left\{\left(I_n^m,R_n^m,\gamma_n^m\right); m\in\{1,\ldots,M\}\right\}.
$$
We then define the sample averages
\begin{align*}
    \hat{I}_n(M)=\frac{1}{M}\sum_{m=1}^MI_n^m,\quad  \hat{R}_n(M)=\frac{1}{M}\sum_{m=1}^MR_n^m,\quad \hat\gamma_n(M)=\frac{1}{M}\sum_{m=1}^M\gamma_n^m.
\end{align*}
We can now use equations \eqref{test_equation_1},\eqref{test_equation_2},\eqref{test_equation_3} to derive the estimators \begin{align}
& \tilde{\lambda}_1(M)= \frac{1}{\hat\gamma_n(M)}\log\left(\frac{n}{\hat{I}_n(M)-\sqrt{\left(\hat{I}_n(M)\right)^2-\frac{2\hat{I}_n(M)}{\hat{R}_n(M)}}}\right),\label{estimate_1}\\
& \tilde{\lambda}_0(M)=\frac{\tilde{\lambda}_1(M)}{1-\frac{1}{n}\hat{I}_n(M)e^{\hat\gamma_n(M) \tilde{\lambda}_1(M)}}, \   \   \  \text{and } \label{estimate_2}\\
& \tilde{\mu}(M)=\frac{\tilde{\lambda}_1(M)}{e^{\hat\gamma_n(M)
\tilde{\lambda}_1(M)}-\frac{n}{\hat{I}_n(M)}}, \label{estimate_3}
\end{align}
via a method of moments approach. We conduct a simulation to evaluate our estimators. We obtain $100$ estimates with $M=100$. We then resample these estimates $100$ times to obtain a $95\%$ bootstrap confidence interval (Table \ref{tab:table1}). In the table, we observe that our estimators have a very small bias and variance.  Our estimators are applicable in practice because we only require data from patients at the time of detection and recurrence. We are not required to collect data from patients during treatment, which is often impractical, particularly for solid tumors. However, it should be noted that our estimators do require knowledge of the initial tumor burden $n.$ Another drawback of our estimators is that we require multiple independent observations, i.e, $M>1$. In addition, our estimators assume that mutants have a death rate of $0$. We believe that, it is possible to derive consistent estimators (in the large $n$ limit) based on a single sample, i.e., $M=1$, without assuming mutants have zero death rate. In a forthcoming work we investigate these more general estimators.

\begin{table}[h!]
\caption{Estimation of $(\mu n^{-\alpha}, \lambda_0, \lambda_1)$. Model Parameters: $r_0=1$, $d_0=1.2$, $\lambda_0=-0.2$, $\lambda_1=0.2$, $\mu=0.5$, $\alpha=0.6$, $n=100000$. }
\label{tab:table1}
\begin{center}
	\begin{tabular}{| l | l | l | l | }
		\hline
		  & $\mu n^{-\alpha}$ & $\lambda_0$ & $\lambda_1$ \\ \hline
		True value  & $5\times10^{-4}$ & $-0.2$ &  $0.2$ \\ \hline
		Estimate value & $5.0686\times10^{-4}$ & $ -0.2030$ &  $0.1996$ \\ \hline
		Bootstrap 95\% C.I. & $[5.035,  5.101]\times 10^{-4}$ & $[-0.2042, -0.2018]$ &  $[0.1995,  0.1997]$ \\ \hline
	\end{tabular}
\end{center}
\end{table}

\section{Results on diversity indices conditioned on early recurrence}\label{sec:results_with_conditioning}

In this section, we examine the number of mutant clones and the Simpson's Index of mutant clones at cancer recurrence conditioned on the event of early recurrence ($\{\gamma_n \le \zeta_n-y\}$). Note that throughout this section $y$ is a positive number independent of $n$. We compare results conditioned on early recurrence to those without conditioning. Our goal is to determine whether early recurrence is primarily caused by a single mutation, leading to an unusually large family of drug-resistant cells (resulting in low clonal diversity), or if it is instead attributed to an abnormally high number of mutations causing drug resistance (resulting in high clonal diversity).

In order to obtain results conditioned on early recurrence, we need to add a strong assumption that sensitive cells have deterministic exponential decay (i.e., $Z_0^n\left(t\right)=z_0^n\left(t\right)$, where, abusing the notation, $z_0^n\left(t\right)=\EE \left[Z_0^n\left(t\right)\right]$ in the original model). Note that this assumption does not affect results obtained in previous sections. Therefore, we use the same notation introduced in previous sections.

\subsection{Number of mutant clones conditioned on early recurrence}\label{sec:num_clones_with_conditioning}

Because of the assumption that sensitive cells have deterministic exponential decay, we can study the distribution of the number of clones at cancer recurrence. For non-negative integer valued random variables $X$ and $Y$, their total variation distance is given by
$$
TV\left(X,Y\right)=\sum_{k=0}^\infty \left|P(X=k)-P(Y=k)\right|.
$$ 

We first show under assumption on $\alpha$ that $I_n\left(\gamma_n\right)$ is close to $I_n\left(\zeta_n\right)$ in total variation distance. 
\begin{proposition}\label{dist_1}
If $\alpha\in \left(\frac{\lambda_1}{\lambda_1-\lambda_0},1\right)$, then 
$$
\lim_{n\to\infty}TV\left(I_n\left(\gamma_n\right),I_n\left(\zeta_n\right)\right)=0.
$$
\end{proposition}
Because sensitive cells have deterministic exponential decay, we know that $\left(I_n\left(t\right)\right)_{t\ge 0}$ is a non-homogeneous Poisson process, which gives us that
\begin{align*}
\PP\left(I_{n}\left(\zeta_n\right)=k\right)=\frac{\lambda_n^k e^{-\lambda_n}}{k!},
\end{align*}
where $\lambda_n=-\frac{\mu n^{1-\alpha}}{\lambda_0}\left(1-e^{\lambda_0\zeta_n}\right)$. Proposition \ref{dist_1} tells us that the distribution of the number of mutant clones at cancer recurrence is close to a Poisson distribution with mean $\lambda_n$.

We then show that the number of clones at the cancer recurrence time conditioned on the event of early recurrence ($\{\gamma_n \le \zeta_n-y\}$) stochastically dominates that without conditioning asymptotically under assumption on $\alpha$.

\begin{figure}[!ht]
\begin{center}
	\includegraphics[width=0.8\textwidth]{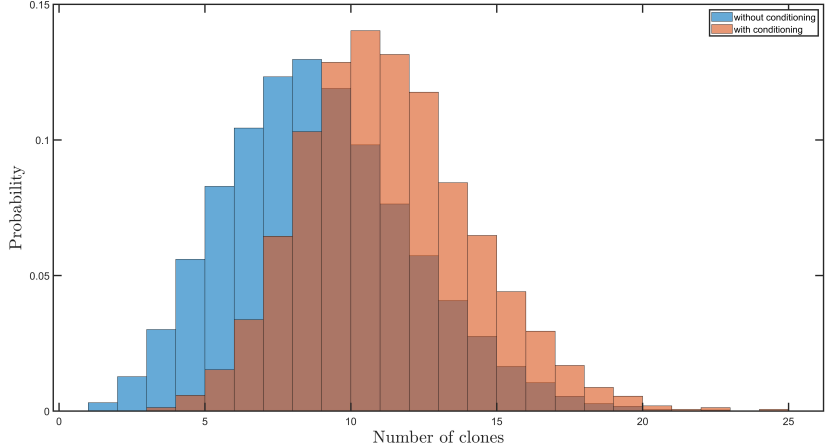}
\end{center}
\caption{Scaled histogram of number of clones present at recurrence for simulations conditioned on early recurrence (light brown) and unconditioned simulations (blue). Histogram is based on $10^4$ simulations and all simulations used the model parameter set $r_0=1$, $d_0=1.2$, $r_1=1$, $d_1=0.8$, $\mu=0.5$, $\alpha=0.6$, $y=1$, $n=1,000$ and $a=1$.}
\label{fig:num_of_clones}
\end{figure}

\begin{proposition}\label{proposition_diversity}
If $\alpha\in \left(\frac{\lambda_1}{\lambda_1-\lambda_0}\vee \frac{1}{2},1\right)$, then
\begin{align*}
\liminf\limits_{n\rightarrow \infty}\inf\limits_{x>0}\left(\PP\left(I_n\left(\gamma_n\right)\ge x|\gamma_n<\zeta_n-y\right)-\PP\left(I_n\left(\gamma_n\right)\ge x\right)\right)\ge 0.
\end{align*}
\end{proposition}
This finding sheds light on the effect of early recurrence on the distribution of mutant clones at cancer recurrence. If a patient experiences an early cancer recurrence, the number of mutant clones in the recurring tumor is very likely to be higher than expected. As a result, the recurrent tumor is more likely to be resistant to second-line medications, which should be taken into account when determining future treatment options. 
Ideally, we would like to generalize the evolution of each mutant clone to a birth-death process (with birth rate $r_1$, death rate $d_1$, net growth rate $\lambda_1 = r_1-d_1$) instead of a Yule process. Unfortunately, we were unable to do so. Instead we conduct a simulation to visually display the distribution of the number of clones present at cancer recurrence for the more general model (see Figure \ref{fig:num_of_clones}).

From Proposition \ref{proposition_diversity}, we know that the number of mutant clones at cancer recurrence conditioned on early recurrence stochastically dominates that without conditioning in the large population limit. Because a stochastically larger random variable has a larger expectation, the number of mutant clones at cancer recurrence conditioned on early recurrence should have a higher expectation than that without conditioning. In Lemma \ref{expectation_with_conditioning}, we obtain the expectation of the scaled number of mutant clones at $\zeta_n-y$ conditioned on early recurrence. For simplicity, we denote by $A_{n,y}$ the event of early recurrence $\{ \gamma_n < \zeta_n-y\}$, and the corresponding conditional probability measure by $\PP_{A_{n,y}}\left(\cdot\right)=\PP\left(\cdot \middle| A_{n,y}\right)$.
\begin{lemma}\label{expectation_with_conditioning}
\begin{align*}
\lim\limits_{n\rightarrow \infty}\frac{1}{n^{1-\alpha}}\EE_{A_{n,y}}\left[I_n\left(\zeta_n-y\right)\right]=\mu \int_{0}^{\infty}\frac{e^{\lambda_1 s}}{ e^{\lambda_1 s}-\theta^*_y}e^{\lambda_0 s}ds,
\end{align*}
where $\theta^*_y$ is defined in \eqref{theta_2}.
\end{lemma}
It is easy to observe that $\mu \int_{0}^{\infty}\frac{e^{\lambda_1 s}}{ e^{\lambda_1 s}-\theta^*_y}e^{\lambda_0s}ds$ increases in $y$, and hence for any $y>0$,
\begin{align*}
\mu \int_{0}^{\infty}\frac{e^{\lambda_1 s}}{ e^{\lambda_1 s}-\theta^*_y}e^{\lambda_0s}ds> -\frac{\mu}{\lambda_0}.
\end{align*}
This result is consistent with Proposition \ref{proposition_diversity}. Moreover, we can see that the recurrence time is an indicator of the clonal diversity at cancer recurrence. 

Lastly, we investigate the number of clones that are generated in a given time period over the course of the treatment. Let $I_{n, \left(t_1, t_2\right)}$ denote the number of clones generated in the time period $\left(t_1, t_2 \right)$. Let $B_{n,y,\left(t_1, t_2\right)}$ denote the number of mutants at time $\zeta_n-y$ which are descendants of those clones generated in the time period $\left(t_1,t_2\right)$. For simplicity, we let
\begin{align}\label{def:bar_I}
\bar{I}_{n, \left(t_1, t_2\right)}&=\EE\left[I_{n, \left(t_1, t_2\right)}\right] = \mu n^{1-\alpha}\int_{t_1}^{t_2}e^{\lambda_0 t}dt,
\end{align}
and 
\begin{align*}
\bar{B}_{n,y,\left(t_1, t_2\right)}=\EE\left[B_{n,y,\left(t_1, t_2\right)}\right].
\end{align*}
Recall that we denote by $A_{n,y}$ the event of early recurrence $\{ \gamma_n < \zeta_n-y\}$, and the corresponding conditional probability measure by $\PP_{A_{n,y}}\left(\cdot\right)=\PP\left(\cdot \middle| A_{n,y}\right)$.

We are interested in the number $I_{n, \left(t_1, t_2\right)}$ conditioned on early recurrence. We first show the following lemma.
\begin{lemma}\label{period_mutants_concentration}
For any $\epsilon>0$,
\begin{align*}
\limsup\limits_{n\rightarrow \infty}\frac{1}{n^{1-\alpha}}\log \PP_{A_{n,y}}\left(\frac{B_{n,y,\left(t_1, t_2\right)}}{\bar
{B}_{n,y,\left(t_1, t_2\right)}}\notin \left(1+\delta^*-\epsilon, 1+\delta^*+\epsilon\right)\right)<0,
\end{align*}
where $\delta^*$ is given by
\begin{align*}
\delta^*=\frac{\int_{t_1}^{t_2}\frac{e^{\lambda_1s}}{\left(e^{\lambda_1 s}-\theta^*_y\right)^2}e^{\lambda_0 s}ds}{\int_{t_1}^{t_2}e^{-\left(\lambda_1-\lambda_0\right)s}ds}-1>0,
\end{align*}
and $\theta^*_y$ is defined in \eqref{theta_2}.
\end{lemma}
This result tells us that the number of mutants at time $\zeta_n-y$ which are descendants of those clones generated in the time period $\left(t_1,t_2\right)$ is concentrated around $\left(1+\delta^*\right)\bar{B}_{n,y,\left(t_1, t_2\right)}$. 

We then analyze the number $I_{n, \left(t_1, t_2\right)}$ conditioned on the event
\begin{align*}
A_{n,y,\epsilon}^{t_1,t_2} = \{B_{n,y,\left(t_1, t_2\right)}\in \left(\left(1+\delta^*-\epsilon\right)\bar
{B}_{n,y,\left(t_1, t_2\right)}, \left(1+\delta^*+\epsilon\right)\bar
{B}_{n,y,\left(t_1, t_2\right)}\right)\}.
\end{align*}
We have the following lemma.
\begin{lemma}\label{period_clones_concentration}
Assume that $t_2-t_1< -\frac{1}{\lambda_1}\log\left(\frac{1}{2-\theta^*_y}\right)$. For any $\sigma>0$, there exists $\bar{\epsilon}>0$ such that when $0<\epsilon<\bar{\epsilon}$,
\begin{align*}
\limsup\limits_{n\rightarrow \infty}\frac{1}{n^{1-\alpha}}\log \PP_{A_{n,y,\epsilon}^{t_1,t_2}}\left(\frac{I_{n,\left(t_1, t_2\right)}}{\bar{I}_{n,\left(t_1, t_2\right)}}\notin \left(1+\kappa^*-\sigma, 1+\kappa^*+\sigma\right)\right)<0,
\end{align*}
where
\begin{align*}
\kappa^*=\frac{\int_{t_1}^{t_2}\frac{e^{\lambda_1 s}}{e^{\lambda_1 s}-\theta^*_y}e^{\lambda_0 s}ds}{\int_{t_1}^{t_2}e^{\lambda_0 s}ds}-1.
\end{align*}
\end{lemma}
This result tells us that conditioned on the event that the number of mutants is concentrated around $\left(1+\delta^*\right)\bar{B}_{n,y,\left(t_1, t_2\right)}$, the number of clones generated in the time period $\left(t_1,t_2\right)$ is concentrated around $\left(1+\kappa^*\right)\bar{I}_{n,\left(t_1, t_2\right)}$.

From Lemma \ref{period_mutants_concentration} and \ref{period_clones_concentration}, we can obtain the desired result.
\begin{proposition}\label{period_clones_result}
Assume that $t_2-t_1< -\frac{1}{\lambda_1}\log\left(\frac{1}{2-\theta^*_y}\right)$. For any $\sigma>0$,
\begin{align*}
\limsup\limits_{n\rightarrow \infty}\frac{1}{n^{1-\alpha}}\log \PP_{A_{n,y}}\left(\frac{I_{n,\left(t_1, t_2\right)}}{\bar
{I}_{n,\left(t_1, t_2\right)}}\notin \left(1+\kappa^*-\sigma, 1+\kappa^*+\sigma\right)\right)<0,
\end{align*}
where $\kappa^*$ is defined in Lemma \ref{period_clones_concentration}.
\end{proposition}

Note that $\delta^*$ and $\kappa^*$ depend on $t_1$ and $t_2$, and we omit the dependence in their notations for simplicity. If we fix $t_1$ and let $t_2-t_1$ goes to zero, then $1+\delta^*$ goes to 
$$
\frac{e^{2\lambda_1 t_1}}{\left(e^{\lambda_1 t_1}-\theta^*_y\right)^2},
$$ and $1+\kappa^*$ goes to $\sqrt{1+\delta^*}$. Hence, the number of clones generated in an infinitesimal time period $\left(t_1, t_1+dt\right)$ conditioned on early recurrence is approximately 
$\frac{e^{\lambda_1 t_1}}{e^{\lambda_1 t_1}-\theta^*_y}\bar{I}_{n, \left(t_1, t_1+dt\right)}$. 

Lemma \ref{period_mutants_concentration} tells us that conditioned on early recurrence, the number of mutants at time $\zeta_n-y$ which are descendants of those clones generated in a small time period $\left(t_1,t_2\right)$ is larger than that without conditioning by the  approximate factor
$$
\frac{e^{2\lambda_1 t_1}}{\left(e^{\lambda_1 t_1}-\theta^*_y\right)^2}.
$$ 
Moreover, we can specify the contribution to such an increase in the number of mutants from (1) the increase in the number of clones, and (2) the increase in the average clone size. Our results (Lemma \ref{period_mutants_concentration} and Proposition \ref{period_clones_result}) indicate that both factors contribute equally to the increase in the number of mutants such that the number of clones generated in the time period $\left(t_1,t_2\right)$ conditioned on early recurrence is larger than that without conditioning by the approximate factor
$$
\frac{e^{\lambda_1 t_1}}{e^{\lambda_1 t_1}-\theta^*_y},
$$ 
and the average clone size conditioned on early recurrence is larger than that without conditioning by approximately $\frac{e^{\lambda_1 t_1}}{e^{\lambda_1 t_1}-\theta^*_y}$ times as well.

\subsection{Simpson's Index conditioned on early recurrence}\label{sec:simpson_with_conditioning}

In this section, we investigate the Simpson's Index of mutant clones at time $\zeta_n-y$ conditioned on early recurrence. Recall that we denote by $A_{n,y}$ the event of early recurrence $\{\gamma_n < \zeta_n-y\}$, and the corresponding conditional probability measure by $\PP_{A_{n,y}}\left(\cdot\right)=\PP\left(\cdot \middle| A_{n,y}\right)$. Let $R_{n,y}$ be the Simpson's Index of mutant clones at time $\zeta_n-y$. We have the following theorem.

\begin{theorem}\label{Simpson_conditioning_I}
\begin{align}\label{eq:LimitSI_cond}
\lim\limits_{n\rightarrow \infty} n^{1-\alpha}\EE_{A_{n,y}}\left[R_{n,y}\right]& =\frac{2\left(\lambda_1-\lambda_0\right)^2}{\mu} e^{-2\lambda_1 y}\int_0^{\infty}\frac{e^{-\left(2\lambda_1-\lambda_0\right)s}}{\left(1-\theta^*_y e^{-\lambda_1 s}\right)^3}ds\doteq S_c(y),
\end{align}
where $\theta^*_y$ is defined in \eqref{theta_2}.
\end{theorem}

\begin{figure}[ht]\label{plot_simpson_conditioning}
\hspace*{0cm}  
\includegraphics[width=1.05\textwidth]{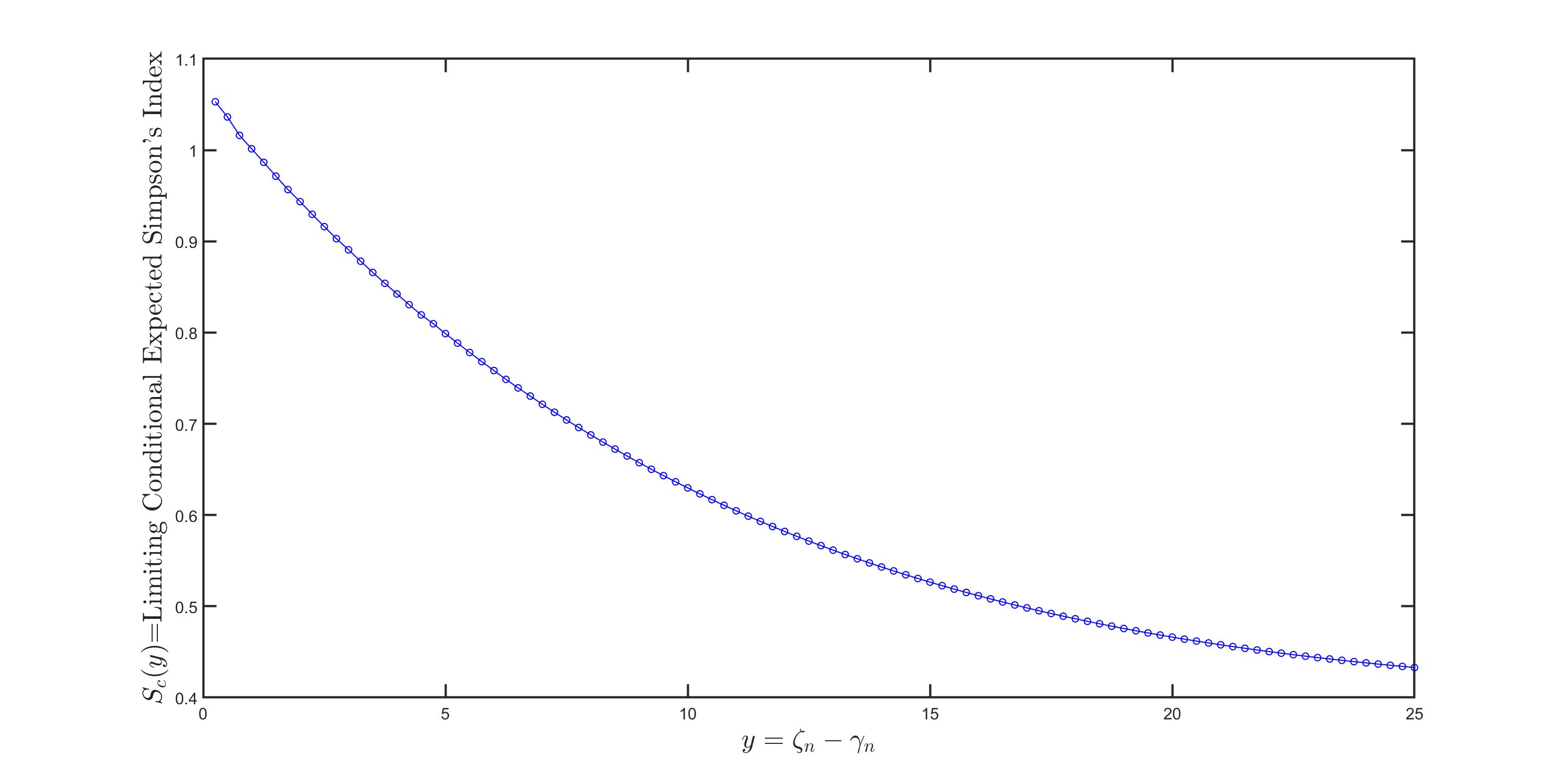}
\caption{Plot of $S_c(y)$ versus $y$, $S_c(y)$ is defined in \eqref{eq:LimitSI_cond}. Model Parameters: $r_0=1$, $d_0=1.2$, $\lambda_0=-0.2$, $\lambda_1=0.2$, $\mu=0.5$, $\alpha=0.6$.}
\end{figure}

This result tells us that the Simpson's Index of mutant clones at $\zeta_n-y$ conditioned on early recurrence is of order $\Theta\left(n^{\alpha-1}\right)$ as well. The limit is associated with the value of $y$. We conjecture that the limit $S_c(y)$ decreases in $y$, and it is always smaller than that obtained in Proposition \ref{simpson_without_conditioning}. However, we were not able to prove this result. Our conjecture is supported by numerical results, an example of which is given in Figure \ref{plot_simpson_conditioning}. 

Our analysis of both clonal diversity indices reveals that early recurrent tumors are more likely to exhibit higher clonal diversity. This suggests that early recurrence is primarily driven by a larger number of mutations, rather than the successful expansion of a single mutation. Furthermore, our findings highlight that the time of recurrence serves as an important indicator of clonal diversity at cancer recurrence.

\section{Summary}\label{Sec_summary}

In this work, we have examined the clonal diversity of mutant clones at cancer recurrence. We obtained the expectation of the number of mutant clones and the Simpson's Index at cancer recurrence with and without conditioning on early recurrence. We showed that the number of mutant clones at cancer recurrence conditioned on early recurrence stochastically dominates that without conditioning in the large population limit. In addition we have derived an expression for the large population limit of Simpson's Index conditioned on early recurrence. Our findings suggest that an earlier recurrent tumor is more likely to have a higher clonal diversity. Furthermore, our results suggest that early recurrence is most likely to be generated by a larger number of mutations, as opposed to the successful growth of a single mutation. In this work we use the method of moments to develop estimators for the mutation rate, the net growth rate for drug-sensitive cells, and drug-resistant cells. We show in numerical results that these estimators are able to accurately estimate model parameters. 


\newpage

\clearpage

\section{Proofs of results in Section \ref{result_diversity_indices}}\label{Sec:proof_result_diversity_indices}

\subsection{Proof of Lemma \ref{expectation_without_conditioning}}
\begin{proof}
This result can be derived from direct calculation and thus the proof is omitted.
\end{proof}
\\

\subsection{Proof of Proposition \ref{simpson_without_conditioning}}\label{proposition_simpson_without_conditioning}

\begin{proof}\\
Note that in this proof for ease of notation, we will use the notation $R_n \equiv R_n(\zeta_n)$. Recall that $I_n\left(\zeta_n\right)$ represents the number of clones generated in the time period $\left(0,\zeta_n\right)$. We define
\begin{align}
\tilde{R}_n=\sum_{i=1}^{I_n\left(\zeta_n\right)}\left(\frac{X_{i,n}}{n}\right)^2.
\end{align}
We can compute that 
\begin{align}
n^{1-\alpha}\EE\left[\tilde{R}_n\right]&=n^{1-\alpha}\EE\left[\sum_{i=1}^{I_n\left(\zeta_n\right)}\left(\frac{X_{i,n}}{n}\right)^2\right] \nonumber\\
&=n^{-1-\alpha}\EE\left[\EE\left[\EE\left[\sum_{i=1}^{I_n\left(\zeta_n\right)} X_{i,n}^2\middle| I_n\left(\zeta_n\right)\right]\middle| Z_0^n\right]\right] \nonumber\\
&\overset{\text{(a)}}{=}n^{-1-\alpha}\EE\left[\int_0^{\zeta_n}\mu n^{-\alpha}Z_0^n\left(s\right)\EE\left[\left(Z\left(\zeta_n-s\right)\right)^2\right]ds\right] \nonumber\\
&=n^{-1-\alpha}\int_0^{\zeta_n}\mu n^{1-\alpha}e^{\lambda_0 s}\EE\left[\left(Z\left(\zeta_n-s\right)\right)^2\right]ds \nonumber\\
&\overset{\text{(b)}}{=}n^{-1-\alpha}\int_0^{\zeta_n}\mu n^{1-\alpha}e^{\lambda_0 s}\left(2e^{2\lambda_1\left(\zeta_n-s\right)}-e^{\lambda_1 \left(\zeta_n-s\right)}\right)ds \nonumber\\
&\overset{\text{(c)}}{\rightarrow} \frac{2\left(\lambda_1-\lambda_0\right)^2}{\mu \left(2\lambda_1-\lambda_0\right)} \text{ as } n\rightarrow \infty, \label{eqn:convergence_R_tilde}
 \end{align}
where we use the uniformity of arrival times for a Poisson process in step (a), $\EE\left[\left(Z\left(t\right)\right)^2\right]=2e^{2\lambda_1 t}-e^{\lambda_1 t}$ in step (b), and $e^{\lambda_1 \zeta_n}\sim \frac{\lambda_1-\lambda_0}{\mu}n^{\alpha}$ in step (c). Also note that the process $Z(\cdot)$ in step (a) is defined in page 4 before \eqref{generating_function} with $d_1=0$. We then compute the difference between $\EE\left[\tilde{R}_n\right]$ and $\EE\left[R_n\right]$. Recall that we define $R_n=0$ when $Z_1^n\left(\zeta_n\right)=0$. It allows us to work on the event that $Z_1^n\left(\zeta_n\right)>0$. Hence, for the rest of this section, we condition on the event $\rho_n = \{Z_1^n\left(\zeta_n\right)>0\}$, and denote by $\EE_{\rho_n}$ the conditional expectation. We obtain that
\begin{align*}
n^{1-\alpha}\EE_{\rho_n}\left[\left\lvert \tilde{R}_n-R_n\right\rvert\right]&=n^{1-\alpha}\EE_{\rho_n}\left[\left\lvert\sum_{i=1}^{I_n\left(\zeta_n\right)}\left(\frac{X_{i,n}}{Z_1^n\left(\zeta_n\right)}\right)^2-\sum_{i=1}^{I_n\left(\zeta_n\right)}\left(\frac{X_{i,n}}{n}\right)^2\right\rvert\right]\\
&=n^{1-\alpha}\EE_{\rho_n}\left[\left\lvert\frac{\sum_{i=1}^{I_n\left(\zeta_n\right)}\left(X_{i,n}\right)^2}{n^2}\left(\frac{n^2}{Z_1^n\left(\zeta_n\right)^2}-1\right)\right\rvert\right]\\
&\le n^{1-\alpha}\EE_{\rho_n}\left[\left(\sum\limits_{i=1}^{I_n\left(\zeta_n\right)}\left(X_{i,n}\right)^2\right)^2/n^4\right]^{1/2}\EE_{\rho_n}\left[\left(\frac{n^2}{Z_1^n\left(\zeta_n\right)^2}-1\right)^2\right]^{1/2}.
\end{align*}
We can show that (for some positive constants $c_1$ and $c_2$)
\begin{align*}
& \quad \EE\left[\left(\sum\limits_{i=1}^{I_n\left(\zeta_n\right)}\left(X_{i,n}\right)^2\right)^2/n^4\right]\\
& = \frac{1}{n^4}\int_0^{\zeta_n}\mu n^{1-\alpha}e^{\lambda_0 s}\EE\left[\left(Z\left(\zeta_n-s\right)\right)^4\right]ds\\
& \quad \quad +\frac{1}{n^4}\EE\left[\frac{\EE\left[I_n\left(\zeta_n\right)^2-I_n\left(\zeta_n\right)\middle| Z_0^n\right]}{\EE\left[I_n\left(\zeta_n\right)\middle| Z_0^n\right]^2}\left(\int_0^{\zeta_n}\mu n^{1-\alpha}Z_0^n\left(s\right)\EE\left[\left(Z\left(\zeta_n-s\right)\right)^2\right]ds\right)^2\right]\\
& = \frac{1}{n^4}\int_0^{\zeta_n}\mu n^{1-\alpha}e^{\lambda_0 s}\EE\left[\left(Z\left(\zeta_n-s\right)\right)^4\right]ds\\
& \quad \quad +\frac{1}{n^4}\EE\left[\left(\int_0^{\zeta_n}\mu n^{-\alpha}Z_0^n\left(s\right)\EE\left[\left(Z\left(\zeta_n-s\right)\right)^2\right]ds\right)^2\right]\\
& \sim c_1 n^{-2+2\alpha},
\end{align*}
where we use the fact that $\frac{\EE\left[I_n\left(\zeta_n\right)^2-I_n\left(\zeta_n\right)\middle| Z_0^n\right]}{\EE\left[I_n\left(\zeta_n\right)\middle| Z_0^n\right]^2}=1$, $\EE\left[\left(Z\left(t\right)\right)^4\right]\sim 24e^{4\lambda_1 t}$, and 
\begin{align*}
& \quad \EE\left[\left(\int_0^{\zeta_n}Z_0^n\left(s\right)\EE\left[\left(Z\left(\zeta_n-s\right)\right)^2\right]ds\right)^2\right]\\
& =  \EE\left[\left(\int_0^{\zeta_n}Z_0^n\left(s\right)\left(2e^{2\lambda_1 \left(\zeta_n-s\right)}-e^{\lambda_1 \left(\zeta_n-s\right)}\right)ds\right)^2\right]\\
& = \int_0^{\zeta_n}\int_0^{\zeta_n}\EE\left[Z_0^n\left(s\right)\left(2e^{2\lambda_1 \left(\zeta_n-s\right)}-e^{\lambda_1 \left(\zeta_n-s\right)}\right)Z_0^n\left(t\right)\left(2e^{2\lambda_1 \left(\zeta_n-t\right)}-e^{\lambda_1 \left(\zeta_n-t\right)}\right)\right]dtds\\
& =\int_0^{\zeta_n}\int_s^{\zeta_n}\EE\left[Z_0^n\left(s\right)Z_0^n\left(s\right)e^{\lambda_0(t-s)}\left(2e^{2\lambda_1 \left(\zeta_n-s\right)}-e^{\lambda_1 \left(\zeta_n-s\right)}\right)\left(2e^{2\lambda_1 \left(\zeta_n-t\right)}-e^{\lambda_1 \left(\zeta_n-t\right)}\right)\right]dtds\\
& \quad + \int_0^{\zeta_n}\int_0^{s}\EE\left[Z_0^n\left(t\right)Z_0^n\left(t\right)e^{\lambda_0(s-t)}\left(2e^{2\lambda_1 \left(\zeta_n-s\right)}-e^{\lambda_1 \left(\zeta_n-s\right)}\right)\left(2e^{2\lambda_1 \left(\zeta_n-t\right)}-e^{\lambda_1 \left(\zeta_n-t\right)}\right)\right]dtds\\
& \sim c_2 n^{2+4\alpha}.
\end{align*}
We can then obtain that
\begin{align*}
n^{1-\alpha}\EE_{\rho_n}\left[\left(\sum\limits_{i=1}^{I_n\left(\zeta_n\right)}\left(X_{i,n}\right)^2\right)^2/n^4\right]^{1/2}\rightarrow \sqrt{c_1}.
\end{align*}
It remains to analyze
\begin{align*}
\EE_{\rho_n}\left[\left(\frac{n^2}{Z_1^n\left(\zeta_n\right)^2}-1\right)^2\right]=\EE_{\rho_n}\left[\frac{n^4}{Z_1^n\left(\zeta_n\right)^4}-\frac{2n^2}{Z_1^n\left(\zeta_n\right)^2}+1\right].
\end{align*}
We need to analyze the probability that $Z_1^n\left(\zeta_n\right)$ deviates from its mean $n$. Recall that
\begin{align*}
\phi_{t}\left(\theta\right)=\frac{e^{\theta}e^{-\lambda_1 t}}{1-\left(1-e^{-\lambda_1 t}\right)e^{\theta}}.
\end{align*}
For $\epsilon\in \left(0,1\right)$, and $\theta>0$, we have
\begin{align}
\PP\left(Z_1^n\left(\zeta_n\right)<\left(1-\epsilon\right)n\right)&=\PP\left(e^{-\theta Z_1^n\left(\zeta_n\right)}>e^{-\theta\left(1-\epsilon\right)n}\right) \nonumber\\
& \le \min\limits_{\theta>0}e^{\theta\left(1-\epsilon\right)n}\EE\left[e^{-\theta Z_1^n\left(\zeta_n\right)}\right] \nonumber\\
&=\min\limits_{\theta>0}e^{\theta\left(1-\epsilon\right)n} \EE\left[\exp\left(\frac{\mu}{n^{\alpha}}\int_{0}^{\zeta_n}Z_0^n\left(s\right)\left(\phi_{\zeta_n-s}\left(-\theta\right)-1\right)ds\right)\right]. \label{exp_term}
\end{align}
The expectation term in (\ref{exp_term}) can be decomposed into a mean behavior term and a fluctuation term:
\begin{align*}
& \quad \EE\left[\exp\left(\frac{\mu}{n^{\alpha}}\int_{0}^{\zeta_n}Z_0^n\left(s\right)\left(\phi_{\zeta_n-s}\left(-\theta\right)-1\right)ds\right)\right]\\
& =\exp\left(\frac{\mu}{n^{\alpha-1}}\int_{0}^{\zeta_n}e^{\lambda_0 s}\left(\phi_{\zeta_n-s}\left(-\theta\right)-1\right)ds\right) \\
& \quad \times \EE \left[\exp \left(\frac{\mu}{n^{\alpha}}\int_{0}^{\zeta_n}\left(\phi_{\zeta_n-s}\left(-\theta\right)-1\right)\left(Z_0^n\left(s\right)-ne^{\lambda_0 s}\right)ds\right)\right].
\end{align*}
By Proposition 1 of \cite{BP2020}, we can safely discard the fluctuation term, and focus on the remaining terms in the exponential expression of (\ref{exp_term}):
\begin{align*}
& \quad \theta\left(1-\epsilon\right)n+\frac{\mu}{n^{\alpha}}\int_{0}^{\zeta_n}ne^{\lambda_0 s}\left(\phi_{\zeta_n-s}\left(-\theta\right)-1\right)ds\\
&=\theta\left(1-\epsilon\right)n+\frac{\mu}{n^{\alpha}}\int_{0}^{\zeta_n}ne^{\lambda_0 s}\left(\frac{e^{-\theta}e^{-\lambda_1 \left(\zeta_n-s\right)}}{1-\left(1-e^{-\lambda_1 \left(\zeta_n-s\right)}\right)e^{-\theta}}-1\right)ds.
\end{align*}
Let $\theta=\delta e^{-\lambda_1 \zeta_n}$, we can obtain
\begin{align*}
& \quad n^{\alpha-1}\left(\theta\left(1-\epsilon\right)n+\frac{\mu}{n^{\alpha}}\int_{0}^{\zeta_n}ne^{\lambda_0 s}\left(\frac{e^{-\theta}e^{-\lambda_1 \left(\zeta_n-s\right)}}{1-\left(1-e^{-\lambda_1 \left(\zeta_n-s\right)}\right)e^{-\theta}}-1\right)ds\right)\\
&=n^{\alpha-1}\delta e^{-\lambda_1 \zeta_n}\left(1-\epsilon\right)n+n^{\alpha-1}\frac{\mu}{n^{\alpha}}\int_{0}^{\zeta_n}ne^{\lambda_0 s}\left(\frac{e^{-\delta e^{-\lambda_1 \zeta_n}}e^{-\lambda_1 \left(\zeta_n-s\right)}}{1-\left(1-e^{-\lambda_1 \left(\zeta_n-s\right)}\right)e^{-\delta e^{-\lambda_1 \zeta_n}}}-1\right)ds\\
& \rightarrow \delta \frac{\mu}{\lambda_1-\lambda_0}\left(1-\epsilon\right)-\mu\int_{0}^{\infty}e^{\lambda_0 s}\frac{\delta}{\delta+e^{\lambda_1 s}}ds\\
&=\delta \mu \left(\frac{1-\epsilon}{\lambda_1-\lambda_0}-\int_{0}^{\infty}\frac{e^{\lambda_0 s}}{\delta+e^{\lambda_1 s}}ds\right)\\
&=\delta \mu \left(\int_0^{\infty}\frac{\left(1-\epsilon\right)e^{\lambda_0 s}}{e^{\lambda_1 s}}ds-\int_{0}^{\infty}\frac{e^{\lambda_0 s}}{\delta+e^{\lambda_1 s}}ds\right).
\end{align*}
By comparing the integrand, we have that
\begin{align*}
\frac{1-\epsilon}{e^{\lambda_1 s}}-\frac{1}{\delta+e^{\lambda_1 s}}<0 & \Leftrightarrow \left(1-\epsilon\right)\delta-\epsilon e^{\lambda_1 s}<0\\
& \Leftrightarrow \delta < \frac{\epsilon e^{\lambda_1 s}}{1-\epsilon}\\
& \Leftarrow \delta < \frac{\epsilon}{1-\epsilon}.
\end{align*}
Hence, we conclude that for some $c>0$
\begin{align*}
\limsup\limits_{n\rightarrow \infty}\frac{1}{n^{1-\alpha}}\log \PP\left(Z_1^n\left(\zeta_n\right)<\left(1-\epsilon\right)n\right)\le -c.
\end{align*}
Since $\PP\left(\rho_n\right)\rightarrow 1$, the above result also holds when conditioned on the event $\rho_n$. Hence, we can obtain that for any $\epsilon\in \left(0, 1\right)$,
\begin{align*}
\quad \EE_{\rho_n}\left[\frac{n^4}{Z_1^n\left(\zeta_n\right)^4}\right] & = \EE_{\rho_n}\left[\frac{n^4}{Z_1^n\left(\zeta_n\right)^4}\vert Z_1^n\left(\zeta_n\right)<\left(1-\epsilon\right)n\right]\PP_{\rho_n}\left(Z_1^n\left(\zeta_n\right)<\left(1-\epsilon\right)n\right)\\
& \quad \quad + \EE_{\rho_n}\left[\frac{n^4}{Z_1^n\left(\zeta_n\right)^4}\vert Z_1^n\left(\zeta_n\right)\ge \left(1-\epsilon\right)n\right]\PP_{\rho_n}\left(Z_1^n\left(\zeta_n\right)\ge\left(1-\epsilon\right)n\right)\\
& \le n^4\PP_{\rho_n}\left(Z_1^n\left(\zeta_n\right)<\left(1-\epsilon\right)n\right)+\frac{1}{\left(1-\epsilon\right)^4}\\
& \rightarrow \frac{1}{\left(1-\epsilon\right)^4}.
\end{align*}
Since $\epsilon$ can be arbitrarily small, we have
\begin{align*}
\limsup\limits_{n\rightarrow \infty}\EE_{\rho_n}\left[\frac{n^4}{Z_1^n\left(\zeta_n\right)^4}\right] \le 1.
\end{align*}
By Jensen's inequality
\begin{align*}
\quad \liminf\limits_{n\rightarrow \infty}\EE_{\rho_n}\left[\frac{n^4}{Z_1^n\left(\zeta_n\right)^4}\right] \ge \liminf\limits_{n\rightarrow \infty} \frac{n^4}{\left(\EE_{\rho_n}\left[Z_1^n\left(\zeta_n\right)\right]\right)^4}= 1.
\end{align*}
Similar results can be obtained for $\EE_{\rho_n}\left[\frac{2n^2}{Z_1^n\left(\zeta_n\right)^2}\right]$. We conclude that 
\begin{align*}
\EE_{\rho_n}\left[\frac{n^4}{Z_1^n\left(\zeta_n\right)^4}-\frac{2n^2}{Z_1^n\left(\zeta_n\right)^2}+1\right]\rightarrow 0,
\end{align*}
and hence
\begin{align}\label{eqn:difference_zero}
    n^{1-\alpha}\EE_{\rho_n}\left[\left\lvert \tilde{R}_n-R_n\right\rvert\right]\rightarrow 0.
\end{align}
The desired result follows by combining \eqref{eqn:difference_zero} and \eqref{eqn:convergence_R_tilde}.
\qed
\end{proof}
\\

\section{Proofs of results in Section \ref{sec:num_clones_with_conditioning}}\label{Sec:proof_num_clones_with_conditioning}

\subsection{Proof of Proposition \ref{dist_1}}
\begin{proof}\\
Note that Proposition \ref{dist_1} does not rely on the assumption that sensitive cells have deterministic decay. Hence we provide the proof for the more general setting where $Z_0^n\left(t\right)$ is stochastic. We first observe that 
\begin{align}
& \quad \PP\left(I_n\left(\gamma_n\right)=k\right) \nonumber\\
& =\PP\left(I_n\left(\gamma_n\right)=k, I_n\left(\zeta_n\right)=k\right)+\PP\left(I_n\left(\gamma_n\right)=k, I_n\left(\zeta_n\right)\neq k\right)\nonumber\\
& =\PP\left(I_n\left(\zeta_n\right)=k\right)\nonumber \\
& \quad -\PP\left(I_n\left(\gamma_n\right)\neq k, I_n\left(\zeta_n\right)=k\right) \label{term_1_1}\\
& \quad +\PP\left(I_n\left(\gamma_n\right)=k, I_n\left(\zeta_n\right)\neq k\right) \label{term_1_2}.
\end{align}
For term (\ref{term_1_1}),
\begin{align*}
\PP\left(I_n\left(\gamma_n\right)\neq k, I_n\left(\zeta_n\right)=k\right)&=\PP\left(I_n\left(\gamma_n\right)\neq k, I_n\left(\zeta_n\right)=k, \gamma_n\notin \left({\zeta_n}-\delta, {\zeta_n}+\delta\right)\right)\\
& \quad +\PP\left(I_n\left(\gamma_n\right)\neq k, I_n\left(\zeta_n\right)=k, \gamma_n\in \left({\zeta_n}-\delta, {\zeta_n}+\delta\right)\right),
\end{align*}
where $\delta>0$. For the first probability, by Theorem $1$ of \cite{BP2020}, there exists $N_1>0$ such that when $n\ge N_1$,
\begin{align*}
& \quad \sum_{k=1}^{\infty}\PP\left(I_n\left(\gamma_n\right)\neq k, I_n\left(\zeta_n\right)=k, \gamma_n\notin \left({\zeta_n}-\delta, {\zeta_n}+\delta\right)\right)\\
& \le \PP\left(\gamma_n\notin \left({\zeta_n}-\delta, {\zeta_n}+\delta\right)\right)\\
& \le \frac{\epsilon}{4}.
\end{align*}
For the second probability, 
\begin{align*}
& \quad \sum_{k=1}^{\infty}\PP\left(I_n\left(\gamma_n\right)\neq k, I_n\left(\zeta_n\right)=k, \gamma_n\in \left({\zeta_n}-\delta, {\zeta_n}+\delta\right)\right)\\
& \le \PP\left(I_n\left(\gamma_n\right)\neq I_n\left(\zeta_n\right), \gamma_n\in \left({\zeta_n}-\delta, {\zeta_n}+\delta\right)\right)\\
& \le \PP\left(\exists t\in \left({\zeta_n}-\delta, {\zeta_n}+\delta\right) \text{ such that } I_n\left(t\right)\neq I_n\left(\zeta_n\right), \gamma_n\in \left({\zeta_n}-\delta, {\zeta_n}+\delta\right)\right)\\
& \le \PP\left(\text{mutation occurs in }\left({\zeta_n}-\delta, {\zeta_n}+\delta\right)\right).
\end{align*}
The expected number of mutations between times $\zeta_n-\delta$ and $\zeta_n+\delta$ can be written as
\begin{align*}
\int_{\zeta_n-\delta}^{\zeta_n+\delta} \mu n^{-\alpha} \EE\left[Z_0^n\left(t\right)\right]dt & \le 2\delta \cdot \mu n^{-\alpha} \cdot ne^{\lambda_0\left(\zeta_n-\delta\right)}\\
& \le Cn^{1-\alpha+\frac{\lambda_0}{\lambda_1}\alpha},
\end{align*}
where $C$ is a constant. Then by Markov's Inequality, there exists $N_2>0$ such that when $n\ge N_2$,
\begin{align*}
\PP\left(\text{mutation occurs in }\left({\zeta_n}-\delta, {\zeta_n}+\delta\right)\right)\le \frac{\epsilon}{4}
\end{align*}
Therefore, when $n\ge \max\left(N_1, N_2\right)$, $\sum_{k=1}^{\infty}\PP\left(I_n\left(\gamma_n\right)\neq k, I_n\left(\zeta_n\right)=k\right)\le \frac{\epsilon}{2}$. We can show a similar result for $\sum_{k=1}^{\infty}\PP\left(I_n\left(\gamma_n\right)=k, I_n\left(\zeta_n\right)\neq k\right)$ with the same reasoning which completes the proof. 
\qed
\end{proof}
\\

Prior to presenting the proof of Proposition \ref{proposition_diversity}, it is necessary to establish several preliminary results. We first consider the distribution of $I_n\left(\gamma_n\right)$ conditioned on the event of early recurrence. We show under assumption on $\alpha$ that $I_n\left(\gamma_n\right)$ is close to $I_n\left(\zeta_n-y\right)$ in total variation distance conditioned on the event of early recurrence.
\begin{proposition}\label{dist_2}
If $\alpha\in \left(\frac{\lambda_1}{\lambda_1-\lambda_0}\vee \frac{1}{2},1\right)$, then 
$$	
\lim_{n\to\infty}TV\left(I_n(\gamma_n)|\gamma_n<\zeta_n-y,I_n\left(\zeta_n-y\right)|\gamma_n<\zeta_n-y\right)=0.
$$
\end{proposition}
Note that via Bayes rule we can write
$$
\PP\left(I_n\left(\zeta_n-y\right)=k|\gamma_n<\zeta_n-y\right)=Q_{k,n} \PP\left(I_n\left(\zeta_n-y\right)=k\right),
$$
where	
\begin{equation*}
Q_{k,n}=\frac{\PP\left(\gamma_n<\zeta_n-y|I_{n}\left(\zeta_n-y\right)=k\right)}{\PP\left(\gamma_n<\zeta_n-y\right)}.
\end{equation*} 
By stochastic dominance, we can show that $Q_{k,n}$ increases in $k$. We then obtain the following result.
\begin{lemma}\label{lemma_diversity}
For all $x$,
\begin{align*}
\PP\left(I_n\left(\zeta_n-y\right)\ge x|\gamma_n<\zeta_n-y\right)\ge\PP\left(I_n\left(\zeta_n-y\right)\ge x\right).
\end{align*}
\end{lemma}
This result tells us that the number of clones at $\zeta_n-y$ conditioned on early recurrence stochastically dominates that without conditioning.


\subsection{Proof of Proposition \ref{dist_2}}
\begin{proof}\\
Note that Proposition \ref{dist_2} does not rely on the assumption that sensitive cells have deterministic decay. Hence we provide the proof for the more general setting where $Z_0^n\left(t\right)$ is stochastic. We first show the following lemma.
\begin{lemma}\label{lemma_3}
	Assume $\alpha\in \left(\frac{\lambda_1}{\lambda_1-\lambda_0},1\right)$. For any $M>0$,
	\begin{equation*}
	    \limsup\limits_{n\rightarrow \infty}\frac{1}{n^{\alpha}}\log\PP\left(\sup\limits_{t\in \left(\zeta_n-y-\delta, \zeta_n-y\right)}Z_0^n\left(s\right)>Mn^{\alpha}\right)<0.
	\end{equation*}
\end{lemma}
\begin{proof}\\
From \eqref{generating_function} and \eqref{bar_theta}, we know that for a fixed $0<\theta<\log\left(\frac{d_0}{r_0}\right)$, $\EE[e^{\theta Z_0^n\left(s\right)}]$ exists for $s\in \left(0, \zeta_n\right)$. 
We first show that for $M_1>0$, $\limsup\limits_{n\rightarrow \infty}\frac{1}{n^{\alpha}}\log\PP\left(Z_0^n\left(\zeta_n-y-\delta\right)\ge M_1n^{\alpha}\right)<0$. We observe that
\begin{align*}
\frac{1}{n^{\alpha}}\log\PP\left(Z_0^n\left(\zeta_n-y-\delta\right)\ge M_1n^{\alpha}\right)&=\frac{1}{n^{\alpha}}\log\PP\left(e^{\theta Z_0^n\left(\zeta_n-y-\delta\right)}\ge e^{M_1\theta n^{\alpha} }\right)\\
&\le \frac{1}{n^{\alpha}}\log\frac{\EE[e^{\theta Z_0^n\left(\zeta_n-y-\delta\right)}]}{e^{M_1\theta n^{\alpha}}}\\
& =n^{1-\alpha}\log\left(\frac{d_0\left(e^{\theta}-1\right)-e^{-\lambda_0 \left(\zeta_n-y-\delta\right)}\left(r_0e^{\theta}-d_0\right)}{r_0\left(e^{\theta}-1\right)-e^{-\lambda_0 \left(\zeta_n-y-\delta\right)}\left(r_0e^{\theta}-d_0\right)}\right)-M_1\theta.
\end{align*}
We can show that 
$$
n^{1-\alpha}\log\left(\frac{d_0\left(e^{\theta}-1\right)-e^{-\lambda_0 \left(\zeta_n-y-\delta\right)}\left(r_0e^{\theta}-d_0\right)}{r_0\left(e^{\theta}-1\right)-e^{-\lambda_0 \left(\zeta_n-y-\delta\right)}\left(r_0e^{\theta}-d_0\right)}\right)=O\left(n^{1+\frac{\lambda_0\alpha}{\lambda_1}-\alpha}\right).
$$ 
Since $\alpha>\frac{\lambda_1}{\lambda_1-\lambda_0}$,we have $\limsup\limits_{n\rightarrow \infty}\frac{1}{n^{\alpha}}\log\PP\left(Z_0^n\left(\zeta_n-y-\delta\right)\ge M_1n^{\alpha}\right)<0$. Next, we observe that for $M_2>M_1$, 
we have 
\begin{align*}
& \quad \frac{1}{n^{\alpha}}\log\PP\left(\sup\limits_{t\in \left(\zeta_n-y-\delta, \zeta_n-y\right)}Z_0^n\left(s\right)> M_2n^{\alpha}\middle|Z_0^n\left(\zeta_n-y-\delta\right)<M_1 n^{\alpha}\right)\\
& \le \frac{1}{n^{\alpha}}\log\frac{1-\left(d_0/r_0\right)^{M_1n^{\alpha}}}{1-\left(d_0/r_0\right)^{M_2n^{\alpha}}}\\
& \rightarrow \left(M_1-M_2\right)\log\left(\frac{d_0}{r_0}\right)\\
& <0,
\end{align*} 
where we apply a gambler's ruin argument to obtain the first inequality. The desired result follows by applying Lemma 1.2.15 in \cite{LD} which is restated in the following lemma.
\begin{lemma}\label{laplace_principle}
Let $M$ be a fixed integer. Then, for every $a^i_{\epsilon}\ge 0$,
\begin{align*}
\limsup\limits_{\epsilon\rightarrow 0} \epsilon\log \left(\sum_{i=1}^M a^i_{\epsilon}\right)=\max_{i\in \{1,...,M\}}\limsup\limits_{\epsilon\rightarrow 0} \epsilon \log a^i_{\epsilon}.
\end{align*}
\end{lemma}
\qed
\end{proof}
\\

We know that the total variation distance stated in Proposition \ref{dist_2} can be bounded above by
\begin{align*}
& \quad \PP\left(I_n\left(\gamma_n\right) \neq I_n\left(\zeta_n-y\right)| \gamma_n<\zeta_n-y\right)\\
= & \frac{\PP\left(I_n\left(\gamma_n\right) \neq I_n\left(\zeta_n-y\right), \gamma_n<\zeta_n-y\right)}{\PP\left(\gamma_n<\zeta_n-y\right)}\\
\le & \frac{\PP\left(\gamma_n<\zeta_n-y-\delta\right)}{\PP\left(\gamma_n<\zeta_n-y\right)}+\frac{\PP\left(I_n\left(\gamma_n\right) \neq I_n\left(\zeta_n-y\right), \zeta_n-y-\delta<\gamma_n<\zeta_n-y\right)}{\PP\left(\gamma_n<\zeta_n-y-\delta\right)}.
\end{align*}
The first term goes to zero by Theorem \ref{recurrence_LD}. The second term is bounded above by the conditional probability $\PP\left(I_n\left(\gamma_n\right) \neq I_n\left(\zeta_n-y\right) | \gamma_n\in \left(\zeta_n-y-\delta, \zeta_n-y\right)\right)$. By Lemma \ref{lemma_3}, Theorem \ref{recurrence_LD}, and the assumption that $\alpha>\frac{1}{2}$, we have
$$\PP\left(I_n\left(\gamma_n\right)\neq I_n\left(\zeta_n-y\right), \sup\limits_{t\in \left(\zeta_n-y-\delta, \zeta_n-y\right)}Z_0^n\left(s\right)>Mn^{\alpha}\middle|\gamma_n\in \left(\zeta_n-y-\delta, \zeta_n-y\right)\right)\rightarrow 0.
$$
We also have 
\begin{align*}
& \PP\left(I_n\left(\gamma_n\right)\neq I_n\left(\zeta_n-y\right), \sup\limits_{t\in \left(\zeta_n-y-\delta, \zeta_n-y\right)}Z_0^n\left(s\right)\le Mn^{\alpha}\middle|\gamma_n\in \left(\zeta_n-y-\delta, \zeta_n-y\right)\right)\\
\le & \PP\left(I_n\left(\gamma_n\right)\neq I_n\left(\zeta_n-y\right) \middle| \sup\limits_{t\in \left(\zeta_n-y-\delta, \zeta_n-y\right)}Z_0^n\left(s\right)\le Mn^{\alpha},\gamma_n\in \left(\zeta_n-y-\delta, \zeta_n-y\right)\right)\\
\le & 1-e^{-\mu M\delta}.
\end{align*}
Since $\delta$ can be arbitrarily small, the desired result follows.

\qed
\end{proof}
\\


\subsection{Proof of Lemma \ref{lemma_diversity}}
\begin{proof}\\
We first show that $Q_{k,n}$ increases in $k$. If we define $Z_{1,k}^n\left(t\right)$ as the number of mutants at time $t$ conditioned on the event of $\{I_{n}\left(\zeta_n-y\right)=k\}$, then it is easy to find a coupling for $Z_{1,k}^n\left(t\right)$ and $Z_{1,k+1}^n\left(t\right)$ such that for $t\in \left(0, \zeta_n-y\right)$, $Z_{1,k+1}^n\left(t\right)\ge Z_{1,k}^n\left(t\right)$, which indicates that
\begin{align*}
\PP\left(\gamma_n<\zeta_n-y|I_{n}\left(\zeta_n-y\right)=k+1\right) \ge \PP\left(\gamma_n<\zeta_n-y|I_{n}\left(\zeta_n-y\right)=k\right).
\end{align*}
It then follows that $Q_{k,n}$ increases in $k$. Since
\begin{align}\label{eqn: equality_condition}
\PP\left(I_n\left(\zeta_n-y\right)=k\middle|\gamma_n<\zeta_n-y\right)=Q_{k,n} \PP\left(I_n\left(\zeta_n-y\right)=k\right),
\end{align}
we now claim that $\sum\limits_{k=0}^{K}\PP\left(I_n\left(\zeta_n-y\right)=k\middle|\gamma_n<\zeta_n-y\right) \le \sum\limits_{k=0}^{K}  \PP\left(I_n\left(\zeta_n-y\right)=k\right)$ for any non-negative integer $K$ (which is equivalent to Lemma \ref{lemma_diversity}). Suppose (for a contradiction) that there exists $K>0$ such that 
\begin{align}\label{eqn: equality_condition_2}
\sum\limits_{k=0}^{K}\PP\left(I_n\left(\zeta_n-y\right)=k\middle|\gamma_n<\zeta_n-y\right) > \sum\limits_{k=0}^{K}  \PP\left(I_n\left(\zeta_n-y\right)=k\right).
\end{align}
Then we must have $Q_{K,n}>1$ due to \eqref{eqn: equality_condition} and the fact that $Q_{k,n}$ increases in $k$. Moreover, for any $k\ge K$, $Q_{k,n}>1$, which indicates that $\PP\left(I_n\left(\zeta_n-y\right)=k\middle|\gamma_n<\zeta_n-y\right) > \PP\left(I_n\left(\zeta_n-y\right)=k\right)$ for any $k\ge K$ by \eqref{eqn: equality_condition}. Combining this result with \eqref{eqn: equality_condition_2}, we can obtain that
$\sum\limits_{k=0}^{\infty}\PP\left(I_n\left(\zeta_n-y\right)=k|\gamma_n<\zeta_n-y\right)>\sum\limits_{k=0}^{\infty}  \PP\left(I_n\left(\zeta_n-y\right)=k\right)=1$ (a contradiction). The desired result then follows.
\qed
\end{proof}
\\


By utilizing Lemma \ref{lemma_diversity}, Proposition \ref{dist_1} and Proposition \ref{dist_2}, we are able to prove Proposition \ref{proposition_diversity}.

\subsection{Proof of Proposition \ref{proposition_diversity}}

\begin{proof}\\
Given $x>0$, we have
\begin{align*}
& \quad \PP\left(I_n\left(\gamma_n\right)\ge x\right)-\PP\left(I_n\left(\gamma_n\right)\ge x\middle|\gamma_n<\zeta_n-y\right)\\
& = \PP\left(I_n\left(\gamma_n\right)\ge x\right)-\PP\left(I_n\left(\zeta_n-y\right)\ge x\right)\\
& \quad \quad + \PP\left(I_n\left(\zeta_n-y\right)\ge x|\gamma_n<\zeta_n-y\right)-\PP\left(I_n\left(\gamma_n\right)\ge x|\gamma_n<\zeta_n-y\right)\\
& \quad \quad +\PP\left(I_n\left(\zeta_n-y\right)\ge x\right)-\PP\left(I_n\left(\zeta_n-y\right)\ge x|\gamma_n<\zeta_n-y\right)\\
& \le TV\left(I_n\left(\gamma_n\right),I_n\left(\zeta_n\right)\right)+TV\left(I_n\left(\zeta_n-y\right),I_n\left(\zeta_n\right)\right)+TV\left(I_n(\gamma_n)|\gamma_n<\zeta_n-y,I_n\left(\zeta_n-y\right)|\gamma_n<\zeta_n-y\right)\\
& \quad \quad +\PP\left(I_n\left(\zeta_n-y\right)\ge x\right)-\PP\left(I_n\left(\zeta_n-y\right)\ge x|\gamma_n<\zeta_n-y\right).
\end{align*}
We know from Lemma \ref{lemma_diversity} that for all $x$, 
\begin{align*}
\PP\left(I_n\left(\zeta_n-y\right)\ge x|\gamma_n<\zeta_n-y\right)\ge\PP\left(I_n\left(\zeta_n-y\right)\ge x\right).
\end{align*}
By Proposition \ref{dist_1} and Proposition \ref{dist_2}, we have
\begin{align*}
& \lim_{n\to\infty}TV\left(I_n\left(\gamma_n\right),I_n\left(\zeta_n\right)\right)=0, \text{ and}\\
& \lim_{n\to\infty}TV\left(I_n(\gamma_n)|\gamma_n<\zeta_n-y,I_n\left(\zeta_n-y\right)|\gamma_n<\zeta_n-y\right)=0.
\end{align*}
Hence, it remains to show that 
\begin{align}\label{dist_3}
& \lim_{n\to\infty}TV\left(I_n\left(\zeta_n-y\right),I_n\left(\zeta_n\right)\right)=0.
\end{align}
We know that 
\begin{align*}
TV\left(I_n\left(\zeta_n-y\right),I_n\left(\zeta_n\right)\right)& \le \PP\left(\text{mutation occurs in }\left({\zeta_n}-y, {\zeta_n}\right)\right).
\end{align*}
By a similar argument to that in the proof of Lemma \ref{lemma_3}, we can show that 
$$
\PP\left(\sup\limits_{s\in \left(\zeta_n-y, \zeta_n\right)}Z_0^n\left(s\right)>Mn^{\alpha^*}\right)
$$ 
decays exponentially fast for $\frac{\lambda_1}{\lambda_1-\lambda_0}<\alpha^*<\alpha$, and the desired result follows.
\qed
\end{proof}
\\

\subsection{Proof of Lemma \ref{expectation_with_conditioning}}\label{lemma_expectation_with_conditioning}

We omit the detailed proof here as Lemma \ref{expectation_with_conditioning} is implied by Lemma \ref{dominated_event} in Section \ref{Sec:proof_simpson_with_conditioning}.

\subsection{Proof of Lemma \ref{period_mutants_concentration}}\label{lemma_period_mutants_concentration}
 
\begin{proof}\\
For $0<\delta_1<\delta_2$,
\begin{align}
\PP_{A_{n,y}}\left(\frac{B_{n,y,\left(t_1, t_2\right)}}{\bar
{B}_{n,y,\left(t_1, t_2\right)}}\in \left(1+\delta_1, 1+\delta_2\right)\right)= \frac{\PP\left(\frac{B_{n,y,\left(t_1, t_2\right)}}{\bar
{B}_{n,y,\left(t_1, t_2\right)}}\in \left(1+\delta_1, 1+\delta_2\right), A_{n,y}\right)}{\PP\left(A_{n,y}\right)}. \label{period_mutants_conditioning}
\end{align}
We have the following upper and lower bound for the numerator in \eqref{period_mutants_conditioning}:
{\footnotesize
\begin{align}
&  \PP\left(\frac{B_{n,y,\left(t_1, t_2\right)}}{\bar
{B}_{n,y,\left(t_1, t_2\right)}}\in \left(1+\delta_1, 1+\delta_2\right), A_{n,y}\right) \nonumber\\
\le & \PP\left(\frac{B_{n,y,\left(t_1, t_2\right)}}{\bar
{B}_{n,y,\left(t_1, t_2\right)}}\in \left(1+\delta_1, 1+\delta_2\right), B_{n,y,\left(0, \zeta_n-y\right)} - B_{n,y,\left(t_1, t_2\right)}> n-\left(1+\delta_2\right)\bar
{B}_{n,y,\left(t_1, t_2\right)}\right), \label{period_mutants_upper}
\end{align}}
and
{\footnotesize
\begin{align}
& \PP\left(\frac{B_{n,y,\left(t_1, t_2\right)}}{\bar
{B}_{n,y,\left(t_1, t_2\right)}}\in \left(1+\delta_1, 1+\delta_2\right), A_{n,y}\right) \nonumber\\
\ge & \PP\left(\frac{B_{n,y,\left(t_1, t_2\right)}}{\bar
{B}_{n,y,\left(t_1, t_2\right)}}\in \left(1+\delta_1, 1+\delta_2\right), B_{n,y,\left(0, \zeta_n-y\right)} - B_{n,y,\left(t_1, t_2\right)}> n-\left(1+\delta_1\right)\bar
{B}_{n,y,\left(t_1, t_2\right)}\right). \label{period_mutants_lower}
\end{align}}
By the Gartner-Ellis Theorem (\cite{Hollander}) and a similar calculation to that in the proof of Theorem 2 in \cite{BP2020} (calculation of the moment generating function of $B_{n,y,\left(t_1, t_2\right)}$), we can obtain the large deviations rate for the probability presented in \eqref{period_mutants_lower} by considering the following large deviation rates for $\delta,\varepsilon > 0$: 
\begin{align*}
L_1\left(\delta\right)& \triangleq - \lim\limits_{n\rightarrow \infty}\frac{1}{n^{1-\alpha}}\log \PP\left(\frac{B_{n,y,\left(t_1, t_2\right)}}{\bar
{B}_{n,y,\left(t_1, t_2\right)}}\in \left(1+\delta, 1+\delta+\varepsilon\right)\right)\\
& =\sup_{\theta\in \left(0,1\right)}\left[\mu\theta\left(1+\delta\right)\int_{t_1}^{t_2} e^{-\left(\lambda_1-\lambda_0\right)s}ds - \mu\int_{t_1}^{t_2}\frac{\theta}{e^{\lambda_1 s}-\theta}e^{\lambda_0 s}ds\right],
\end{align*}
and
\begin{align}
L_2\left(\delta\right)& \triangleq - \lim\limits_{n\rightarrow \infty}\frac{1}{n^{1-\alpha}}\log \PP\left(B_{n,y,\left(0, \zeta_n-y\right)}-B_{n,y,\left(t_1, t_2\right)}> n-\left(1+\delta\right)\bar
{B}_{n,y,\left(t_1, t_2\right)}\right)\label{L_2_value_1}\\
& =\sup_{\theta\in \left(0,1\right)}\left[\frac{\mu \theta e^{\lambda_1 y}}{\lambda_1-\lambda_0}-\mu\int_{0}^{t_1}\frac{\theta}{e^{\lambda_1 s}-\theta}e^{\lambda_0 s}ds - \mu\int_{t_2}^{\infty}\frac{\theta}{e^{\lambda_1 s}-\theta}e^{\lambda_0 s}ds -\mu\theta\left(1+\delta\right)\int_{t_1}^{t_2} e^{-\left(\lambda_1-\lambda_0\right)s}ds\right].\label{L_2_value_2}
\end{align}
Since the generation of clones in different non-overlapping time periods are independent (the evolution of each clone is also independent), the large deviations rate for the probability presented in \eqref{period_mutants_lower} is given by
\begin{align*}
L_1\left(\delta_1\right)+L_2\left(\delta_1\right).
\end{align*}
We then restate the Envelope Theorem (see page 158 of \cite{Bordernotes}) which can be used to analyze $L_1\left(\delta_1\right)+L_2\left(\delta_1\right)$.
\begin{theorem}[Envelope Theorem]\label{envelope theorem}
Let $X$ be a metric space and $P$ an open subset of $\RR^n$. Let $w: X\times P\rightarrow \RR$ and assume $\frac{\partial w}{\partial p}$ exists and is continuous in $X\times P$. For each $p\in P$, let $x^*\left(p\right)$ maximize $w\left(x,p\right)$ over $X$. Set
\begin{align*}
    V\left(p\right)=w\left(x^*\left(p\right),p\right).
\end{align*}
Assume that $x^*: P\rightarrow X$ is a continuous function. Then $V$ is continuously differentiable and 
\begin{align*}
DV\left(p\right)=\frac{\partial w \left(x,p\right)}{\partial p},
\end{align*}
where the derivative is evaluated at the point $\left(x^*\left(p\right),p\right)$.
\end{theorem}
By the Envelope Theorem, we can obtain that $L_1\left(\delta\right)$ is convex, increasing in $\delta$, and $L_2\left(\delta\right)$ is convex, decreasing in $\delta$. Hence, by standard convex analysis, we can obtain the optimal $\delta^*$ to the following optimization problem:
\begin{align}
\min_{\delta>0} \left(L_1\left(\delta\right)+L_2\left(\delta\right) \right). \label{period_mutants_optimization}
\end{align}
In particular,
\begin{align*}
\frac{dL_1\left(\delta\right)}{d\delta}& =\mu \theta_{1, \delta}\int_{t_1}^{t_2} e^{-\left(\lambda_1-\lambda_0\right)s}ds, \text{ and }\\
\frac{dL_2\left(\delta\right)}{d\delta}& =-\mu \theta_{2, \delta}\int_{t_1}^{t_2} e^{-\left(\lambda_1-\lambda_0\right)s}ds,
\end{align*}
where $\theta_{1, \delta}$ ($\theta_{2, \delta}$) is the maximizer to the optimization problem within the expression of $L_1\left(\delta\right)$ ($L_2\left(\delta\right)$).
Therefore, the minimum of $\left(L_1\left(\delta\right)+L_2\left(\delta\right) \right)$ is achieved when $\theta_{1, \delta}=\theta_{2, \delta}$. By some calculation, we can obtain that $\theta_{1, \delta^*}=\theta_{2, \delta^*}=\theta_y^*$, where $\delta^*$ is given by
\begin{align*}
\delta^*=\frac{\int_{t_1}^{t_2}\frac{e^{\lambda_1s}}{\left(e^{\lambda_1 s}-\theta^*_y\right)^2}e^{\lambda_0 s}ds}{\int_{t_1}^{t_2}e^{-\left(\lambda_1-\lambda_0\right)s}ds}-1,
\end{align*}
and $\theta^*_y$ is defined in \eqref{theta_2}. From \eqref{period_mutants_lower}, we know that
\begin{align*}
& \liminf\limits_{n\rightarrow \infty}\frac{1}{n^{1-\alpha}}\log\PP\left(\frac{B_{n,y,\left(t_1, t_2\right)}}{\bar
{B}_{n,y,\left(t_1, t_2\right)}}\in \left(1+\delta^*-\epsilon, 1+\delta^*+\epsilon\right), A_{n,y}\right)\\
\ge & \liminf\limits_{n\rightarrow \infty}\frac{1}{n^{1-\alpha}}\log\PP\left(\frac{B_{n,y,\left(t_1, t_2\right)}}{\bar
{B}_{n,y,\left(t_1, t_2\right)}}\in \left(1+\delta^*, 1+\delta^*+\epsilon\right), A_{n,y}\right)\\
\ge & -\left(L_1\left(\delta^*\right)+L_2\left(\delta^*\right)\right).
\end{align*}
To obtain the desired result, it suffices to show that
\begin{align*}
\liminf\limits_{n\rightarrow \infty}\frac{1}{n^{1-\alpha}}\log\PP\left(\frac{B_{n,y,\left(t_1, t_2\right)}}{\bar
{B}_{n,y,\left(t_1, t_2\right)}}\notin \left(1+\delta^*-\epsilon, 1+\delta^*+\epsilon\right), A_{n,y}\right) < -\left(L_1\left(\delta^*\right)+L_2\left(\delta^*\right)\right).
\end{align*}
We first investigate the event
$$
\Biggl\{\frac{B_{n,y,\left(t_1, t_2\right)}}{\bar
{B}_{n,y,\left(t_1, t_2\right)}}\in \left(1+\delta^*+\epsilon, \frac{e^{\lambda_1 y}}{\left(\lambda_1-\lambda_0\right)\int_{t_1}^{t_2}e^{-\left(\lambda_1-\lambda_0\right)s}ds}\right),A_{n,y}\Biggr\},
$$
where
\begin{align*}
\frac{e^{\lambda_1 y}}{\left(\lambda_1-\lambda_0\right)\int_{t_1}^{t_2}e^{-\left(\lambda_1-\lambda_0\right)s}ds} = \lim\limits_{n\rightarrow \infty}\frac{n}{\bar
{B}_{n,y,\left(t_1, t_2\right)}}.
\end{align*}
We note that
\begin{align*}
& \limsup\limits_{n\rightarrow \infty}\frac{1}{n^{1-\alpha}}\PP\left(\frac{B_{n,y,\left(t_1, t_2\right)}}{\bar
{B}_{n,y,\left(t_1, t_2\right)}}\in \left(\frac{e^{\lambda_1 y}}{\left(\lambda_1-\lambda_0\right)\int_{t_1}^{t_2}e^{-\left(\lambda_1-\lambda_0\right)s}ds}, \infty\right), A_{n,y}\right)\\
\le & \limsup\limits_{n\rightarrow \infty}\frac{1}{n^{1-\alpha}}\PP\left(\frac{B_{n,y,\left(t_1, t_2\right)}}{\bar
{B}_{n,y,\left(t_1, t_2\right)}}\in \left(\frac{e^{\lambda_1 y}}{\left(\lambda_1-\lambda_0\right)\int_{t_1}^{t_2}e^{-\left(\lambda_1-\lambda_0\right)s}ds}, \infty\right)\right)\\
\le & -L_1\left(\frac{e^{\lambda_1 y}}{\left(\lambda_1-\lambda_0\right)\int_{t_1}^{t_2}e^{-\left(\lambda_1-\lambda_0\right)s}ds}-1\right)-L_2\left(\frac{e^{\lambda_1 y}}{\left(\lambda_1-\lambda_0\right)\int_{t_1}^{t_2}e^{-\left(\lambda_1-\lambda_0\right)s}ds}-1\right)\\
< & -\left(L_1\left(\delta^*\right)+L_2\left(\delta^*\right)\right),
\end{align*}
where the second inequality is due to the fact that 
$$
L_2\left(\frac{e^{\lambda_1 y}}{\left(\lambda_1-\lambda_0\right)\int_{t_1}^{t_2}e^{-\left(\lambda_1-\lambda_0\right)s}ds}-1\right)=0,
$$ 
and the last inequality is because $\delta^*$ is the optimal solution to \eqref{period_mutants_optimization}. Therefore,
we could safely omit the event 
$$
\Biggl\{\frac{B_{n,y,\left(t_1, t_2\right)}}{\bar
{B}_{n,y,\left(t_1, t_2\right)}}\in \left(\frac{e^{\lambda_1 y}}{\left(\lambda_1-\lambda_0\right)\int_{t_1}^{t_2}e^{-\left(\lambda_1-\lambda_0\right)s}ds}, \infty\right),A_{n,y}\Biggr\},
$$
and focus on the interval
\begin{align*}
\left(1+\delta^*+\epsilon, \frac{e^{\lambda_1 y}}{\left(\lambda_1-\lambda_0\right)\int_{t_1}^{t_2}e^{-\left(\lambda_1-\lambda_0\right)s}ds}\right).    
\end{align*}

For an arbitrarily large but fixed integer $M$, we divide the interval 
$$
\left(1+\delta^*+\epsilon, \frac{e^{\lambda_1 y}}{\left(\lambda_1-\lambda_0\right)\int_{t_1}^{t_2}e^{-\left(\lambda_1-\lambda_0\right)s}ds}\right)
$$ 
into $M$ sub-intervals with equal length $\delta_M$. 
Since $M$ is a fixed number, we could analyze the event that $\frac{B_{n,y,\left(t_1, t_2\right)}}{\bar
{B}_{n,y,\left(t_1, t_2\right)}}$ falls in each sub-interval, and then apply Lemma \ref{laplace_principle}. 
By Lemma \ref{laplace_principle}, to obtain the desired result, it suffices to show that for all $i\in \{0,...,M-1\}$,
\begin{align}
& \limsup\limits_{n\rightarrow \infty}\frac{1}{n^{1-\alpha}}\log \PP\left(\frac{B_{n,y,\left(t_1, t_2\right)}}{\bar
{B}_{n,y,\left(t_1, t_2\right)}}\in \left(1+\delta^*+\epsilon+i \delta_M, 1+\delta^*+\epsilon+\left(i+1\right)\delta_M\right), A_{n,y}\right) \nonumber\\
< & -\left(L_1\left(\delta^*\right)+L_2\left(\delta^*\right)\right). \label{inequality_L1L2_B}
\end{align}
From \eqref{period_mutants_upper}, it suffices to show that
\begin{align}
L_1\left(\delta^*+\epsilon+i\delta_M\right)+L_2\left(\delta^*+\epsilon+\left(i+1\right)\delta_M\right)>L_1\left(\delta^*\right)+L_2\left(\delta^*\right),\label{ineq:L1L2}
\end{align}
Since $L_2\left(\delta\right)$ is continuous for 
$$
\delta\in \left[\delta^*+\epsilon, \frac{e^{\lambda_1 y}}{\left(\lambda_1-\lambda_0\right)\int_{t_1}^{t_2}e^{-\left(\lambda_1-\lambda_0\right)s}ds}-1\right]
$$ 
and thus also uniform continuous, and $\delta^*$ is the optimal solution to \eqref{period_mutants_optimization}, we conclude that for sufficiently large $M$, \eqref{ineq:L1L2} holds
for all $i\in \{0,...,M-1\}$ which validates \eqref{inequality_L1L2_B}. 

By a similar argument, we can deal with the event
$$
\Biggl\{\frac{B_{n,y,\left(t_1, t_2\right)}}{\bar
{B}_{n,y,\left(t_1, t_2\right)}}\in \left(0, 1+\delta^*-\epsilon\right),A_{n,y}\Biggr\},
$$
which completes the proof.
\qed
\end{proof}
\\


\subsection{Proof of Lemma \ref{period_clones_concentration}}\label{lemma_period_clones_concentration}
\begin{proof}\\
We first notice that for $0<\kappa_1<\kappa_2<\delta_1<\delta_2$,
{\footnotesize
\begin{align}
&  \PP\left(\frac{B_{n,y,\left(t_1, t_2\right)}}{\bar
{B}_{n,y,\left(t_1, t_2\right)}}\in \left(1+\delta_1, 1+\delta_2\right), \frac{I_{n, \left(t_1, t_2\right)}}{\bar{I}_{n, \left(t_1, t_2\right)}}\in \left(1+\kappa_1, 1+\kappa_2\right)\right) \nonumber\\
\le &  \PP\left(\frac{B_{n,y,\left(t_1, t_2\right)}}{\bar
{B}_{n,y,\left(t_1, t_2\right)}}\in \left(1+\delta_1, \infty\right), \frac{I_{n, \left(t_1, t_2\right)}}{\bar{I}_{n, \left(t_1, t_2\right)}}\in \left(1+\kappa_1, 1+\kappa_2\right)\right) \nonumber\\
\le & \PP\left(\frac{B_{n,y,\left(t_1, t_2\right)}}{\bar
{B}_{n,y,\left(t_1, t_2\right)}}\in \left(1+\delta_1, \infty\right)\middle| I_{n, \left(t_1, t_2\right)}= \lfloor\left(1+\kappa_2\right)\bar{I}_{n, \left(t_1, t_2\right)}\rfloor\right)\nonumber\\
& \quad \quad \quad \times \PP\left(\frac{I_{n, \left(t_1, t_2\right)}}{\bar{I}_{n, \left(t_1, t_2\right)}}\in \left(1+\kappa_1, 1+\kappa_2\right)\right),\label{period_clones_upper}
\end{align}}
where the last inequality is due to stochastic dominance. With a similar reasoning, we also have
{\footnotesize
\begin{align}
& \PP\left(\frac{B_{n,y,\left(t_1, t_2\right)}}{\bar
{B}_{n,y,\left(t_1, t_2\right)}}\in \left(1+\delta_1, 1+\delta_2\right), \frac{I_{n, \left(t_1, t_2\right)}}{\bar{I}_{n, \left(t_1, t_2\right)}}\in \left(1+\kappa_1, 1+\kappa_2\right)\right)\nonumber\\
= & \PP\left(\frac{B_{n,y,\left(t_1, t_2\right)}}{\bar
{B}_{n,y,\left(t_1, t_2\right)}}\in \left(1+\delta_1, \infty\right), \frac{I_{n, \left(t_1, t_2\right)}}{\bar{I}_{n, \left(t_1, t_2\right)}}\in \left(1+\kappa_1, 1+\kappa_2\right)\right)\nonumber\\
& -\PP\left(\frac{B_{n,y,\left(t_1, t_2\right)}}{\bar
{B}_{n,y,\left(t_1, t_2\right)}}\in \left(1+\delta_2, \infty\right), \frac{I_{n, \left(t_1, t_2\right)}}{\bar{I}_{n, \left(t_1, t_2\right)}}\in \left(1+\kappa_1, 1+\kappa_2\right)\right)\nonumber\\
\ge & \PP\left(\frac{B_{n,y,\left(t_1, t_2\right)}}{\bar
{B}_{n,y,\left(t_1, t_2\right)}}\in \left(1+\delta_1, \infty\right)\middle| I_{n, \left(t_1, t_2\right)}= \lfloor\left(1+\kappa_1\right)\bar{I}_{n, \left(t_1, t_2\right)}\rfloor\right)\nonumber\\
& \quad \quad \quad \times \PP\left(\frac{I_{n, \left(t_1, t_2\right)}}{\bar{I}_{n, \left(t_1, t_2\right)}}\in \left(1+\kappa_1, 1+\kappa_2\right)\right) \nonumber\\
& - \PP\left(\frac{B_{n,y,\left(t_1, t_2\right)}}{\bar
{B}_{n,y,\left(t_1, t_2\right)}}\in \left(1+\delta_2, \infty\right)\middle| I_{n, \left(t_1, t_2\right)}= \lfloor\left(1+\kappa_2\right)\bar{I}_{n, \left(t_1, t_2\right)}\rfloor\right)\nonumber\\
& \quad \quad \quad \times \PP\left(\frac{I_{n, \left(t_1, t_2\right)}}{\bar{I}_{n, \left(t_1, t_2\right)}}\in \left(1+\kappa_1, 1+\kappa_2\right)\right). \label{period_clones_lower} 
\end{align}}
By considering the sum of $\lfloor\left(1+\kappa_1\right)\bar{I}_{n, \left(t_1, t_2\right)}\rfloor$ i.i.d. random variables, each of which is the number of descendants of one mutated cell, we can compute (using the Cramér's theorem) that
{\footnotesize
\begin{align*}
& L_1\left(\delta_1, \kappa_1\right) \\
= & -\lim\limits_{n\rightarrow \infty}\frac{1}{n^{1-\alpha}}\log \PP\left(\frac{B_{n,y,\left(t_1, t_2\right)}}{\bar
{B}_{n,y,\left(t_1, t_2\right)}}\in \left(1+\delta_1, 1+\delta_2\right)\middle| I_{n, \left(t_1, t_2\right)}= \lfloor\left(1+\kappa_1\right)\bar{I}_{n, \left(t_1, t_2\right)}\rfloor\right)\\
= &\sup_{\theta\in \left(0,1\right)}\left[\mu\theta\left(1+\delta_1\right)\int_{t_1}^{t_2} e^{-\left(\lambda_1-\lambda_0\right)s}ds - \mu\left(1+\kappa_1\right)\int_{t_1}^{t_2}e^{\lambda_0 s}ds \log\left(\frac{\int_{t_1}^{t_2}\frac{e^{\lambda_1 s}}{e^{\lambda_1 s}-\theta}e^{\lambda_0 s}ds}{\int_{t_1}^{t_2}e^{\lambda_0 s}ds}\right)\right].
\end{align*}}
Because the distribution of $I_{n, \left(t_1, t_2\right)}$ is Poisson with a mean of $\bar{I}_{n, \left(t_1, t_2\right)}$, given in \eqref{def:bar_I}, we can compute that
\begin{align*}
L_2\left(\kappa_1\right) &= -\lim\limits_{n\rightarrow \infty}\frac{1}{n^{1-\alpha}}\log \PP\left(\frac{I_{n, \left(t_1, t_2\right)}}{\bar{I}_{n, \left(t_1, t_2\right)}}\in \left(1+\kappa_1, 1+\kappa_2\right)\right)\\
& =\mu \int_{t_1}^{t_2}e^{\lambda_0 s}ds\left(\left(1+\kappa_1\right)\log\left(1+\kappa_1\right)-\left(1+\kappa_1\right)+1\right).
\end{align*}
Consider the optimization problem:
\begin{align}\label{period_clones_optimization}
\min_{\kappa>0} \left(L_1\left(\delta^*, \kappa\right)+L_2\left(\kappa\right)\right).
\end{align}
By the Envelope Theorem, we can show that $L_1\left(\delta^*, \kappa\right)$ is convex, decreasing in $\kappa$ and $L_2\left(\kappa\right)$ is convex, increasing in $\kappa$. Hence, by standard convex analysis, if $t_2-t_1< -\frac{1}{\lambda_1}\log\left(\frac{1}{2-\theta^*_y}\right)$ (this condition guarantees that the optimal solution lies in $\left(0,\delta^*\right)$), the optimal value $\kappa^*$ is given by
\begin{align*}
\kappa^*=\frac{\int_{t_1}^{t_2}\frac{e^{\lambda_1 s}}{e^{\lambda_1 s}-\theta^*_y}e^{\lambda_0 s}ds}{\int_{t_1}^{t_2}e^{\lambda_0 s}ds}-1,
\end{align*}
where $\theta^*_y$ is defined in \eqref{theta_2}.

By a similar argument to that in the proof of Lemma \ref{period_mutants_concentration}, we can focus on the event
\begin{align*}
\Biggl\{\frac{B_{n,y,\left(t_1, t_2\right)}}{\bar
{B}_{n,y,\left(t_1, t_2\right)}}\in \left(1+\delta^*-\epsilon, 1+\delta^*+\epsilon\right), \frac{I_{n, \left(t_1, t_2\right)}}{\bar{I}_{n, \left(t_1, t_2\right)}}\in \left(1, 1+\delta^*-\epsilon\right)\Biggr\}.    
\end{align*}
From \eqref{period_clones_lower}, we have
{\footnotesize
\begin{align}
& \PP\left(\frac{B_{n,y,\left(t_1, t_2\right)}}{\bar
{B}_{n,y,\left(t_1, t_2\right)}}\in \left(1+\delta^*-\epsilon, 1+\delta^*+\epsilon\right), \frac{I_{n, \left(t_1, t_2\right)}}{\bar{I}_{n, \left(t_1, t_2\right)}}\in \left(1+\kappa^*-\sigma, 1+\kappa^*+\sigma\right)\right) \label{def_main_event}\\
\ge & \PP\left(\frac{B_{n,y,\left(t_1, t_2\right)}}{\bar
{B}_{n,y,\left(t_1, t_2\right)}}\in \left(1+\delta^*-\epsilon, \infty\right)\middle| I_{n, \left(t_1, t_2\right)}= \lfloor\left(1+\kappa^*-\sigma\right)\bar{I}_{n, \left(t_1, t_2\right)}\rfloor\right)\nonumber\\
& \quad \quad \quad \times \PP\left(\frac{I_{n, \left(t_1, t_2\right)}}{\bar{I}_{n, \left(t_1, t_2\right)}}\in \left(1+\kappa^*-\sigma, 1+\kappa^*+\sigma\right)\right) \nonumber\\
& - \PP\left(\frac{B_{n,y,\left(t_1, t_2\right)}}{\bar
{B}_{n,y,\left(t_1, t_2\right)}}\in \left(1+\delta^*+\epsilon, \infty\right)\middle| I_{n, \left(t_1, t_2\right)}= \lfloor\left(1+\kappa^*+\sigma\right)\bar{I}_{n, \left(t_1, t_2\right)}\rfloor\right)\nonumber\\
& \quad \quad \quad \times \PP\left(\frac{I_{n, \left(t_1, t_2\right)}}{\bar{I}_{n, \left(t_1, t_2\right)}}\in \left(1+\kappa^*-\sigma, 1+\kappa^*+\sigma\right)\right).\nonumber
\end{align}}
To use this lower bound, we need to make sure that the negative part decays exponentially faster. By the calculation of $L_1$, we can obtain that for a positive number $N>\frac{\sigma \left(1+\delta^*\right)}{\epsilon\left(1+\kappa^*\right)}$, we have $\frac{1+\delta^*+\epsilon}{1+\delta^*-\epsilon}>\frac{1+\kappa^*+\frac{\sigma}{N}}{1+\kappa^*-\frac{\sigma}{N}}$ which indicates that $L_1\left(\delta^*-\epsilon, \kappa^*-\frac{\sigma}{N}\right)< L_1\left(\delta^*+\epsilon, \kappa^*+\frac{\sigma}{N}\right)$. Therefore,
{\footnotesize
\begin{align*}
& \liminf\limits_{n\rightarrow \infty}\frac{1}{n^{1-\alpha}}\log \PP\left(\frac{B_{n,y,\left(t_1, t_2\right)}}{\bar
{B}_{n,y,\left(t_1, t_2\right)}}\in \left(1+\delta^*-\epsilon, 1+\delta^*+\epsilon\right), \frac{I_{n, \left(t_1, t_2\right)}}{\bar{I}_{n, \left(t_1, t_2\right)}}\in \left(1+\kappa^*-\sigma, 1+\kappa^*+\sigma\right)\right)\nonumber\\
\ge & \liminf\limits_{n\rightarrow \infty}\frac{1}{n^{1-\alpha}}\log \PP\left(\frac{B_{n,y,\left(t_1, t_2\right)}}{\bar
{B}_{n,y,\left(t_1, t_2\right)}}\in \left(1+\delta^*-\epsilon, 1+\delta^*+\epsilon\right), \frac{I_{n, \left(t_1, t_2\right)}}{\bar{I}_{n, \left(t_1, t_2\right)}}\in \left(1+\kappa^*-\frac{\sigma}{N}, 1+\kappa^*+\frac{\sigma}{N}\right)\right)\nonumber\\
\ge & -\left(L_1\left(\delta^*-\epsilon, \kappa^*-\frac{\sigma}{N}\right)+L_2\left(\kappa^*-\frac{\sigma}{N}\right)\right).
\end{align*}}
We first divide the interval $\left(0, \kappa^*-\sigma\right)$ into $M$ sub-intervals with equal length $\delta_M$ for an arbitrarily large but fixed integer $M>0$. By Lemma \ref{laplace_principle}, to obtain the desired result, the first step is to show that for all $i\in \{0,...,M-1\}$,
{\small
\begin{align*}
& \limsup\limits_{n\rightarrow \infty}\frac{1}{n^{1-\alpha}}\log \PP\left(\frac{I_{n,\left(t_1, t_2\right)}}{\bar
{I}_{n,\left(t_1, t_2\right)}}\in \left(1+i \delta_M, 1+\left(i+1\right)\delta_M\right), A_{n,y,\epsilon}^{t_1,t_2}\right)\\
< & -\left(L_1\left(\delta^*-\epsilon, \kappa^*-\frac{\sigma}{N}\right)+L_2\left(\kappa^*-\frac{\sigma}{N}\right)\right).
\end{align*}}
From \eqref{period_clones_upper}, it suffices to show that
\begin{align}
L_1\left(\delta^*-\epsilon, \left(i+1\right)\delta_M\right)+L_2\left(i\delta_M\right)>L_1\left(\delta^*-\epsilon, \kappa^*-\frac{\sigma}{N}\right)+L_2\left(\kappa^*-\frac{\sigma}{N}\right).\label{final_inequality_time_period}
\end{align}
We first analyze the left hand side of \eqref{final_inequality_time_period}. 
Since $L_1\left(\delta^*, \kappa\right)$ is continuous for $\kappa\in \left[0, \delta^*-\epsilon\right]$, and thus also uniformly continuous, we have
\begin{align}
\lim\limits_{M\rightarrow \infty}\sup\limits_{i\in \{0,...,M-1\}}\left|L_1\left(\delta^*, i\delta_M\right)-L_1\left(\delta^*, \left(i+1\right)\delta_M\right)\right|=0.\label{L_1_uniform_continuous}
\end{align}
Let $\theta^*_y\left(\delta_1, \kappa_1\right)$ be the optimal value for the inner optimization problem appearing in the calculation of $L_1\left(\delta_1, \kappa_1\right)$. We can see that $\theta^*_y\left(\delta_1, \kappa_1\right)$ decreases in $\kappa_1$. Hence, by the Envelope Theorem,
\begin{align*}
\frac{\partial L_1\left(\delta_1, \kappa_1\right)}{\partial \delta_1}=\mu\theta^*_y\left(\delta_1, \kappa_1\right)\int_{t_1}^{t_2} e^{-\left(\lambda_1-\lambda_0\right)s}ds\le \mu\theta^*_y\left(\delta_1, 0\right)\int_{t_1}^{t_2} e^{-\left(\lambda_1-\lambda_0\right)s}ds,
\end{align*}
where the upper bound does not depend on $\kappa_1$. Therefore, 
\begin{align}
\lim\limits_{\epsilon\rightarrow 0}\sup\limits_{\kappa\in \left(0, \delta^*-\epsilon\right)}\left|L_1\left(\delta^*-\epsilon, \kappa\right)-L_1\left(\delta^*, \kappa\right)\right|=0.\label{left_hand_side_result}
\end{align}
By \eqref{L_1_uniform_continuous} and \eqref{left_hand_side_result}, we know that 
\begin{align*}
\sup\limits_{i\in \{0,...,M-1\}}\left|L_1\left(\delta^*-\epsilon, \left(i+1\right)\delta_M\right)+L_2\left(i\delta_M\right)-L_1\left(\delta^*, i\delta_M\right)-L_2\left(i\delta_M\right)\right|
\end{align*}
can be arbitrarily small for sufficiently small $\epsilon$ and sufficiently large $M$. We then analyze the right hand side of \eqref{final_inequality_time_period}. It is easy to see that
\begin{align*}
\left|L_1\left(\delta^*-\epsilon, \kappa^*-\frac{\sigma}{N}\right)+L_2\left(\kappa^*-\frac{\sigma}{N}\right)-L_1\left(\delta^*, \kappa^*\right)-L_2\left(\kappa^*\right)\right|   
\end{align*}
can be arbitrarily small for sufficiently small $\epsilon$ and sufficiently large $N$, which completes the first step which shows that the probability of the following event,
$$
\Biggl\{\frac{B_{n,y,\left(t_1, t_2\right)}}{\bar
{B}_{n,y,\left(t_1, t_2\right)}}\in \left(1+\delta^*-\epsilon, 1+\delta^*+\epsilon\right), \frac{I_{n, \left(t_1, t_2\right)}}{\bar{I}_{n, \left(t_1, t_2\right)}}\in \left(1, 1+\kappa^*-\sigma\right)\Biggr\},
$$
decays faster than that of the event of interest (the probability of which is evaluated in \eqref{def_main_event}.

By a similar argument, we can deal with the following event (the second step)
$$
\Biggl\{\frac{B_{n,y,\left(t_1, t_2\right)}}{\bar
{B}_{n,y,\left(t_1, t_2\right)}}\in \left(1+\delta^*-\epsilon, 1+\delta^*+\epsilon\right), \frac{I_{n, \left(t_1, t_2\right)}}{\bar{I}_{n, \left(t_1, t_2\right)}}\in \left(1+\kappa^*+\sigma, 1+\delta^*-\epsilon\right)\Biggr\},
$$
which completes the proof.
\qed
\end{proof}
\\


\subsubsection{Proof of Proposition \ref{period_clones_result}}\label{proposition_period_clones_result}

\begin{proof}\\
We first observe that
{\footnotesize
\begin{align*}
& \PP_{A_{n,y}}\left(\frac{I_{n,\left(t_1, t_2\right)}}{\bar{I}_{n,\left(t_1, t_2\right)}}\notin \left(1+\kappa^*-\sigma, 1+\kappa^*+\sigma\right)\right)\\
= & \PP_{A_{n,y}}\left(\frac{I_{n,\left(t_1, t_2\right)}}{\bar{I}_{n,\left(t_1, t_2\right)}}\notin \left(1+\kappa^*-\sigma, 1+\kappa^*+\sigma\right), {A_{n,y,\epsilon}^{t_1,t_2}}^c\right)\\
& \quad + \PP_{A_{n,y}}\left(\frac{I_{n,\left(t_1, t_2\right)}}{\bar{I}_{n,\left(t_1, t_2\right)}}\notin \left(1+\kappa^*-\sigma, 1+\kappa^*+\sigma\right), A_{n,y,\epsilon}^{t_1,t_2}\right).
\end{align*}}
By Lemma \ref{period_mutants_concentration}, we know that the first probability decays exponentially fast, and thus we focus on the second probability. We have
\begin{align*}
& \PP_{A_{n,y}}\left(\frac{I_{n,\left(t_1, t_2\right)}}{\bar{I}_{n,\left(t_1, t_2\right)}}\notin \left(1+\kappa^*-\sigma, 1+\kappa^*+\sigma\right), A_{n,y,\epsilon}^{t_1,t_2}\right)\\
= & \frac{\PP\left(\frac{I_{n,\left(t_1, t_2\right)}}{\bar{I}_{n,\left(t_1, t_2\right)}}\notin \left(1+\kappa^*-\sigma, 1+\kappa^*+\sigma\right),A_{n,y,\epsilon}^{t_1,t_2}\right)}{\PP\left(A_{n,y}\right)}\\
& \quad \quad \quad \quad \times \PP\left(A_{n,y}\middle| \frac{I_{n,\left(t_1, t_2\right)}}{\bar{I}_{n,\left(t_1, t_2\right)}}\notin \left(1+\kappa^*-\sigma, 1+\kappa^*+\sigma\right), A_{n,y,\epsilon}^{t_1,t_2}\right).
\end{align*}
Because
\begin{align*}
\PP\left(A_{n,y}\middle| \frac{I_{n,\left(t_1, t_2\right)}}{\bar{I}_{n,\left(t_1, t_2\right)}}\notin \left(1+\kappa^*-\sigma, 1+\kappa^*+\sigma\right), A_{n,y,\epsilon}^{t_1,t_2}\right)\le \PP\left(A_{n,y}\middle| B_{n,y,\left(t_1, t_2\right)}= \lceil \left(1+\delta^*+\epsilon\right)\bar
{B}_{n,y,\left(t_1, t_2\right)} \rceil \right),
\end{align*}
and
\begin{align*}
\PP\left(A_{n,y}\middle| \frac{I_{n,\left(t_1, t_2\right)}}{\bar{I}_{n,\left(t_1, t_2\right)}}\notin \left(1+\kappa^*-\sigma, 1+\kappa^*+\sigma\right), A_{n,y,\epsilon}^{t_1,t_2}\right)\ge \PP\left(A_{n,y}\middle| B_{n,y,\left(t_1, t_2\right)}= \lfloor \left(1+\delta^*-\epsilon\right)\bar
{B}_{n,y,\left(t_1, t_2\right)} \rfloor \right),
\end{align*}
by Theorem \ref{recurrence_LD}, we have
\begin{align*}
\left| \limsup\limits_{n\rightarrow \infty}\frac{1}{n^{1-\alpha}}\PP\left(A_{n,y}\middle| \frac{I_{n,\left(t_1, t_2\right)}}{\bar{I}_{n,\left(t_1, t_2\right)}}\notin \left(1+\kappa^*-\sigma, 1+\kappa^*+\sigma\right), A_{n,y,\epsilon}^{t_1,t_2}\right)- \limsup\limits_{n\rightarrow \infty}\frac{1}{n^{1-\alpha}}\PP\left(A_{n,y}\middle| A_{n,y,\epsilon}^{t_1,t_2}\right)\right| \rightarrow 0
\end{align*}
as $\epsilon$ goes to zero. From the proof of Lemma \ref{period_clones_concentration}, we know that for sufficiently small $\epsilon$, the large deviations rate of 
\begin{align*}
\PP\left(\frac{I_{n,\left(t_1, t_2\right)}}{\bar{I}_{n,\left(t_1, t_2\right)}}\notin \left(1+\kappa^*-\sigma, 1+\kappa^*+\sigma\right),A_{n,y,\epsilon}^{t_1,t_2}\right)
\end{align*}
is larger than that of 
\begin{align*}
\PP\left(\frac{I_{n,\left(t_1, t_2\right)}}{\bar{I}_{n,\left(t_1, t_2\right)}}\in \left(1+\kappa^*-\sigma, 1+\kappa^*+\sigma\right),A_{n,y,\epsilon}^{t_1,t_2}\right),
\end{align*}
and their difference is bounded away from zero. The desired result then follows from the fact that
\begin{align*}
\frac{\PP\left(\frac{I_{n,\left(t_1, t_2\right)}}{\bar{I}_{n,\left(t_1, t_2\right)}}\notin \left(1+\kappa^*-\sigma, 1+\kappa^*+\sigma\right),A_{n,y,\epsilon}^{t_1,t_2}\right)}{\PP\left(A_{n,y,\epsilon}^{t_1,t_2}\right)} \frac{\PP\left(A_{n,y,\epsilon}^{t_1,t_2}\right)}{\PP\left(A_{n,y}\right)}  \PP\left(A_{n,y}\middle| A_{n,y,\epsilon}^{t_1,t_2}\right)
\end{align*}
decays exponentially fast as the first term decays exponentially fast by Lemma \ref{period_clones_concentration} and $\frac{\PP\left(A_{n,y,\epsilon}^{t_1,t_2}\right)}{\PP\left(A_{n,y}\right)}\PP\left(A_{n,y}\middle| A_{n,y,\epsilon}^{t_1,t_2}\right)\approx 1$ by Lemma \ref{period_mutants_concentration}.  
\qed
\end{proof}
\\


\section{Proof of results in Section \ref{sec:simpson_with_conditioning}}\label{Sec:proof_simpson_with_conditioning}


\subsection{Proof of Theorem \ref{Simpson_conditioning_I}}\label{theorem_Simpson_conditioning_I}
\begin{proof}\\
The key idea of the proof is to carefully choose a sequence of the most likely events $E_{n,y,\epsilon_1} \subset A_{n,y}$, such that 
$$
\frac{\PP\left(A_{n,y} \setminus E_{n,y,\epsilon_1}\right)}{\PP\left(A_{n,y}\right)}
$$ 
decays exponentially fast. Then it suffices to calculate the Simpson's Index conditioned on $E_{n,y,\epsilon_1}$. Our choice of $E_{n,y,\epsilon_1}$ will make the analysis much easier. 

Recall that we denote by $I_n\left(\zeta_n-y\right)$ the number of clones generated in the time period $\left(0, \zeta_n-y \right)$. In \cite{BP2020} (see Section 3.2), we obtain the most likely number of clones given early recurrence has occurred. We restate the result in the following lemma.
\begin{lemma}\label{most_likely_num}
{\small
\begin{align*}
& \quad \mbox{argmax}_{a\in (0,e^{\lambda_1 y}]}\lim_{n\rightarrow \infty}\frac{1}{n^{1-\alpha}}\log \PP\left(\gamma_n\le \zeta_n-y \middle| I_n\left(\zeta_n-y\right)=\lfloor -a\frac{\mu}{\lambda_0}n^{1-\alpha} \rfloor\right)\PP\left(I_n\left(\zeta_n-y\right)=\lfloor -a\frac{\mu}{\lambda_0}n^{1-\alpha} \rfloor\right)\\
& = -\lambda_0\int_{0}^{\infty}\frac{e^{\lambda_1 s}}{e^{\lambda_1 s}-\theta^*_y}e^{\lambda_0 s}ds,
\end{align*}}
where $\theta^*_y$ is defined in \eqref{theta_2}.
\end{lemma}
From Lemma \ref{most_likely_num}, we know that given early recurrence has happened, the number of clones is approximately $\mu n^{1-\alpha}\int_{0}^{\infty}\frac{e^{\lambda_1 s}}{e^{\lambda_1 s}- \theta^*_y}e^{\lambda_0 s}ds$. Hence we consider the following set of events for $\epsilon_1>0$:
\begin{align*}
E_{n,y,\epsilon_1}=\{\gamma_n < \zeta_n-y, I_n\left(\zeta_n-y\right)\in O_{n,y,\epsilon_1}\},
\end{align*}
where
\begin{align*}
O_{n,y,\epsilon_1}=\left(\left(1-\epsilon_1\right)\mu n^{1-\alpha}\int_{0}^{\infty}\frac{e^{\lambda_1 s}}{e^{\lambda_1 s}- \theta^*_y}e^{\lambda_0 s}ds, \left(1+\epsilon_1\right)\mu n^{1-\alpha}\int_{0}^{\infty}\frac{e^{\lambda_1 s}}{ e^{\lambda_1 s}-\theta^*_y}e^{\lambda_0 s}ds\right).
\end{align*}
We choose $E_{n,y,\epsilon_1}$ in this way so that the number of clones is concentrated. For simplicity, we let $\PP_{E_{n,y,\epsilon_1}}\left(\cdot\right)=\PP\left(\cdot \middle| E_{n,y,\epsilon_1}\right)$. We then analyze $\PP\left(A_{n,y}\setminus E_{n,y,\epsilon_1}\right)$ in the next lemma.
\begin{lemma}\label{dominated_event}
\begin{align*}
\limsup\limits_{n\rightarrow \infty}\frac{1}{n^{1-\alpha}}\log\PP\left(A_{n,y}\setminus E_{n,y,\epsilon_1}\right) < -L\left(y\right),
\end{align*}
where
\begin{align*}
L\left(y\right)& \triangleq - \lim\limits_{n\rightarrow \infty}\frac{1}{n^{1-\alpha}}\log\PP\left(A_{n,y}\right) =\sup\limits_{\theta\in\left(0,1\right)}\left[\frac{\theta \mu e^{\lambda_1 y}}{\lambda_1-\lambda_0}-\mu\theta\int_{0}^{\infty}\frac{e^{\lambda_0 s}}{e^{\lambda_1s}-\theta}ds\right].
\end{align*}
\end{lemma}
\begin{proof}
See Section \ref{lemma_dominated_event}.
\end{proof}
\\

By Lemma \ref{dominated_event}, $\frac{\PP\left(A_{n,y} \setminus E_{n,y,\epsilon_1}\right)}{\PP\left(A_{n,y}\right)}$ decays exponentially fast as desired. Since Simpson's Index is at most $1$, we have
\begin{align*}
\lim\limits_{n\rightarrow \infty} \left|n^{1-\alpha}\EE_{A_{n,y}}\left[R_{n,y}\right]-  n^{1-\alpha}\EE_{E_{n,y,\epsilon_1}}\left[R_{n,y}\right]\right|=0.    
\end{align*}

Notice that in the time period $\left(0, \zeta_n-y\right)$, each clone is generated according to a Poisson process. We denote by $X_{n,y,i}$ the size of the i-th clone at $\zeta_n-y$. Note that the mutant clones are ordered at random, not in chronological order by when the mutation occurred. For ease of exposition, we drop the subscript $n$ and $y$. Since the generation times of clones as well as their evolution are independent conditioned on the total number of clones generated in the time period $\left(0, \zeta_n-y\right)$, we can obtain that
\begin{align}\label{mean_clone_size}
\EE\left[X_i\right]=e^{\lambda_1 \left(\zeta_n-y\right)}\frac{\int_0^{\zeta_n-y}e^{-\left(\lambda_1-\lambda_0\right)t}dt}{\int_0^{\zeta_n-y}e^{\lambda_0 t}dt}.
\end{align} 
Hence, it is convenient to condition on the number of clones and then apply the law of total expectation:
\begin{align*}
n^{1-\alpha}\EE_{E_{n,y,\epsilon_1}}\left[R_{n,y}\right]=n^{1-\alpha}\sum_{k\in O_{n,y,\epsilon_1}}\EE\left[R_{n,y}\middle| \gamma_n < \zeta_n-y, I_n\left(\zeta_n-y\right)=k\right]\frac{\PP\left(\gamma_n < \zeta_n-y, I_n\left(\zeta_n-y\right)=k\right)}{\PP\left(E_{n,y,\epsilon_1}\right)}.
\end{align*}
Define $\tilde{R}_{n,y}=\frac{\sum_{i=1}^{I_n\left(\zeta_n-y\right)}X_i^2}{n^2}$. 
Because $Z_1^n\left(\zeta_n-y\right)$ concentrates around $n$ conditioned on $A_{n,y}$ (this can be shown easily with Theorem \ref{recurrence_LD} and the assumption that $d_1=0$)\footnote{For any $\epsilon>0$, $\PP\left(Z_1^n\left(\zeta_n-y\right)>\left(1+\epsilon\right)n \mid A_{n,y}\right)$ decays exponentially fast by Theorem \ref{recurrence_LD} and $\PP\left(Z_1^n\left(\zeta_n-y\right)<n \mid A_{n,y}\right) = 0$ by the assumption that $d_1=0$.}, we have
\begin{align*}
\lim\limits_{n\rightarrow \infty}\left|n^{1-\alpha}\EE\left[R_{n,y}\middle| \gamma_n < \zeta_n-y\right]-n^{1-\alpha}\EE\left[\tilde{R}_{n,y}\middle| \gamma_n < \zeta_n-y\right]\right|=0.
\end{align*}
Hence, it suffices to analyze
\begin{align}
& n^{1-\alpha}\EE\left[\tilde{R}_{n,y}\middle| \gamma_n < \zeta_n-y, I_n\left(\zeta_n-y\right)=k\right] \nonumber\\ 
= &\frac{k}{n^{1+\alpha}}\EE\left[X_i^2\middle| \sum_{i=1}^{k}X_i\ge n\right] \nonumber\\
= &\frac{ke^{2\lambda_1\left(\zeta_n-y\right)}}{n^{1+\alpha}}\EE\left[\left(e^{-\lambda_1\left(\zeta_n-y\right)}X_i\right)^2\middle| \sum_{i=1}^{k}e^{-\lambda_1\left(\zeta_n-y\right)}X_i\ge e^{-\lambda_1\left(\zeta_n-y\right)}n\right].\label{eqn:tildeR_value}
\end{align}
The next step is to find a bound on \begin{align}\label{conditioned_second_moment}
\EE\left[\left(e^{-\lambda_1\left(\zeta_n-y\right)}X_i\right)^2\middle| \sum_{i=1}^{k}e^{-\lambda_1\left(\zeta_n-y\right)}X_i\ge e^{-\lambda_1\left(\zeta_n-y\right)}n\right],
\end{align} 
for $k\in O_{n,y,\epsilon_1}$. We note that $\sum\limits_{i=1}^{k}e^{-\lambda_1\left(\zeta_n-y\right)}X_i\ge e^{-\lambda_1\left(\zeta_n-y\right)}n$ is equivalent to $\frac{\sum_{i=1}^{k}e^{-\lambda_1\left(\zeta_n-y\right)}X_i}{k}\ge q_{k,n} \EE\left[e^{-\lambda_1\left(\zeta_n-y\right)}X_i\right]$, where
\begin{align*}
q_{k,n}=\frac{n}{k\EE\left[X_i\right]}.
\end{align*}
Since
\begin{align*}
\frac{\lambda_1-\lambda_0}{\mu}n^{\alpha}\le e^{\lambda_1 \zeta_n}\le \frac{\lambda_1-\lambda_0}{\mu}n^{\alpha}+1,
\end{align*}
we have
\begin{align*}
\lim\limits_{n\rightarrow \infty}\frac{n/\left(\left(1-\epsilon_1\right)\mu n^{1-\alpha} \int_{0}^{\infty}\frac{e^{\lambda_1 s}}{e^{\lambda_1 s}- \theta^*_y}e^{\lambda_0 s}ds\right)}{\EE\left[X_i\right]}=- \frac{e^{\lambda_1 y}}{\left(1-\epsilon_1\right) \lambda_0 \int_{0}^{\infty}\frac{e^{\lambda_1 s}}{e^{\lambda_1 s}- \theta^*_y}e^{\lambda_0 s}ds},
\end{align*}
and
\begin{align*}
\lim\limits_{n\rightarrow \infty}\frac{n/\left(\left(1+\epsilon_1\right)\mu n^{1-\alpha} \int_{0}^{\infty}\frac{e^{\lambda_1 s}}{e^{\lambda_1 s}- \theta^*_y}e^{\lambda_0 s}ds\right)}{\EE\left[X_i\right]}= - \frac{e^{\lambda_1 y}}{\left(1+\epsilon_1\right) \lambda_0 \int_{0}^{\infty}\frac{e^{\lambda_1 s}}{e^{\lambda_1 s}- \theta^*_y}e^{\lambda_0 s}ds}.
\end{align*}
Hence, for arbitrarily small but fixed $\epsilon_2>0$, there exists $N_{\epsilon_2}>0$ such that when $n>N_{\epsilon_2}$, $k\in O_{n,y,\epsilon_1}$ implies that
\begin{align}
q_{k,n}\in \tilde{O}_{y,\epsilon_1, \epsilon_2}=\left[\frac{\left(1-\epsilon_2\right)}{\left(1+\epsilon_1\right)}q^*, \frac{\left(1+\epsilon_2\right)}{\left(1-\epsilon_1\right)}q^*\right]. \label{bound:q_k_n}
\end{align}
where 
\begin{align}
q^*=-\frac{e^{\lambda_1 y}}{\lambda_0\int_{0}^{\infty}\frac{e^{\lambda_1 s}}{e^{\lambda_1 s}- \theta^*_y}e^{\lambda_0 s}ds}.\label{def:q^star}
\end{align}
Since $q^*>1$ by the definition of $\theta^*_y$ (see Lemma 1 in \cite{BP2020}), we will make $\epsilon_1$, $\epsilon_2$ sufficiently small such that $\frac{\left(1-\epsilon_2\right)}{\left(1+\epsilon_1\right)}q^*>1$ in the rest of the proof. Let
\begin{align}
Q_{k,n} & \overset{\Delta}{=} q_{k,n}\EE\left[e^{-\lambda_1 \left(\zeta_n-y\right)}X_i\right]\nonumber\\
& =q_{k,n}\frac{\int_0^{\zeta_n-y}e^{-\left(\lambda_1-\lambda_0\right)t}dt}{\int_0^{\zeta_n-y}e^{\lambda_0 t}dt},\label{def:Q_k_n}
\end{align}
where the equality follows from \eqref{mean_clone_size}. Note that $Q_{k,n}$ depends on $y$, but we omit this dependence for simplicity. Define
\begin{align*}
Q^*& \overset{\Delta}{=} \lim\limits_{n\rightarrow \infty}q^*\frac{\int_0^{\zeta_n-y}e^{-\left(\lambda_1-\lambda_0\right)t}dt}{\int_0^{\zeta_n-y}e^{\lambda_0 t}dt}=-\frac{\lambda_0}{\lambda_1-\lambda_0}q^*.
\end{align*}
We then calculate the log moment-generating function for $e^{-\lambda_1 \left(\zeta_n-y\right)}X_i$ and its derivatives. For $\theta \in \left(0,1\right)$, let
\begin{align*}
\Lambda_n\left(\theta\right) &\triangleq \log\EE\left[e^{\theta e^{-\lambda_1 \left(\zeta_n-y\right)}X_i}\right]\\
& = \log\left(\frac{\int_0^{\zeta_n-y}e^{\lambda_0 t}\frac{e^{\lambda_1 t}e^{-\lambda_1\left(\zeta_n-y\right)}e^{\theta e^{-\lambda_1\left(\zeta_n-y\right)}}}{e^{\lambda_1 t}e^{-\lambda_1\left(\zeta_n-y\right)}e^{\theta e^{-\lambda_1\left(\zeta_n-y\right)}}-e^{\theta e^{-\lambda_1\left(\zeta_n-y\right)}}+1} dt}{\int_0^{\zeta_n-y}e^{\lambda_0 t}dt}\right)\\
& \rightarrow \log\left(-\lambda_0 \int_0^{\infty}e^{\lambda_0 t}\frac{1}{1-\theta e^{-\lambda_1 t}}dt\right) \text{ as } n\rightarrow \infty.
\end{align*}
We denote the limit by $\Lambda\left(\theta\right)$. Taking the first derivative, we have
\begin{align*}
\Lambda_n'\left(\theta\right)& =\frac{\EE\left[e^{-\lambda_1 \left(\zeta_n-y\right)}X_i e^{\theta e^{-\lambda_1 \left(\zeta_n-y\right)}X_i}\right]}{\EE\left[e^{\theta e^{-\lambda_1 \left(\zeta_n-y\right)}X_i}\right]}\\
& =\frac{\int_0^{\zeta_n-y}e^{\lambda_0 t}\frac{e^{\lambda_1 t}e^{-2\lambda_1\left(\zeta_n-y\right)}e^{-\theta e^{-\lambda_1\left(\zeta_n-y\right)}}}{\left(e^{\lambda_1 t}e^{-\lambda_1\left(\zeta_n-y\right)}-1+e^{-\theta e^{-\lambda_1\left(\zeta_n-y\right)}}\right)^2}dt}{\int_0^{\zeta_n-y}e^{\lambda_0 t}\frac{e^{\lambda_1 t}e^{-\lambda_1\left(\zeta_n-y\right)}}{e^{\lambda_1 t}e^{-\lambda_1\left(\zeta_n-y\right)}-1+e^{-\theta e^{-\lambda_1\left(\zeta_n-y\right)}}} dt}\\
& \rightarrow \frac{\int_0^{\infty}e^{\lambda_0 t}\frac{e^{-\lambda_1t}}{\left(1-\theta e^{-\lambda_1 t}\right)^2}dt}{\int_0^{\infty}e^{\lambda_0 t}\frac{1}{1-\theta e^{-\lambda_1 t}}dt} \text{ as } n\rightarrow \infty.
\end{align*}
We denote the limit by $\Lambda'\left(\theta\right)$. Take the second derivative, we have
\begin{align*}
\Lambda_n''\left(\theta\right)& =\frac{\EE\left[e^{-2\lambda_1 \left(\zeta_n-y\right)}X_i^2 e^{\theta e^{-\lambda_1 \left(\zeta_n-y\right)}X_i}\right]\EE\left[e^{\theta e^{-\lambda_1 \left(\zeta_n-y\right)}X_i}\right]-\left(\EE\left[e^{-\lambda_1 \left(\zeta_n-y\right)}X_i e^{\theta e^{-\lambda_1 \left(\zeta_n-y\right)}X_i}\right]\right)^2}{\left(\EE\left[e^{\theta e^{-\lambda_1 \left(\zeta_n-y\right)}X_i}\right]\right)^2}\\
& > 0, \text{ and }
\end{align*}
\begin{align*}
\Lambda_n''\left(\theta\right)&=\frac{\int_0^{\zeta_n-y}e^{\lambda_0 t}\frac{e^{\lambda_1 t}e^{-3\lambda_1\left(\zeta_n-y\right)}e^{-\theta e^{-\lambda_1\left(\zeta_n-y\right)}}\left(e^{-\theta e^{-\lambda_1\left(\zeta_n-y\right)}}-e^{\lambda_1 t}e^{-\lambda_1\left(\zeta_n-y\right)}+1\right)}{\left(e^{\lambda_1 t}e^{-\lambda_1\left(\zeta_n-y\right)}-1+e^{-\theta e^{-\lambda_1\left(\zeta_n-y\right)}}\right)^3}}{\int_0^{\zeta_n-y}e^{\lambda_0 t}\frac{e^{\lambda_1 t}e^{-\lambda_1\left(\zeta_n-y\right)}}{e^{\lambda_1 t}e^{-\lambda_1\left(\zeta_n-y\right)}-1+e^{-\theta e^{-\lambda_1\left(\zeta_n-y\right)}}} dt}-\left(\Lambda_n'\left(\theta\right)\right)^2\\
& \rightarrow \frac{\int_0^{\infty}e^{\lambda_0 t}\frac{2e^{-2\lambda_1 t}}{\left(1-\theta e^{-\lambda_1 t}\right)^3}dt\int_0^{\infty}e^{\lambda_0 t}\frac{1}{1-\theta e^{-\lambda_1 t}}dt-\left(\int_0^{\infty}e^{\lambda_0 t}\frac{e^{-\lambda_1t}}{\left(1-\theta e^{-\lambda_1 t}\right)^2}dt\right)^2}{\left(\int_0^{\infty}e^{\lambda_0 t}\frac{1}{1-\theta e^{-\lambda_1 t}}dt\right)^2} \text{ as } n\rightarrow \infty.
\end{align*}
We denote the limit by $\Lambda''\left(\theta\right)$. Taking the third derivative, we have
\begin{align*}
\Lambda_n'''\left(\theta\right)=\frac{f_4\left(\theta\right)f_1^2\left(\theta\right)-3f_3\left(\theta\right)f_2\left(\theta\right)f_1\left(\theta\right)+2f_2^3\left(\theta\right)}{f_1^3\left(\theta\right)},
\end{align*}
where
\begin{align}\label{f_1}
f_1\left(\theta\right)=\int_0^{\zeta_n-y}e^{\lambda_0 t}\frac{e^{\lambda_1 t}e^{-\lambda_1\left(\zeta_n-y\right)}}{e^{\lambda_1 t}e^{-\lambda_1\left(\zeta_n-y\right)}-1+e^{-\theta e^{-\lambda_1\left(\zeta_n-y\right)}}} dt,
\end{align}
\begin{align}\label{f_2}
f_2\left(\theta\right)=\int_0^{\zeta_n-y}e^{\lambda_0 t}\frac{e^{\lambda_1 t}e^{-2\lambda_1\left(\zeta_n-y\right)}e^{-\theta e^{-\lambda_1\left(\zeta_n-y\right)}}}{\left(e^{\lambda_1 t}e^{-\lambda_1\left(\zeta_n-y\right)}-1+e^{-\theta e^{-\lambda_1\left(\zeta_n-y\right)}}\right)^2}dt, 
\end{align}
\begin{align}\label{f_3}
f_3\left(\theta\right)=\int_0^{\zeta_n-y}e^{\lambda_0 t}\frac{e^{\lambda_1 t}e^{-3\lambda_1\left(\zeta_n-y\right)}e^{-\theta e^{-\lambda_1\left(\zeta_n-y\right)}}\left(e^{-\theta e^{-\lambda_1\left(\zeta_n-y\right)}}-e^{\lambda_1 t}e^{-\lambda_1\left(\zeta_n-y\right)}+1\right)}{\left(e^{\lambda_1 t}e^{-\lambda_1\left(\zeta_n-y\right)}-1+e^{-\theta e^{-\lambda_1\left(\zeta_n-y\right)}}\right)^3}dt,
\end{align}
\begin{align}\label{f_4}
f_4\left(\theta\right)=\int_0^{\zeta_n-y}e^{\lambda_0 t}\frac{e^{\lambda_1 t}e^{-4\lambda_1\left(\zeta_n-y\right)}e^{-\theta e^{-\lambda_1\left(\zeta_n-y\right)}}g_4\left(\theta\right)}{\left(e^{\lambda_1 t}e^{-\lambda_1\left(\zeta_n-y\right)}-1+e^{-\theta e^{-\lambda_1\left(\zeta_n-y\right)}}\right)^4}dt, \text{ and }
\end{align}
\begin{align*}
g_4\left(\theta\right)& =4e^{-\theta e^{-\lambda_1\left(\zeta_n-y\right)}}-4e^{\lambda_1 t}e^{-\theta e^{-\lambda_1\left(\zeta_n-y\right)}}e^{-\lambda_1\left(\zeta_n-y\right)}+e^{-2\theta e^{-\lambda_1\left(\zeta_n-y\right)}}\\
&\quad \quad \quad -2e^{\lambda_1 t}e^{-\lambda_1\left(\zeta_n-y\right)}+e^{2\lambda_1 t}e^{-2\lambda_1\left(\zeta_n-y\right)}+1.
\end{align*}
By observing that $f_1\left(\theta\right)$ through $f_4\left(\theta\right)$ have well defined positive limits as $n$ approaches infinity, for any $0<\epsilon<1$, we can find an $M_{\epsilon}>0$ such that for sufficiently large $n$, 
\begin{align}\label{bound_third}
\left|\Lambda_n'''\left(\theta\right)\right|\le M_{\epsilon}  
\end{align}
for $\theta\in \left(0, 1-\epsilon\right)$. We then show the following lemma.
\begin{lemma}\label{unique_sol}
For sufficiently large $n$ and sufficiently small $\epsilon_1$ and $\epsilon_2$, we can find a unique solution $\eta_{k,n}\in \left(0,1\right)$ to
\begin{align*}
\Lambda_n'\left(\eta\right)=Q_{k,n}.
\end{align*}
\end{lemma}
\begin{proof}
See Section \ref{lemma_unique_sol}.
\end{proof}

We then show that $\theta^*_y$ is the solution to the equation in Lemma \ref{unique_sol} when we take limits on both sides.
\begin{lemma}\label{limit_eta_q}
$\Lambda'\left(\theta^*_y \right)=Q^*$.
\end{lemma}
\begin{proof}
See Section \ref{lemma_limit_eta_q}.
\end{proof}

By Lemma \ref{limit_eta_q} and the previous calculation of derivatives of $\Lambda_n\left(\theta\right)$, we can obtain very tight bounds for $\eta_{k,n}$, $Q_{k,n}$, $\Lambda_n\left(\eta_{k,n}\right)$, $\Lambda_n'\left(\eta_{k,n}\right)$ and $\Lambda_n''\left(\eta_{k,n}\right)$ for sufficiently large $n$, where $\eta_{k,n}$ is the solution to $\Lambda_n'\left(\eta\right)=Q_{k,n}$.
\begin{lemma}\label{uniform_bounds_parameters_1}
Define 
$$
Q_{k,n}\left[\epsilon_1\right]=\left(\inf\limits_{k\in O_{n,y,\epsilon_1}}Q_{k,n}, \sup\limits_{k\in O_{n,y,\epsilon_1}}Q_{k,n}\right),
$$
and define $\eta_{k,n}\left[\epsilon_1\right]$, $\Lambda_n\left(\eta_{k,n}\right)\left[\epsilon_1\right]$, $\Lambda_n'\left(\eta_{k,n}\right)\left[\epsilon_1\right]$, $\Lambda_n''\left(\eta_{k,n}\right)\left[\epsilon_1\right]$ in the same way. For any $\epsilon>0$, there exists $\delta_{\epsilon}>0$ and $N_{\epsilon}>0$ such that when $\epsilon_1<\delta_{\epsilon}$ and $n>N_{\epsilon}$,
\begin{align*}
& Q_{k,n}\left[\epsilon_1\right]\subset \left(Q^*-\epsilon, Q^*+\epsilon\right)\\
& \eta_{k,n}\left[\epsilon_1\right]\subset\left(\theta^*_y-\epsilon, \theta^*_y+\epsilon\right)\\
& \Lambda_n\left(\eta_{k,n}\right)\left[\epsilon_1\right]\subset\left(\Lambda\left(\theta^*_y\right)-\epsilon, \Lambda\left(\theta^*_y\right)+\epsilon\right)\\
& \Lambda_n'\left(\eta_{k,n}\right)\left[\epsilon_1\right]\subset\left(\Lambda'\left(\theta^*_y\right)-\epsilon, \Lambda'\left(\theta^*_y\right)+\epsilon\right)\\
& \Lambda_n''\left(\eta_{k,n}\right)\left[\epsilon_1\right]\subset\left(\Lambda''\left(\theta^*_y\right)-\epsilon, \Lambda''\left(\theta^*_y\right)+\epsilon\right).
\end{align*}
\end{lemma}
\begin{proof}
See Section \ref{lemma_uniform_bounds_parameters_1}.
\end{proof}

We then show the following proposition which gives a bound on \eqref{conditioned_second_moment}.
\begin{proposition}\label{uniform_bound}
Let $\nu_{n,y}$ be the law of $e^{-\lambda_1\left(\zeta_n-y\right)}X_i$, and consider the probability measure $\tilde{\nu}_{n,y,k}$ defined by $d\tilde{\nu}_{n,y,k}/d\nu_{n,y} \left(x\right)=e^{\eta_{k,n} x - \Lambda_n\left(\eta_{k,n}\right)}$. For any $\epsilon>0$, there exists $N_{\epsilon}>0$ and $\delta_{\epsilon}>0$ such that when $n>N_{\epsilon}$, $\epsilon_1<\delta_{\epsilon}$, and $k\in O_{n,y,\epsilon_1}$,
{\small
\begin{align*}
& \quad \EE\left[\left(e^{-\lambda_1\left(\zeta_n-y\right)}X_i\right)^2\middle| \frac{\sum_{i=1}^{k}e^{-\lambda_1\left(\zeta_n-y\right)}X_i}{k}\ge Q_{k,n}\right]\in \left(\int_0^{\infty}x^2 d\tilde{\nu}_{n,y,k}\left(x\right)-\epsilon,\int_0^{\infty}x^2 d\tilde{\nu}_{n,y,k}\left(x\right)+\epsilon\right).
\end{align*}}
\end{proposition}
\begin{proof}
See Section \ref{proposition_uniform_bound}.
\end{proof}
\\

Notice that (see page 111 of \cite{LD})
\begin{align*}
\int_0^{\infty}x^2 d\tilde{\nu}_{n,y,k}\left(x\right)&=\Lambda_n''\left(\eta_{k,n}\right)+Q_{k,n}^2.
\end{align*}
By Lemma \ref{uniform_bounds_parameters_1}, we know that for any $\epsilon>0$, the following result holds for sufficiently small $\epsilon_1$, sufficiently large $n$ and $k\in O_{n,y,\epsilon_1}$:
\begin{align*}
\int_0^{\infty}x^2 d\tilde{\nu}_{n,y,k}\in \left(\Lambda''\left(\theta^*_y\right)+{Q^*}^2-\epsilon, \Lambda''\left(\theta^*_y\right)+{Q^*}^2+\epsilon\right).
\end{align*}
Hence
{\footnotesize
\begin{align*}
& \quad \lim\limits_{n\rightarrow \infty} n^{1-\alpha}\EE_{A_{n,y}}\left[R_{n,y}\right]\\
& =\lim\limits_{n\rightarrow \infty}n^{1-\alpha}\EE_{E_{n,y,\epsilon_1}}\left[R_{n,y}\right]\\
& = \lim\limits_{n\rightarrow \infty}n^{1-\alpha}\sum_{k\in O_{n,y,\epsilon_1}}\EE\left[R_{n,y}\middle| \gamma_n < \zeta_n-y, I_n\left(\zeta_n-y\right)=k\right]\frac{\PP\left(\gamma_n < \zeta_n-y, I_n\left(\zeta_n-y\right)=k\right)}{\PP\left(E_{n,y,\epsilon_1}\right)}\\
& =\lim\limits_{n\rightarrow \infty}n^{1-\alpha}\sum_{k\in O_{n,y,\epsilon_1}}\EE\left[\tilde{R}_{n,y}\middle| \gamma_n < \zeta_n-y, I_n\left(\zeta_n-y\right)=k\right]\frac{\PP\left(\gamma_n < \zeta_n-y, I_n\left(\zeta_n-y\right)=k\right)}{\PP\left(E_{n,y,\epsilon_1}\right)}\\
&=\lim\limits_{n\rightarrow \infty}n^{1-\alpha}\sum_{k\in O_{n,y,\epsilon_1}}\frac{ke^{2\lambda_1\left(\zeta_n-y\right)}}{n^{2}}\EE\left[\left(e^{-\lambda_1\left(\zeta_n-y\right)}X_i\right)^2\middle| \sum_{i=1}^{k}e^{-\lambda_1\left(\zeta_n-y\right)}X_i\ge e^{-\lambda_1\left(\zeta_n-y\right)}n\right]\frac{\PP\left(\gamma_n < \zeta_n-y, I_n\left(\zeta_n-y\right)=k\right)}{\PP\left(E_{n,y,\epsilon_1}\right)},
\end{align*}}
where the last equality follows from \eqref{eqn:tildeR_value}. By Proposition \ref{uniform_bound}, the definition of $O_{n,y,\epsilon_1}$ and $\tilde{O}_{y,\epsilon_1, \epsilon_2}$, and the fact that $\epsilon_1$ and $\epsilon_2$ can be arbitrarily small, we have
\begin{align*}
& \quad \lim\limits_{n\rightarrow \infty} n^{1-\alpha}\EE_{A_{n,y}}\left[R_{n,y}\right]\\
&=\lim\limits_{n\rightarrow \infty}n^{1-\alpha}\frac{\mu n^{1-\alpha}\int_0^{\infty}\frac{e^{\lambda_1 s}}{e^{\lambda_1 s}-\theta^*_y}e^{\lambda_0 s}ds}{n^2}e^{2\lambda_1\left(\zeta_n - y\right)}\left(\Lambda''\left(\theta^*_y\right)+{Q^*}^2\right)\\
&=\lim\limits_{n\rightarrow \infty}n^{1-\alpha}\frac{\mu n^{1-\alpha}\int_0^{\infty}\frac{e^{\lambda_1 s}}{e^{\lambda_1 s}-\theta^*_y}e^{\lambda_0 s}ds}{n^2}e^{2\lambda_1\left(\zeta_n - y\right)}\frac{2\int_0^{\zeta_n-y}\frac{e^{-\left(2\lambda_1 -\lambda_0\right)s}}{\left(1-\theta^*_y e^{-\lambda_1 s}\right)^3}ds}{\int_0^{\zeta_n-y}\frac{e^{\lambda_0 s}}{1-\theta^*_y e^{-\lambda_1 s}}ds}\\
& =\frac{2\left(\lambda_1-\lambda_0\right)^2}{\mu} e^{-2\lambda_1 y}\int_0^{\infty}\frac{e^{-\left(2\lambda_1-\lambda_0\right)s}}{\left(1-\theta^*_y e^{-\lambda_1 s}\right)^3}ds.
\end{align*}
\qed
\end{proof}

\subsubsection{Proof of Lemma \ref{dominated_event}}\label{lemma_dominated_event}


\begin{proof}\\
Let $z= -\lambda_0 \int_{0}^{\infty}\frac{e^{\lambda_1 s}}{e^{\lambda_1 s}-\theta^*_y}e^{\lambda_0 s}ds$, which is the optimal value obtained in Lemma \ref{most_likely_num}. We first notice that
\begin{align*}
A_{n,y}\setminus E_{n,y,\epsilon_1}=\{\gamma_n < \zeta_n-y, I_n\left(\zeta_n-y\right)\in O_1\}\cup \{\gamma_n < \zeta_n-y, I_n\left(\zeta_n-y\right)\in O_2\},
\end{align*}
where
\begin{align*}
O_1=\left(0,-\left(1-\epsilon_1\right)z\frac{\mu}{\lambda_0}n^{1-\alpha}\right), \text{ and } O_2=\left(-\left(1+\epsilon_1\right)z\frac{\mu}{\lambda_0}n^{1-\alpha}, \infty\right).
\end{align*}
We consider $O_2$ first. It is easy to show that the probability of $\{\gamma_n < \zeta_n-y, I_n\left(\zeta_n-y\right)\ge  -e^{\lambda_1 y}\frac{\mu}{\lambda_0}n^{1-\alpha}\}$ decays exponentially faster than that of $A_{n,y}$. Hence, we focus on the interval $\left(\left(1+\epsilon_1\right)z, e^{\lambda_1 y}\right)$, and divide it into $M$ sub-intervals with equal length. Let $\delta_M=\frac{e^{\lambda_1 y}-\left(1+\epsilon_1\right)z}{M}$ be the length of each sub-interval. Then for $i\in \{1,2,...,M\}$, we consider the sub-interval
\begin{align*}
O_{i,M}=\left(\left(1+\epsilon_1\right)z+\left(i-1\right)\delta_M, \left(1+\epsilon_1\right)z+i\delta_M\right).
\end{align*}
By stochastic dominance and the Gartner-Ellis Theorem \cite{Hollander}, we can show that
{\footnotesize
\begin{align*}
& \quad \limsup_{n\rightarrow \infty}\frac{1}{n^{1-\alpha}}\log \PP\left(\gamma_n\le \zeta_n-y, -\frac{\lambda_0}{\mu}n^{\alpha-1} I_n\left(\zeta_n-y\right)\in O_{i,M}\right)\\
& \le \limsup_{n\rightarrow \infty}\frac{1}{n^{1-\alpha}}\log \left(\PP\left(\gamma_n\le \zeta_n-y \middle| I_n\left(\zeta_n-y\right)= \lfloor -\frac{\mu}{\lambda_0}n^{1-\alpha}\left(\left(1+\epsilon_1\right)z+i\delta_M\right) \rfloor\right)\PP\left(-\frac{\lambda_0}{\mu}n^{\alpha-1} I_n\left(\zeta_n-y\right)\in O_{i,M}\right)\right)\\
& =-\sup\limits_{\theta\in \left(0,1\right)}\left[\frac{\theta \mu e^{\lambda_1 y}}{\lambda_1-\lambda_0}- \lfloor \left(1+\epsilon_1\right)z+i\delta_M \rfloor \frac{\mu}{\lambda_0}\log\left( -\lambda_0 \int_{0}^{\infty}\frac{e^{\lambda_1 s}}{e^{\lambda_1 s}-\theta}e^{\lambda_0 s}ds\right)\right]\\
& \quad \quad +\frac{\mu}{\lambda_0}\left(\lfloor \left(1+\epsilon_1\right)z+\left(i-1\right)\delta_M \rfloor\log\left(\lfloor \left(1+\epsilon_1\right)z+\left(i-1\right)\delta_M \rfloor\right)-\lfloor \left(1+\epsilon_1\right)z+\left(i-1\right)\delta_M \rfloor+1\right).
\end{align*}}
Let 
\begin{align*}
f_y\left(x\right)=\sup\limits_{\theta\in \left(0,1\right)}\left[\frac{\theta \mu e^{\lambda_1 y}}{\lambda_1-\lambda_0}+ x \frac{\mu}{\lambda_0}\log\left( -\lambda_0 \int_{0}^{\infty}\frac{e^{\lambda_1 s}}{e^{\lambda_1 s}-\theta}e^{\lambda_0 s}ds\right)\right]-\frac{\mu}{\lambda_0}\left(x\log\left(x\right)-x+1\right).
\end{align*}
From \cite{BP2020} and the Envelope Theorem, we know that $f_y\left(x\right)$ is strictly convex in $x$, and
\begin{align*}
& z=\mbox{argmin}_{x\in \left(0,e^{\lambda_1 y}\right]}f_y\left(x\right).
\end{align*}
Moreover, $\theta^*_y$ is the corresponding optimal solution to the inner optimization problem. By some algebra we can obtain that $f_y\left(z\right)=L\left(y\right)$. Hence, we conclude that
\begin{align*}
& L\left(y\right)-\sup\limits_{\theta\in \left(0,1\right)}\left[\frac{\theta \mu e^{\lambda_1 y}}{\lambda_1-\lambda_0}- \lfloor \left(1+\epsilon_1\right)z+i\delta_M \rfloor \frac{\mu}{\lambda_0}\log\left( -\lambda_0 \int_{0}^{\infty}\frac{e^{\lambda_1 s}}{e^{\lambda_1 s}-\theta}e^{\lambda_0 s}ds\right)\right]\\
& \quad \quad +\frac{\mu}{\lambda_0}\left(\lfloor \left(1+\epsilon_1\right)z+i\delta_M \rfloor\log\left(\lfloor \left(1+\epsilon_1\right)z+i\delta_M \rfloor\right)-\lfloor \left(1+\epsilon_1\right)z+i\delta_M \rfloor+1\right)
\end{align*}
is strictly less than $0$ for any $i$ and $M$. Therefore, there exists $\epsilon_2>0$ and an integer $M_{\epsilon_2}>0$ such that for all $i\in \{1,2,...,M_{\epsilon_2}\}$,
{\footnotesize
\begin{align*}
& \quad L\left(y\right)-\sup\limits_{\theta\in \left(0,1\right)}\left[\frac{\theta \mu e^{\lambda_1 y}}{\lambda_1-\lambda_0}+ \lfloor \left(1+\epsilon_1\right)z+i\delta_{M_{\epsilon_2}} \rfloor \frac{\mu}{\lambda_0}\log\left( -\lambda_0 \int_{0}^{\infty}\frac{e^{\lambda_1 s}}{e^{\lambda_1 s}-\theta}e^{\lambda_0 s}ds\right)\right]\\
& \quad \quad +\frac{\mu}{\lambda_0}\left(\lfloor \left(1+\epsilon_1\right)z+\left(i-1\right)\delta_{M_{\epsilon_2}} \rfloor\log\left(\lfloor \left(1+\epsilon_1\right)z+\left(i-1\right)\delta_{M_{\epsilon_2}} \rfloor\right)-\lfloor \left(1+\epsilon_1\right)z+\left(i-1\right)\delta_{M_{\epsilon_2}} \rfloor+1\right)\\
& < -\epsilon_2.
\end{align*}}
Hence, we conclude that 
\begin{align*}
& \quad \limsup_{n\rightarrow \infty}\frac{1}{n^{1-\alpha}}\log \PP\left(\gamma_n\le \zeta_n-y, -\frac{\lambda_0}{\mu}n^{\alpha-1}I_n\left(\zeta_n-y\right)\in O_{i,M}\right)\\
& < -L\left(y\right).
\end{align*}
We can get a similar result for the interval $\left(0, \left(1-\epsilon_1\right)z\right)$ and the desired result follows.
\qed
\end{proof}

\subsubsection{Proof of Lemma \ref{unique_sol}}\label{lemma_unique_sol}

\begin{proof}\\
Since $\Lambda_n''\left(\theta\right)>0$, $\Lambda_n'\left(\theta\right)$ increases in $\theta$. By \eqref{bound:q_k_n}, \eqref{def:q^star}, \eqref{def:Q_k_n}, and the fact that $q^*>1$, we can obtain that $\Lambda_n'\left(0\right)=\frac{\int_0^{\zeta_n-y}e^{-\left(\lambda_1-\lambda_0\right)t}dt}{\int_0^{\zeta_n-y}e^{\lambda_0 t}dt}<Q_{k,n}$ when $\epsilon_1$ and $\epsilon_2$ are sufficiently small, and $n$ is sufficiently large. Since
\begin{align*}
& \lim\limits_{\theta\rightarrow 1}\Lambda'\left(\theta\right)=\infty, \text{ and }\\
& \lim\limits_{n\rightarrow \infty}\Lambda_n'\left(\theta\right)=\Lambda'\left(\theta\right),
\end{align*}
we can find a sufficiently small $\epsilon>0$ such that for sufficiently large $n$,
\begin{align*}
\Lambda_n'\left(1-\epsilon\right)>Q_{k,n},
\end{align*}
which leads to the desired result. 
\qed
\end{proof}

\subsubsection{Proof of Lemma \ref{limit_eta_q}}\label{lemma_limit_eta_q}

\begin{proof}\\
By the definition of $\Lambda'\left(\theta \right)$ and $Q^*$, and \eqref{theta_2}, it is easy to verify that $\Lambda'\left(\theta^*_y \right)=Q^*$.
\qed
\end{proof}

\subsubsection{Proof of Lemma \ref{uniform_bounds_parameters_1}}\label{lemma_uniform_bounds_parameters_1}
\begin{proof}\\
The result for $Q_{k,n}$ follows immediately from the definition for $Q_{k,n}$ and $Q^*$. We know that $\eta_{k,n}={\Lambda_n'}^{-1}\left(Q_{k,n}\right)$. Since $\Lambda_n'\left(\theta\right)$ is a continuous and strictly increasing function, its inverse is also continuous and strictly increasing. The result for $\eta_{k,n}$ then follows immediately from the result for $Q_{k,n}$. The rest of the results for $\Lambda_n\left(\eta_{k,n}\right)$, and $\Lambda_n'\left(\eta_{k,n}\right)$ follow immediately from their monotonicity. The result of  $\Lambda_n''\left(\eta_{k,n}\right)$ follows from \eqref{bound_third}.
\qed

\end{proof}

\subsubsection{Proof of Proposition \ref{uniform_bound}}\label{proposition_uniform_bound}

\begin{proof}\\
We follow the main idea of the proof of Theorem 3.7.4 (Bahadur and Rao) of \cite{LD}. Define
\begin{align*}
Y_i=\frac{e^{-\lambda_1\left(\zeta_n-y\right)}X_i-Q_{k,n}}{\sqrt{\Lambda_n''\left(\eta_{k,n}\right)}}.
\end{align*}
It is easy to verify that $\EE_{\tilde{\nu}_{n,y,k}}\left[Y_i\right]=0$, $\EE_{\tilde{\nu}_{n,y,k}}\left[Y_i^2\right]=1$, and $\EE_{\tilde{\nu}_{n,y,k}}\left[Y_i^3\right]$ exists. We let
\begin{align*}
& \alpha_{3,n,y,q_{k,n}}=\EE_{\tilde{\nu}_{n,y,k}}\left[Y_i^3\right],\\
& \beta_{3,n,y,q_{k,n}}=\EE_{\tilde{\nu}_{n,y,k}}\left[\left|Y_i\right|^3\right].
\end{align*}
By Lemma \ref{limit_eta_q}, Lemma \ref{uniform_bounds_parameters_1} and the previous calculation of derivatives of $\Lambda_n\left(\theta\right)$, we can obtain useful bounds for $\alpha_{3,n,y,q_{k,n}}$ and $\beta_{3,n,y,q_{k,n}}$. Let the probability measure $\tilde{\nu}_{n,y}$ be defined by $d\tilde{\nu}_{n,y}/d\nu_{n,y} \left(x\right)=e^{\theta_y^* x - \Lambda\left(\eta_y^*\right)}$.
\begin{lemma}\label{3rd_moment_bound}
Define 
$$
\alpha_{3,n,y,q_{k,n}}\left[\epsilon_1\right]=\left(\inf\limits_{k\in O_{n,y,\epsilon_1}}\alpha_{3,n,y,q_{k,n}}, \sup\limits_{k\in O_{n,y,\epsilon_1}}\alpha_{3,n,y,q_{k,n}}\right).
$$ 
For any $\epsilon>0$, there exists $\delta_{\epsilon}>0$ and $N_{\epsilon}>0$ such that when $\epsilon_1<\delta_{\epsilon}$ and $n>N_{\epsilon}$,
\begin{align*}
& \alpha_{3,n,y,q_{k,n}}\left[\epsilon_1\right]\subset \left(\alpha_3^*-\epsilon, \alpha_3^*+\epsilon\right),
\end{align*}
where
\begin{align*}
\alpha_3^*=\lim\limits_{n\rightarrow\infty}\EE_{\tilde{\nu}_{n,y}}\left[\left(\frac{e^{-\lambda_1\left(\zeta_n-y\right)}X_i-Q^*}{\sqrt{\Lambda_n''\left(\theta_y^*\right)}}\right)^3\right].
\end{align*}
Moreover, 
\begin{align*}
& \limsup\limits_{\epsilon_1\rightarrow 0}   \limsup\limits_{n\rightarrow \infty}\sup\limits_{k\in O_{n,y,\epsilon_1}}\beta_{3,n,y,q_{k,n}}\le \lim\limits_{n\rightarrow \infty}\sqrt{\EE_{\tilde{\nu}_{n,y}}\left[\left(\frac{e^{-\lambda_1\left(\zeta_n-y\right)}X_i-Q^*}{\sqrt{\Lambda_n''\left(\theta_y^*\right)}}\right)^6\right]}, \text{ and } \\
& \beta_{3,n,y,q_{k,n}}\ge 1.
\end{align*}
\end{lemma}
\begin{proof}
See Section \ref{lemma_3rd_moment_bound}.
\end{proof}
\\

For ease of exposition, we omit the subscript $n$, $y$, and $q_{k,n}$, and just refer to them as $\alpha_3$ and $\beta_3$. Let $F_{k,n}$ be the distribution function of 
$$
k^{-\frac{1}{2}}\sum\limits_{i=1}^kY_i
$$ 
when $e^{-\lambda_1\left(\zeta_n-y\right)}X_i$ are i.i.d. with marginal law $\tilde{\nu}_{n,y,k}$. Let $\phi\left(x\right)$ be the pdf and $\Phi\left(x\right)$ be the cdf of a standard normal distribution. We know that
\begin{align*}
\phi'\left(x\right)=-\frac{x}{\sqrt{2\pi}}e^{-\frac{x^2}{2}}.
\end{align*}
A key step in our proof is to show a modified version of the Berry-Esseen expansion for lattice distributions (Theorem 3 in Chapter IV of \cite{Esseen1945} and Theorem 1 in \S 43 of \cite{Limit_distribution}). Compared to the original version, we need to consider a parameter-dependent span (notice that $Y_i$ has a span of $\frac{e^{-\lambda_1\left(\zeta_n-y\right)}}{\sqrt{\Lambda_n''\left(\eta_{k,n}\right)}}$) instead of a fixed one. Moreover, we need the convergence to be uniform for $k\in O_{n,y,\epsilon_1}$ as $n$ goes to infinity. We first introduce some notation. Let 
\begin{align*}
& h_{k,n}=\frac{e^{-\lambda_1\left(\zeta_n-y\right)}}{ \sqrt{\Lambda_n''\left(\eta_{k,n}\right)}},\\
& \tau_{k,n}=\frac{2\pi}{h_{k,n}},\\
& x_{k,n}=-\frac{Q_{k,n}}{\sqrt{\Lambda_n''\left(\eta_{k,n}\right)}}, \text{ and }\\
& S\left(x\right)=\left[x\right]-x+\frac{1}{2}.
\end{align*}
where $\left[x\right]$ gives the integral part of $x$. Notice that $h_{k,n}$ is the span of $Y_i$ and $\tau_{k,n}$ is the period of its 
. Let 
\begin{align*}
S_{k,n}\left(x\right)=h_{k,n}S\left(\frac{x\sqrt{k}-x_{k,n}k}{h_{k,n}}\right), \text{ and }
\end{align*}
\begin{align*}
D_{k,n}\left(x\right)=\phi\left(x\right)\frac{S_{k,n}\left(x\right)}{\sqrt{k}}.
\end{align*}
\begin{lemma}\label{Berry_Esseen_lat}
For any $\epsilon>0$, there exists $N_{\epsilon}>0$ and $\delta_{\epsilon}>0$ such that when $n>N_{\epsilon}$ and $\epsilon_1<\delta_{\epsilon}$,
\begin{align*}
\sup\limits_{k\in O_{n,y,\epsilon_1}}\Bigg\{\sqrt{k}\sup_{x}\left|F_{k,n}\left(x\right)-\Phi\left(x\right)-\frac{\alpha_3}{6\sqrt{k}}\left(1-x^2\right)\phi\left(x\right)-D_{k,n}\left(x\right)\right| \Bigg\}< \epsilon.
\end{align*}
\end{lemma}
\begin{proof}
See Section \ref{lemma_Berry_Esseen_lat}.
\end{proof}

Let
\begin{align*}
& \psi_{k,n}=\eta_{k,n} \sqrt{k\Lambda_n''\left(\eta_{k,n}\right)}, \text{ and }\\
& J_{n,y,k}=\eta_{k,n}\sqrt{\Lambda_n''\left(\eta_{k,n}\right)2\pi k}e^{k\Lambda^*_n\left(Q_{k,n}\right)},
\end{align*}
where $\Lambda^*_n\left(Q_{k,n}\right)=\eta_{k,n} Q_{k,n} -\Lambda_n\left(\eta_{k,n}\right)$. Let

{\footnotesize
\begin{align*}
C_{n,y,k}& =\sqrt{2\pi}\int_0^{\infty}\psi_{k,n} e^{-t} \left(\Phi\left(\frac{t}{\psi_{k,n}}\right)+\frac{\alpha_{3}}{6\sqrt{k}}\left(1-\left(\frac{t}{\psi_{k,n}}\right)^2\right)\phi\left(\frac{t}{\psi_{k,n}}\right)+\phi\left(\frac{t}{\psi_{k,n}}\right)\frac{S_{k,n}\left(\frac{t}{\psi_{k,n}}\right)}{\sqrt{k}}\right)dt\\
& \quad -\sqrt{2\pi}\int_0^{\infty}\psi_{k,n} e^{-t}\left(\Phi\left(0\right)+\frac{\alpha_{3}}{6\sqrt{k}}\phi\left(0\right)+\phi\left(0\right)\frac{S_{k,n}\left(0\right)}{\sqrt{k}}\right)dt.
\end{align*}}
We then show the following lemma which provides a good approximation for
\begin{align*}
\PP\left(\sum\limits_{i=1}^{k}e^{-\lambda_1\left(\zeta_n-y\right)}X_i\ge k Q_{k,n}\right).
\end{align*}
\begin{lemma}\label{J_c_bound}
For any $\epsilon>0$, there exists $N_{\epsilon}>0$ and $\delta_{\epsilon}>0$ such that when $n>N_{\epsilon}$ and $\epsilon_1<\delta_{\epsilon}$,
\begin{align*}
\sup\limits_{k\in O_{n,y,\epsilon_1}}\left| \PP\left(\sum_{i=1}^{k}e^{-\lambda_1\left(\zeta_n-y\right)}X_i\ge k Q_{k,n}\right)J_{n,y,k}-C_{n,y,k} \right| \le \epsilon.
\end{align*}
\end{lemma}
\begin{proof}
See Section \ref{lemma_J_c_bound}.
\end{proof}
\\

We then analyze $C_{n,y,k}$. By a Taylor expansion of $\Phi\left(\frac{t}{\psi_{k,n}}\right)$, we have
{\scriptsize
\begin{align*}
&\Phi\left(\frac{t}{\psi_{k,n}}\right)+\frac{\alpha_{3}}{6\sqrt{k}}\left(1-\left(\frac{t}{\psi_{k,n}}\right)^2\right)\phi\left(\frac{t}{\psi_{k,n}}\right)+\phi\left(\frac{t}{\psi_{k,n}}\right)\frac{S_{k,n}\left(\frac{t}{\psi_{k,n}}\right)}{\sqrt{k}}-\Phi\left(0\right)-\frac{\alpha_{3}}{6\sqrt{k}}\phi\left(0\right)-\phi\left(0\right)\frac{S_{k,n}\left(0\right)}{\sqrt{k}}\\
=&\frac{\alpha_{3}}{6\sqrt{k}}\left(1-\left(\frac{t}{\psi_{k,n}}\right)^2\right)\phi\left(\frac{t}{\psi_{k,n}}\right)-\frac{z\phi\left(z\right)}{2}\left(\frac{t}{\psi_{k,n}}\right)^2+\phi\left(0\right)\left(\frac{t}{\psi_{k,n}}\right)-\frac{\alpha_{3}}{6\sqrt{k}}\phi\left(0\right)+\phi\left(\frac{t}{\psi_{k,n}}\right)\frac{S_{k,n}\left(\frac{t}{\psi_{k,n}}\right)}{\sqrt{k}}-\phi\left(0\right)\frac{S_{k,n}\left(0\right)}{\sqrt{k}},
\end{align*}}
where $z$ is some number between $0$ and $\frac{t}{\psi_{k,n}}$. Dividing the integral in $C_{n,y,k}$ into two parts, we have
{\scriptsize
\begin{align*}
C_{n,y,k}& = \sqrt{2\pi}\int_0^{\infty}\psi_{k,n}e^{-t}\left(\frac{\alpha_{3}}{6\sqrt{k}}\left(1-\left(\frac{t}{\psi_{k,n}}\right)^2\right)\phi\left(\frac{t}{\psi_{k,n}}\right)+\phi\left(0\right)\left(\frac{t}{\psi_{k,n}}\right)+\phi\left(\frac{t}{\psi_{k,n}}\right)\frac{S_{k,n}\left(\frac{t}{\psi_{k,n}}\right)}{\sqrt{k}}\right)dt\\
& \quad - \sqrt{2\pi}\int_0^{\infty}\psi_{k,n}e^{-t}\left(\frac{z\phi\left(z\right)}{2}\left(\frac{t}{\psi_{k,n}}\right)^2+\frac{\alpha_{3}}{6\sqrt{k}}\phi\left(0\right)+\phi\left(0\right)\frac{S_{k,n}\left(0\right)}{\sqrt{k}}\right)dt\\
&= \sqrt{2\pi}\int_0^{T}\psi_{k,n}e^{-t}\left(\frac{\alpha_{3}}{6\sqrt{k}}\left(1-\left(\frac{t}{\psi_{k,n}}\right)^2\right)\phi\left(\frac{t}{\psi_{k,n}}\right)+\phi\left(0\right)\left(\frac{t}{\psi_{k,n}}\right)+\phi\left(\frac{t}{\psi_{k,n}}\right)\frac{S_{k,n}\left(\frac{t}{\psi_{k,n}}\right)}{\sqrt{k}}\right)dt\\
& \quad - \sqrt{2\pi}\int_0^{T}\psi_{k,n}e^{-t}\left(\frac{z\phi\left(z\right)}{2}\left(\frac{t}{\psi_{k,n}}\right)^2+\frac{\alpha_{3}}{6\sqrt{k}}\phi\left(0\right)+\phi\left(0\right)\frac{S_{k,n}\left(0\right)}{\sqrt{k}}\right)dt\\
& \quad + \sqrt{2\pi}\int_T^{\infty}\psi_{k,n}e^{-t}\left(\frac{\alpha_{3}}{6\sqrt{k}}\left(1-\left(\frac{t}{\psi_{k,n}}\right)^2\right)\phi\left(\frac{t}{\psi_{k,n}}\right)+\phi\left(0\right)\left(\frac{t}{\psi_{k,n}}\right)+\phi\left(\frac{t}{\psi_{k,n}}\right)\frac{S_{k,n}\left(\frac{t}{\psi_{k,n}}\right)}{\sqrt{k}}\right)dt\\
& \quad - \sqrt{2\pi}\int_T^{\infty}\psi_{k,n}e^{-t}\left(\frac{z\phi\left(z\right)}{2}\left(\frac{t}{\psi_{k,n}}\right)^2+\frac{\alpha_{3}}{6\sqrt{k}}\phi\left(0\right)+\phi\left(0\right)\frac{S_{k,n}\left(0\right)}{\sqrt{k}}\right)dt.
\end{align*}}
where $T$ can be arbitrarily large but fixed. By the definition of $\phi\left(x\right)$, the second part goes to $0$ as $T$ goes to infinity (note that the convergence is uniform for $k\in O_{n,y,\epsilon_1}$ for sufficiently large $n$ and sufficiently small $\epsilon_1$ by Lemma \ref{uniform_bounds_parameters_1} and Lemma \ref{3rd_moment_bound}). We then analyze the first part. By dominated convergence theorem and Lemma \ref{uniform_bounds_parameters_1},   
{\scriptsize
\begin{align*}
& \quad \lim\limits_{n\rightarrow \infty}\sqrt{2\pi}\int_0^{T}\psi_{k,n}e^{-t}\left(\frac{\alpha_{3}}{6\sqrt{k}}\left(1-\left(\frac{t}{\psi_{k,n}}\right)^2\right)\phi\left(\frac{t}{\psi_{k,n}}\right)+\phi\left(0\right)\left(\frac{t}{\psi_{k,n}}\right)+\phi\left(\frac{t}{\psi_{k,n}}\right)\frac{S_{k,n}\left(\frac{t}{\psi_{k,n}}\right)}{\sqrt{k}}\right)dt\\
& \quad \quad - \lim\limits_{n\rightarrow \infty}\sqrt{2\pi}\int_0^{T}\psi_{k,n}e^{-t}\left(\frac{z\phi\left(z\right)}{2}\left(\frac{t}{\psi_{k,n}}\right)^2+\frac{\alpha_{3}}{6\sqrt{k}}\phi\left(0\right)+\phi\left(0\right)\frac{S_{k,n}\left(0\right)}{\sqrt{k}}\right)dt\\
& =\sqrt{2\pi}\phi\left(0\right)\int_0^{T}te^{-t}dt\\
& =-\left(T+1\right)e^{-T}+1,
\end{align*}}
where the convergence is uniform for $k\in O_{n,y,\epsilon_1}$. Hence, when $\epsilon_1$ is sufficiently small,
\begin{align}\label{C_bound}
\lim\limits_{n\rightarrow \infty}\sup\limits_{k\in O_{n,y,\epsilon_1}}\left|C_{n,y,k}-1\right|=0.
\end{align}

Following the proof of Theorem 3.7.4 of \cite{LD} we have
{\footnotesize
\begin{align*}
& \EE\left[\left(e^{-\lambda_1\left(\zeta_n-y\right)}X_i\right)^2\middle| \sum_{i=1}^{k}X_i\ge n\right]\\
=& 1/\PP\left(\sum\limits_{i=1}^{k}X_i\ge n\right) \times \int_0^{\infty}x^2e^{-\eta_{k,n} x +\Lambda_n\left(\eta_{k,n}\right)}e^{-\left(k-1\right)\Lambda^*_n\left(Q_{k,n}\right)}\int_{\frac{Q_{k,n}-x}{\sqrt{\left(k-1\right)\Lambda_n''\left(\eta_{k,n}\right)}}}^{\infty}e^{-\sqrt{\frac{k-1}{k}}\psi_{k,n} z}dF_{k-1,n}\left(z\right)d\tilde{\nu}_{n,y,k}\left(x\right),
\end{align*}}
where $\Lambda^*_n\left(Q_{k,n}\right)=\eta_{k,n} Q_{k,n} -\Lambda_n\left(\eta_{k,n}\right)$, and the lower bound of the inner integral comes from the following calculation:
\begin{align*}
\left(k-1\right)^{-1/2}\sum\limits_{i=1}^{k-1} Y_i & =\left(k-1\right)^{-1/2}\sum\limits_{i=1}^{k-1} \left(e^{-\lambda_1\left(\zeta_n-y\right)}X_i-Q_{k,n}\right)/\sqrt{\Lambda_n''\left(\eta_{k,n}\right)}\\
& \ge \frac{1}{\sqrt{\left(k-1\right) \Lambda_n''\left(\eta_{k,n}\right)}}\left(Q_{k,n}k-x\right)-\frac{\sqrt{k-1}Q_{k,n}}{\sqrt{ \Lambda_n''\left(\eta_{k,n}\right)}}\\
& =\frac{Q_{k,n}-x}{\sqrt{\left(k-1\right)\Lambda_n''\left(\eta_{k,n}\right)}}.
\end{align*}
By Lemma \ref{J_c_bound} and \ref{C_bound}, we know that for any $\epsilon_3>0$, the following holds for sufficiently small $\epsilon_1$ and sufficiently large $n$:
\begin{align*}
\sup\limits_{k\in O_{n,y,\epsilon_1} }\left|\PP\left(\sum_{i=1}^{k}e^{-\lambda_1\left(\zeta_n-y\right)}X_i\ge k Q_{k,n}\right)J_{n,y,k}-1\right|<\epsilon_3.
\end{align*}
Hence, we will work with
{\scriptsize
\begin{align*}
& \quad \EE\left[\left(e^{-\lambda_1\left(\zeta_n-y\right)}X_i\right)^2\middle| \sum_{i=1}^{k}X_i\ge n\right]J_{n,y,k}\PP\left(\sum_{i=1}^{k}e^{-\lambda_1\left(\zeta_n-y\right)}X_i\ge k Q_{k,n}\right)\\
& =\eta_{k,n}\sqrt{2\pi k\Lambda_n''\left(\eta_{k,n}\right)}e^{k\Lambda^*_n\left(Q_{k,n} \right)}\int_0^{\infty}x^2e^{-\eta_{k,n} x +\Lambda_n\left(\eta_{k,n}\right)}e^{-\left(k-1\right)\Lambda^*_n\left(Q_{k,n}\right)}\int_{\frac{Q_{k,n}-x}{\sqrt{\left(k-1\right)\Lambda_n''\left(\eta_{k,n}\right)}}}^{\infty}e^{-\sqrt{\frac{k-1}{k}}\psi_{k,n} z}dF_{k-1,n}\left(z\right)d\tilde{\nu}_{n,y,k}\left(x\right)\\
&=\eta_{k,n}\sqrt{2\pi k\Lambda_n''\left(\eta_{k,n}\right)}e^{\Lambda^*_n\left(Q_{k,n} \right)}\int_0^{\infty}x^2e^{-\eta_{k,n} x +\Lambda_n\left(\eta_{k,n}\right)}\int_{\frac{Q_{k,n}-x}{\sqrt{\left(k-1\right)\Lambda_n''\left(\eta_{k,n}\right)}}}^{\infty}e^{-\sqrt{\frac{k-1}{k}}\psi_{k,n} z}dF_{k-1,n}\left(z\right)d\tilde{\nu}_{n,y,k}\left(x\right),
\end{align*}}
where we use $\PP\left(\sum_{i=1}^{k}e^{-\lambda_1\left(\zeta_n-y\right)}X_i\ge k Q_{k,n}\right)=\PP\left(\sum\limits_{i=1}^{k}X_i\ge n\right)$ in the first equality. Applying an integration by parts to the inner integral, we have
{\footnotesize
\begin{align}
& \quad \eta_{k,n}\sqrt{2\pi k\Lambda_n''\left(\eta_{k,n}\right)}e^{\Lambda^*_n\left(Q_{k,n} \right)}\int_0^{\infty}x^2e^{-\eta_{k,n} x +\Lambda_n\left(\eta_{k,n}\right)}\int_{\frac{Q_{k,n}-x}{\sqrt{\left(k-1\right)\Lambda_n''\left(\eta_{k,n}\right)}}}^{\infty}e^{-\sqrt{\frac{k-1}{k}}\psi_{k,n} z}dF_{k-1,n}\left(z\right)d\tilde{\nu}_{n,y,k}\left(x\right)\nonumber \\
&=\eta_{k,n}\sqrt{2\pi k\Lambda_n''\left(\eta_{k,n}\right)}e^{\Lambda^*_n\left(Q_{k,n} \right)}\int_0^{\infty}x^2e^{-\eta_{k,n} x +\Lambda_n\left(\eta_{k,n}\right)}\nonumber \\
& \quad \quad \quad \times \int_{\frac{Q_{k,n}-x}{\sqrt{\left(k-1\right)\Lambda_n''\left(\eta_{k,n}\right)}}}^{\infty}\sqrt{\frac{k-1}{k}}\psi_{k,n} e^{-\sqrt{\frac{k-1}{k}}\psi_{k,n} z}\left[F_{k-1,n}\left(z\right)-F_{k-1,n}\left(\frac{Q_{k,n}-x}{\sqrt{\left(k-1\right)\Lambda_n''\left(\eta_{k,n}\right)}}\right)\right]dz d\tilde{\nu}_{n,y,k}\left(x\right).\label{expression:before_change}
\end{align}}
By a change of variable ($t=\sqrt{\frac{k-1}{k}}\psi_{k,n} z$),
{\footnotesize
\begin{align}
\eqref{expression:before_change} & =\frac{\eta_{k,n}\sqrt{2\pi k\Lambda_n''\left(\eta_{k,n}\right)}}{\eta_{k,n}\sqrt{\left(k-1\right)\Lambda_n''\left(\eta_{k,n}\right)}}e^{\Lambda^*_n\left(Q_{k,n} \right)}\int_0^{\infty}x^2e^{-\eta_{k,n} x +\Lambda_n\left(\eta_{k,n}\right)} \nonumber\\
& \quad \quad \quad \times \int_{\eta_{k,n}\left(Q_{k,n} -x\right)}^{\infty}\sqrt{\frac{k-1}{k}}\psi_{k,n} e^{-t}\left[F_{k-1,n}\left(\frac{t}{\sqrt{\frac{k-1}{k}}\psi_{k,n}}\right)-F_{k-1,n}\left(\frac{Q_{k,n}-x}{\sqrt{\left(k-1\right)\Lambda_n''\left(\eta_{k,n}\right)}}\right)\right]dt d\tilde{\nu}_{n,y,k}\left(x\right)\nonumber\\
& =\sqrt{2\pi \frac{k}{k-1}}e^{\Lambda^*_n\left(Q_{k,n} \right)}\int_0^{\infty}x^2e^{-\eta_{k,n} x +\Lambda_n\left(\eta_{k,n}\right)}\\
& \quad \quad \quad \times \int_{\eta_{k,n}\left(Q_{k,n} -x\right)}^{\infty}\sqrt{\frac{k-1}{k}}\psi_{k,n} e^{-t}\left[F_{k-1,n}\left(\frac{t}{\sqrt{\frac{k-1}{k}}\psi_{k,n}}\right)-F_{k-1,n}\left(\frac{Q_{k,n}-x}{\sqrt{\left(k-1\right)\Lambda_n''\left(\eta_{k,n}\right)}}\right)\right]dt d\tilde{\nu}_{n,y,k}\left(x\right). \label{formula_final}
\end{align}}
Let
{\footnotesize
\begin{align*}
L_{1,n,y,k}\left(x\right)=e^{-\eta_{k,n} x}\int_{\eta_{k,n}\left(Q_{k,n} -x\right)}^{\infty}\sqrt{\frac{k-1}{k}}\psi_{k,n} e^{-t}\left[F_{k-1,n}\left(\frac{t}{\sqrt{\frac{k-1}{k}}\psi_{k,n}}\right)-F_{k-1,n}\left(\frac{Q_{k,n}-x}{\sqrt{\left(k-1\right)\Lambda_n''\left(\eta_{k,n}\right)}}\right)\right]dt.
\end{align*}}
Let
{\tiny
\begin{align*}
L_{2,n,y,k}\left(x\right)& =  e^{-\eta_{k,n} x}\int_{\eta_{k,n}\left(Q_{k,n} -x\right)}^{\infty} \sqrt{\frac{k-1}{k}}\psi_{k,n}e^{-t}\Bigg(\Phi\left(\frac{t}{\sqrt{\frac{k-1}{k}}\psi_{k,n}}\right)+\frac{\alpha_{3}}{6\sqrt{k-1}}\left(1-\left(\frac{t}{\sqrt{\frac{k-1}{k}}\psi_{k,n}}\right)^2\right)\phi\left(\frac{t}{\sqrt{\frac{k-1}{k}}\psi_{k,n}}\right)\\
&\quad \quad \quad  +\phi\left(\frac{t}{\sqrt{\frac{k-1}{k}}\psi_{k,n}}\right)\frac{S_{k-1,n}\left(\frac{t}{\sqrt{\frac{k-1}{k}}\psi_{k,n}}\right)}{\sqrt{k-1}}-\Phi\left(\frac{Q_{k,n}-x}{\sqrt{\left(k-1\right)\Lambda_n''\left(\eta_{k,n}\right)}}\right)\\
&\quad \quad \quad  -\frac{\alpha_{3}}{6\sqrt{k-1}}\left(1-\left(\frac{Q_{k,n}-x}{\sqrt{\left(k-1\right)\Lambda_n''\left(\eta_{k,n}\right)}}\right)^2\right)\phi\left(\frac{Q_{k,n}-x}{\sqrt{\left(k-1\right)\Lambda_n''\left(\eta_{k,n}\right)}}\right)-\phi\left(\frac{Q_{k,n}-x}{\sqrt{\left(k-1\right)\Lambda_n''\left(\eta_{k,n}\right)}}\right)\frac{S_{k-1,n}\left(\frac{Q_{k,n}-x}{\sqrt{\left(k-1\right)\Lambda_n''\left(\eta_{k,n}\right)}}\right)}{\sqrt{k-1}}\Bigg)dt.
\end{align*}}
By a change of variable ($z=t+\eta_{k,n} x$), we have
\begin{align*}
L_{1,n,y,k}\left(x\right)=\int_{\eta_{k,n} Q_{k,n}}^{\infty}\sqrt{\frac{k-1}{k}}\psi_{k,n} e^{-z}\left[F_{k-1,n}\left(\frac{z-\eta_{k,n} x}{\sqrt{\frac{k-1}{k}}\psi_{k,n}}\right)-F_{k-1,n}\left(\frac{Q_{k,n}-x}{\sqrt{\left(k-1\right)\Lambda_n''\left(\eta_{k,n}\right)}}\right)\right]dz,
\end{align*}
and
{\footnotesize
\begin{align*}
L_{2,n,y,k}\left(x\right)& = \int_{\eta_{k,n} Q_{k,n}}^{\infty} \sqrt{\frac{k-1}{k}}\psi_{k,n}e^{-z}\Bigg(\Phi\left(\frac{z-\eta_{k,n} x}{\sqrt{\frac{k-1}{k}}\psi_{k,n}}\right)+\frac{\alpha_{3}}{6\sqrt{k-1}}\left(1-\left(\frac{z-\eta_{k,n} x}{\sqrt{\frac{k-1}{k}}\psi_{k,n}}\right)^2\right)\phi\left(\frac{z-\eta_{k,n} x}{\sqrt{\frac{k-1}{k}}\psi_{k,n}}\right)\\
&\quad \quad \quad \quad +\phi\left(\frac{z-\eta_{k,n} x}{\sqrt{\frac{k-1}{k}}\psi_{k,n}}\right)\frac{S_{k-1,n}\left(\frac{z-\eta_{k,n} x}{\sqrt{\frac{k-1}{k}}\psi_{k,n}}\right)}{\sqrt{k-1}}-\Phi\left(\frac{Q_{k,n}-x}{\sqrt{\left(k-1\right)\Lambda_n''\left(\eta_{k,n}\right)}}\right)\\
&\quad \quad \quad \quad-\frac{\alpha_{3}}{6\sqrt{k-1}}\left(1-\left(\frac{Q_{k,n}-x}{\sqrt{\left(k-1\right)\Lambda_n''\left(\eta_{k,n}\right)}}\right)^2\right)\phi\left(\frac{Q_{k,n}-x}{\sqrt{\left(k-1\right)\Lambda_n''\left(\eta_{k,n}\right)}}\right)\\
&\quad \quad \quad-\phi\left(\frac{Q_{k,n}-x}{\sqrt{\left(k-1\right)\Lambda_n''\left(\eta_{k,n}\right)}}\right)\frac{S_{k-1,n}\left(\frac{Q_{k,n}-x}{\sqrt{\left(k-1\right)\Lambda_n''\left(\eta_{k,n}\right)}}\right)}{\sqrt{k-1}}\Bigg)dz.
\end{align*}}
By Lemma \ref{Berry_Esseen_lat}, for any $\epsilon_4>0$ and any $x$, when $\epsilon_1$ is sufficiently small and $n$ is sufficiently large,
\begin{align}\label{L_1_L_2}
\sup\limits_{k\in O_{n,y,\epsilon_1}}|L_{1,n,y,k}\left(x\right)-L_{2,n,y,k}\left(x\right)|< \epsilon_4.
\end{align}

We then analyze $L_{2,n,y,k}\left(x\right)$. Applying a Taylor expansion of $\Phi\left(\frac{y-\eta_{k,n} x}{\sqrt{\frac{k-1}{k}}\psi_{k,n}}\right)$:
{\footnotesize
\begin{align*}
\Phi\left(\frac{z-\eta_{k,n} x}{\sqrt{\frac{k-1}{k}}\psi_{k,n}}\right)&=\Phi\left(\frac{Q_{k,n}-x}{\sqrt{\left(k-1\right)\Lambda_n''\left(\eta_{k,n}\right)}}\right)+\phi\left(\frac{Q_{k,n}-x}{\sqrt{\left(k-1\right)\Lambda_n''\left(\eta_{k,n}\right)}}\right)\left(\frac{z-\eta_{k,n} x}{\sqrt{\frac{k-1}{k}}\psi_{k,n}}-\frac{Q_{k,n}-x}{\sqrt{\left(k-1\right)\Lambda_n''\left(\eta_{k,n}\right)}}\right)\\
& -\frac{z_x\phi\left(z_x\right)}{2}\left(\frac{z-\eta_{k,n} x}{\sqrt{\frac{k-1}{k}}\psi_{k,n}}-\frac{Q_{k,n}-x}{\sqrt{\left(k-1\right)\Lambda_n''\left(\eta_{k,n}\right)}}\right)^2\\
&=\Phi\left(\frac{Q_{k,n}-x}{\sqrt{\left(k-1\right)\Lambda_n''\left(\eta_{k,n}\right)}}\right)+\phi\left(\frac{Q_{k,n}-x}{\sqrt{\left(k-1\right)\Lambda_n''\left(\eta_{k,n}\right)}}\right)\left(\frac{z- \eta_{k,n} Q_{k,n}}{\sqrt{\frac{k-1}{k}}\psi_{k,n}}\right)-\frac{z_x\phi\left(z_x\right)}{2}\left(\frac{y-\eta_{k,n} Q_{k,n}}{\sqrt{\frac{k-1}{k}}\psi_{k,n}}\right)^2,
\end{align*}}
for some $z_x$ between $\frac{z-\eta_{k,n} x}{\sqrt{\frac{k-1}{k}}\psi_{k,n}}$ and $\frac{Q_{k,n}-x}{\sqrt{\left(k-1\right)\Lambda_n''\left(\eta_{k,n}\right)}}$. Dividing the integral into two parts, we have
{\tiny
\begin{align*}
L_{2,n,y,k}\left(x\right)& =  \int_{\eta_{k,n} Q_{k,n}}^{T} \sqrt{\frac{k-1}{k}}\psi_{k,n}e^{-z}\Bigg(\phi\left(\frac{Q_{k,n}-x}{\sqrt{\left(k-1\right)\Lambda_n''\left(\eta_{k,n}\right)}}\right)\left(\frac{z-\eta_{k,n} Q_{k,n}}{\sqrt{\frac{k-1}{k}}\psi_{k,n}}\right)+\frac{\alpha_{3}}{6\sqrt{k-1}}\left(1-\left(\frac{z-\eta_{k,n} x}{\sqrt{\frac{k-1}{k}}\psi_{k,n}}\right)^2\right)\phi\left(\frac{z-\eta_{k,n} x}{\sqrt{\frac{k-1}{k}}\psi_{k,n}}\right)\\
&\quad \quad \quad \quad \quad  -\frac{z_x\phi\left(z_x\right)}{2}\left(\frac{z-\eta_{k,n} Q_{k,n}}{\sqrt{\frac{k-1}{k}}\psi_{k,n}}\right)^2-\frac{\alpha_{3}}{6\sqrt{k-1}}\left(1-\left(\frac{Q_{k,n}-x}{\sqrt{\left(k-1\right)\Lambda_n''\left(\eta_{k,n}\right)}}\right)^2\right)\phi\left(\frac{Q_{k,n}-x}{\sqrt{\left(k-1\right)\Lambda_n''\left(\eta_{k,n}\right)}}\right)\\
&\quad \quad \quad \quad \quad +\phi\left(\frac{z-\eta_{k,n} x}{\sqrt{\frac{k-1}{k}}\psi_{k,n}}\right)\frac{S_{k-1,n}\left(\frac{z-\eta_{k,n} x}{\sqrt{\frac{k-1}{k}}\psi_{k,n}}\right)}{\sqrt{k-1}}-\phi\left(\frac{Q_{k,n}-x}{\sqrt{\left(k-1\right)\Lambda_n''\left(\eta_{k,n}\right)}}\right)\frac{S_{k-1,n}\left(\frac{Q_{k,n}-x}{\sqrt{\left(k-1\right)\Lambda_n''\left(\eta_{k,n}\right)}}\right)}{\sqrt{k-1}}\Bigg)dz\\
& +  \int_{T}^{\infty} \sqrt{\frac{k-1}{k}}\psi_{k,n}e^{-z}\Bigg(\phi\left(\frac{Q_{k,n}-x}{\sqrt{\left(k-1\right)\Lambda_n''\left(\eta_{k,n}\right)}}\right)\left(\frac{z-\eta_{k,n} Q_{k,n}}{\sqrt{\frac{k-1}{k}}\psi_{k,n}}\right)+\frac{\alpha_{3}}{6\sqrt{k-1}}\left(1-\left(\frac{z-\eta_{k,n} x}{\sqrt{\frac{k-1}{k}}\psi_{k,n}}\right)^2\right)\phi\left(\frac{z-\eta_{k,n} x}{\sqrt{\frac{k-1}{k}}\psi_{k,n}}\right)\\
&\quad \quad \quad \quad \quad  -\frac{z_x\phi\left(z_x\right)}{2}\left(\frac{z-\eta_{k,n} Q_{k,n}}{\sqrt{\frac{k-1}{k}}\psi_{k,n}}\right)^2-\frac{\alpha_{3}}{6\sqrt{k-1}}\left(1-\left(\frac{Q_{k,n}-x}{\sqrt{\left(k-1\right)\Lambda_n''\left(\eta_{k,n}\right)}}\right)^2\right)\phi\left(\frac{Q_{k,n}-x}{\sqrt{\left(k-1\right)\Lambda_n''\left(\eta_{k,n}\right)}}\right)\\
&\quad \quad \quad \quad \quad  +\phi\left(\frac{z-\eta_{k,n} x}{\sqrt{\frac{k-1}{k}}\psi_{k,n}}\right)\frac{S_{k-1,n}\left(\frac{z-\eta_{k,n} x}{\sqrt{\frac{k-1}{k}}\psi_{k,n}}\right)}{\sqrt{k-1}}-\phi\left(\frac{Q_{k,n}-x}{\sqrt{\left(k-1\right)\Lambda_n''\left(\eta_{k,n}\right)}}\right)\frac{S_{k-1,n}\left(\frac{Q_{k,n}-x}{\sqrt{\left(k-1\right)\Lambda_n''\left(\eta_{k,n}\right)}}\right)}{\sqrt{k-1}}\Bigg)dz,
\end{align*}}
where $T$ can be arbitrarily large but fixed. Let $x< V$, where $V$ can be arbitrarily large but fixed. Then by a similar analysis to that of $C_{n,y,k}$, we conclude that for any arbitrarily small but fixed $\epsilon_5>0$, there exists a sufficiently large $T$, such that for sufficiently small $\epsilon_1$ and sufficiently large $n$, the following holds for all $x<V$:
\begin{align}\label{L_2_bound}
\sup\limits_{k\in O_{n,y,\epsilon_1}}\left|L_{2,n,y,k}\left(x\right)-\frac{1}{\sqrt{2\pi}}e^{-\eta_{k,n} Q_{k,n}}\right|<\epsilon_5.
\end{align}
Now we have that
{\footnotesize
\begin{align*}
\eqref{formula_final} & =\sqrt{2\pi \frac{k}{k-1}}e^{\Lambda^*_n\left(Q_{k,n} \right)}\int_0^{V}x^2e^{\Lambda_n\left(\eta_{k,n}\right)}L_{1,n,y,k}\left(x\right) d\tilde{\nu}_{n,y,k}\left(x\right)\\
& \quad \quad \quad +\sqrt{2\pi \frac{k}{k-1}}e^{\Lambda^*_n\left(Q_{k,n} \right)}\int_V^{\infty}x^2e^{\Lambda_n\left(\eta_{k,n}\right)}L_{1,n,y,k}\left(x\right) d\tilde{\nu}_{n,y,k}\left(x\right).
\end{align*}}
Because $L_{1,n,y,k}\left(x\right)$ is bounded, $\EE_{\tilde{\nu}_{n,y,k}}\left[Y_i\right]=0$, and $\EE_{\tilde{\nu}_{n,y,k}}\left[Y_i^2\right]=1$, by Lemma \ref{uniform_bounds_parameters_1}, we conclude that the second term goes to zero as $V$ goes to infinity uniformly for $k\in O_{n,y,\epsilon_1}$. For the first term, since $V$ can be arbitrarily large, the desired result follows from Lemma \ref{uniform_bounds_parameters_1}, \eqref{L_1_L_2}, and \eqref{L_2_bound}.
\qed
\end{proof}

\subsubsection{Proof of Lemma \ref{3rd_moment_bound}}\label{lemma_3rd_moment_bound}
\begin{proof}\\
We first notice that
\begin{align*}
\alpha_{3,n,y,q_{k,n}}=\frac{rf_1\left(\eta_{k,n}\right)}{\left(\sqrt{\Lambda_n''\left(\eta_{k,n}\right)}\right)^3}\left(f_4\left(\eta_{k,n}\right)-3f_3^2\left(\eta_{k,n}\right){Q_{k,n}}+3g_2\left(\eta_{k,n}\right){Q_{k,n}}^2-{Q_{k,n}}^3\right),
\end{align*}
where $f_1$ through $f_4$ are defined in \eqref{f_1} through \eqref{f_4}. The result for $\alpha_{3,n,y,q_{k,n}}$ then follows from Lemma \ref{uniform_bounds_parameters_1} and the monotonicity and continuity of $f_1$ through $f_4$. Note that a similar argument can be applied to the sixth moment, and thus we can use Holder's inequality to obtain the upper bound for $\beta_{3,n,y,q_{k,n}}$. The lower bound for $\beta_{3,n,y,q_{k,n}}$ follows from Holder's inequality, and the fact that $\EE_{\tilde{\nu}_{n,y,k}}\left[Y_i^2\right]=1$.
\end{proof}

\subsubsection{Proof of Lemma \ref{Berry_Esseen_lat}}\label{lemma_Berry_Esseen_lat}

\begin{proof}\\
The proof follows that of Theorem 1 in \S 43 of \cite{Limit_distribution}. We adopt its structure and notation. Let $C$ be some positive constant whose value might change from line to line. Let
\begin{align*}
G_{k,n}\left(x\right)=\Phi\left(x\right)+\frac{\alpha_{3}}{6\sqrt{k}}\left(1-x^2\right)\phi\left(x\right)+D_{k,n}\left(x\right),
\end{align*}
and recall that
\begin{align*}
D_{k,n}\left(x\right)=\phi\left(x\right)\frac{S_{k,n}\left(x\right)}{\sqrt{k}}.
\end{align*}
We can obtain the Fourier-Stieltjes transform of $G_{k,n}$ (see page 214 of \cite{Limit_distribution}):
\begin{align*}
g_{k,n}\left(t\right)=e^{-\frac{t^2}{2}}+\frac{\alpha_3 \left(it\right)^3}{6\sqrt{k}}e^{-\frac{t^2}{2}}+d_{k,n}\left(t\right),
\end{align*}
where
\begin{align*}
d_{k,n}\left(t\right)&=-\frac{t}{\tau_{k,n}  \sqrt{2\pi k}}\sum_{v=-\infty}^{\infty}\frac{1}{v}\int_{-\infty}^{\infty}e^{it x-\frac{x^2}{2}+iv\tau_{k,n} \sqrt{k}\left(x-x_{k,n}\sqrt{k}\right)}dx\\
&=-\frac{t}{\tau_{k,n}  \sqrt{ k}}\sum_{v=-\infty}^{\infty}\frac{e^{-i\tau_{k,n}v k x_{k,n}}}{v}e^{-\frac{1}{2}\left(t+\tau_{k,n} v \sqrt{k}\right)^2},
\end{align*}
and the summation is over every integer $v\neq 0$. It is not hard to observe that 
$$
\sup\limits_{k\in O_{n,y,\epsilon_1}}\left|G'_{k,n}\left(x\right)\right|
$$ 
is bounded wherever the derivative exists for sufficiently large $n$ and sufficiently small $\epsilon_1$. We denote by $A$ its bound.
Let $T=n$, then $T\frac{h_{k,n}}{\sqrt{k}}$ (recall that $h_{k,n}$ is of order $O\left(n^{-\alpha}\right)$) can be arbitrarily large for sufficiently large $n$ (this is required to apply Theorem 2 in \S 39 of \cite{Limit_distribution}). Let $f_{k,n}\left(\cdot\right)$ be the  Fourier-Stieltjes transform of $F_{k,n}$. By Theorem 2 in \S 39 of \cite{Limit_distribution}, to prove the desired result, it suffices to show that for sufficiently small $\epsilon_1$ and sufficiently large $n$
\begin{align*}
\sup\limits_{k\in O_{n,y,\epsilon_1}}\sqrt{k}\int_{-T}^T\left|\frac{f_{k,n}\left(t\right)-g_{k,n}\left(t\right)}{t}\right|dt\le \frac{\epsilon}{2}.
\end{align*}
Following \cite{Limit_distribution}, let
\begin{align*}
& e_1=\int_{-T}^{-\frac{\tau_{k,n}}{2}\sqrt{k}}\left|\frac{f_{k,n}\left(t\right)-g_{k,n}\left(t\right)}{t}\right|dt,\\
& e_2=\int_{-\frac{\tau_{k,n}}{2}\sqrt{k}}^{\frac{\tau_{k,n}}{2}\sqrt{k}}\left|\frac{f_{k,n}\left(t\right)-g_{k,n}\left(t\right)}{t}\right|dt, \text{ and }\\
& e_3=\int_{\frac{\tau_{k,n}}{2}\sqrt{k}}^{T}\left|\frac{f_{k,n}\left(t\right)-g_{k,n}\left(t\right)}{t}\right|dt.
\end{align*}

We first analyze $e_2$. Let $T_{k}=\frac{\sqrt{k}}{24\beta_3}$. From Lemma \ref{3rd_moment_bound}, we know that $\beta_3$ is a finite positive number bounded away from zero. By the definition of $\tau_{k,n}$, for sufficiently large $n$, we have $T_{k}<\frac{\tau_{k,n}}{2}\sqrt{k}$. Recall that $F_{k,n}$ is the distribution function of 
$$
k^{-\frac{1}{2}}\sum\limits_{i=1}^kY_i.
$$ 
Let $\varphi_{k,n}\left(t\right)$ be the characteristic function of $Y_i$ when $e^{-\lambda_1\left(\zeta_n-y\right)}X_i$ has distribution $\tilde{\nu}_{n,y,k}$. We have the following result which gives an upper bound for $\left|\varphi_{k,n}\left(t\right)\right|$.
\begin{lemma}\label{marginal_cha_bound}
There exists $c_1>0$ such that for sufficiently large $n$ and sufficiently small $\epsilon_1$, when $T_{k}/ \sqrt{k}\le \left|t\right|\le \tau_{k,n}/2$.
\begin{align*}
\sup\limits_{k\in O_{n,y,\epsilon_1}}\left|\varphi_{k,n}\left(t\right)\right|<e^{-c_1}.
\end{align*}
\end{lemma}
\begin{proof}
See Section \ref{lemma_marginal_cha_bound}.
\end{proof}
\\

For sufficiently large $n$, by Lemma \ref{marginal_cha_bound}, the following holds for $T_{k}\le \left|t\right| \le \tau_{k,n}\sqrt{k}/2$:
\begin{align*}
\left|f_{k,n}\left(t\right)\right|=\left|\varphi_{k,n}\left(\frac{t}{\sqrt{k}}\right)\right|^k < e^{-c_1k}.
\end{align*}
We also note that for some $c_2>0$ and sufficiently large $n$ and $k\in O_{n,y,\epsilon_1}$, the following holds when $T_{k}\le \left|t\right| \le \tau_{k,n}\sqrt{k}/2$:
\begin{align*}
\left|g_{k,n}\left(t\right)\right|\le e^{-c_2 k}.
\end{align*}
Therefore, for some positive constant $c$,
\begin{align*}
e_2 \le \int_{-T_{k}}^{T_{k}}\left|\frac{f_{k,n}\left(t\right)-g_{k,n}\left(t\right)}{t}\right|dt+4\int_{T_{k}}^{\frac{1}{2}\tau_{k,n} \sqrt{k}}\frac{e^{-ck}}{t}dt.
\end{align*}
By Theorem 1(b) in \S41 of \cite{Limit_distribution}, there exists some positive constant $C$ such that for sufficiently large $n$ and $k\in O_{n,y,\epsilon_1}$,
\begin{align*}
\int_{-T_{k}}^{T_{k}}\left|\frac{f_{k,n}\left(t\right)-g_{k,n}\left(t\right)}{t}\right|dt\le \frac{C}{\sqrt{k}}\delta\left(k\right)+\int_{-T_{k}}^{T_{k}}\left|\frac{d_{k,n}\left(t\right)}{t}\right|dt,
\end{align*}
where $\delta\left(k\right)$ only depends on $k$, and $\lim\limits_{k\rightarrow \infty}\delta\left(k\right)=0$. We also notice that for $\left|t\right|\le T_{k}$, there exists some positive constant $C$ such that for sufficiently large $n$ and $k\in O_{n,y,\epsilon_1}$,
\begin{align*}
\left|\frac{d_{k,n}\left(t\right)}{t}\right|\le \frac{C}{\tau_{k,n}\sqrt{k}}e^{-\frac{t^2}{2}}.
\end{align*}
Therefore, for sufficiently large $n$ and $k\in O_{n,y,\epsilon_1}$,
\begin{align*}
e_2\le C\left(\frac{\delta\left(k\right)}{\sqrt{k}}+\log\left(n\right)e^{-ck}+\frac{1}{\tau_{k,n}\sqrt{k}}\right),
\end{align*}
which is of order $o\left(\frac{1}{\sqrt{k}}\right)$. 

We then estimate $\epsilon_3$. For sufficiently large $n$ and $k\in O_{n,y,\epsilon_1}$,
\begin{align*}
e_3& =\int_{\frac{\tau_{k,n}}{2}\sqrt{k}}^{n}\left|\frac{f_{k,n}\left(t\right)-g_{k,n}\left(t\right)}{t}\right|dt\\
& \le \frac{C}{k}+\int_{\frac{\tau_{k,n}}{2}\sqrt{k}}^{n}\left|\frac{f_{k,n}\left(t\right)-d_{k,n}\left(t\right)}{t}\right|dt \quad \quad \text{(Mill's Inequality)}\\
& =\frac{C}{k}+\int_{\frac{\tau_{k,n}}{2}}^{n/\sqrt{k}}\left|\frac{\varphi_{k,n}^k\left(t\right)-d_{k,n}\left( \sqrt{k}t\right)}{t}\right|dt\\
& =\frac{C}{k}+\sum_{j=1}^{J}\int_{\frac{2j-1}{2}\tau_{k,n}}^{\frac{2j+1}{2}\tau_{k,n}}\left|\frac{\varphi_{k,n}^k\left(t\right)-d_{k,n}\left( \sqrt{k}t\right)}{t}\right|dt + \int_{\frac{2J+1}{2}\tau_{k,n}}^{n/\sqrt{k}}\left|\frac{\varphi_{k,n}^k\left(t\right)-d_{k,n}\left( \sqrt{k}t\right)}{t}\right|dt,
\end{align*}
where 
\begin{align*}
J=\left[\frac{n/\sqrt{k}}{\tau_{k,n}}-\frac{1}{2}\right].
\end{align*}
Let 
\begin{align*}
I_{j}=\int_{\frac{2j-1}{2}\tau_{k,n}}^{\frac{2j+1}{2}\tau_{k,n}}\left|\frac{\varphi_{k,n}^k\left(t\right)-d_{k,n}\left( \sqrt{k}t\right)}{t}\right|dt.
\end{align*}
Apply a change of variable $t=z+j\tau_{k,n}$, we have
\begin{align*}
I_{j}& =\int_{-\frac{\tau_{k,n}}{2}}^{\frac{\tau_{k,n}}{2}}\left|\frac{e^{ij\tau_{k,n}x_{k,n}k}\varphi_{k,n}^k\left(z\right)+\frac{z+j\tau_{k,n}}{\tau_{k,n}}\sum\limits_v \frac{1}{v}e^{-iv\tau_{k,n}x_{k,n}k}e^{-\frac{k}{2}\left(z+j\tau_{k,n}+v\tau_{k,n}\right)^2}}{z+j\tau_{k,n}}\right|dz\\
& \le \int_{-\frac{\tau_{k,n}}{2}}^{\frac{\tau_{k,n}}{2}}\left|\frac{\varphi_{k,n}^k\left(z\right)-e^{-\frac{1}{2} k z^2}-\frac{z}{j\tau_{k,n}}e^{-\frac{1}{2} k z^2}}{z+j\tau_{k,n}}\right|dz+Ce^{-\frac{k\tau_{k,n}^2}{8}}.
\end{align*}
where the inequality is due to the fact that the term in the inner summation only makes considerable contribution to the integral when $v=-j$ (note that when $v\neq -j$, 
$$
e^{-\frac{k}{2}\left(z+j\tau_{k,n}+v\tau_{k,n}\right)^2}=O\left(e^{-\frac{k\tau_{k,n}^2}{8}}\right)
$$ 
uniformly in $z$; also see page 216 of \cite{Limit_distribution}). We first note that
\begin{align*}
\int_{-\frac{\tau_{k,n}}{2}}^{\frac{\tau_{k,n}}{2}}\left|\frac{\varphi_{k,n}^k\left(z\right)-e^{-\frac{1}{2} k z^2}}{z+j\tau_{k,n}}\right|dz\le \frac{2}{j\tau_{k,n}}\int_{-\frac{\tau_{k,n}}{2}\sqrt{k}}^{{\frac{\tau_{k,n}}{2}\sqrt{k}}}\left|f_{k,n}\left(z\right)-e^{-\frac{1}{2}z^2}\right|\frac{dz}{\sqrt{k}}.
\end{align*}
By a similar argument to the treatment of $e_2$, we can obtain that for sufficiently large $n$ and $k\in O_{n,y,\epsilon_1}$,
\begin{align*}
\int_{-\frac{\tau_{k,n}}{2}}^{\frac{\tau_{k,n}}{2}}\left|\frac{\varphi_{k,n}^k\left(z\right)-e^{-\frac{1}{2} k z^2}}{z+j\tau_{k,n}}\right|dz\le \frac{C}{jk}.
\end{align*}
We then note that for sufficiently large $n$ and $k\in O_{n,y,\epsilon_1}$,
\begin{align*}
\frac{1}{j\tau_{k,n}}\int_{-\frac{\tau_{k,n}}{2}}^{\frac{\tau_{k,n}}{2}}\left|\frac{z}{z+j\tau_{k,n}}\right|e^{-\frac{1}{2}kz^2}dz\le \frac{C}{j^2k}.
\end{align*}
Recall that $J=\left[\frac{n/\sqrt{k}}{\tau_{k,n}}-\frac{1}{2}\right]$. For sufficiently large $n$, we have $e_3 \le C \log\left(n\right)/k$.
By the same method, $e_1 \le C\log\left(n\right)/k$.
Therefore, for sufficiently large $n$ and $k\in O_{n,y,\epsilon_1}$
\begin{align*}
\int_{-T}^T\left|\frac{f_{k,n}\left(t\right)-g_{k,n}\left(t\right)}{t}\right|dt\le C\left(\frac{\delta\left(k\right)}{\sqrt{k}}+\log\left(n\right)e^{-ck}+\frac{1}{\tau_{k,n}\sqrt{k}}+\frac{\log\left(n\right)}{k}\right).
\end{align*}
An application of Theorem 2 of \S 39 of \cite{Limit_distribution} leads to the desired result.
\qed

\end{proof}

\subsubsection{Proof of Lemma \ref{J_c_bound}}\label{lemma_J_c_bound}

\begin{proof}\\
Following the proof of Theorem 3.7.4 in \cite{LD}, we have
\begin{align*}
\PP\left(\sum_{i=1}^{k}e^{-\lambda_1\left(\zeta_n-y\right)}X_i\ge k Q_{k,n}\right)J_{n,y,k}=\sqrt{2\pi}\int_0^{\infty}\psi_{k,n} e^{-t}\left[F_{k,n}\left(\frac{t}{\psi_{k,n}}\right)-F_{k,n}\left(0\right)\right].
\end{align*}
The desired result follows directly by Lemma \ref{uniform_bounds_parameters_1}, Lemma \ref{3rd_moment_bound} and Lemma \ref{Berry_Esseen_lat}.
\qed
\end{proof}

\subsubsection{Proof of Lemma \ref{marginal_cha_bound}}\label{lemma_marginal_cha_bound}

\begin{proof}\\
\begin{align*}
\varphi_{k,n}\left(t\right)& = \EE_{\tilde{\nu}_{n,y,k}}\left[e^{it \frac{e^{-\lambda_1\left(\zeta_n-y\right)}X_i-Q_{k,n}}{\sqrt{\Lambda_n''\left(\eta_{k,n}\right)}}}\right]\\
& = \EE\left[e^{it \frac{e^{-\lambda_1\left(\zeta_n-y\right)}X_i-Q_{k,n}}{\sqrt{\Lambda_n''\left(\eta_{k,n}\right)}}}e^{\eta_{k,n} e^{-\lambda_1\left(\zeta_n-y\right)}X_i-\Lambda_n\left(\eta_{k,n}\right)}\right]\\
& = e^{-it \frac{Q_{k,n}}{\sqrt{\Lambda_n''\left(\eta_{k,n}\right)}}}e^{-\Lambda_n\left(\eta_{k,n}\right)} \EE\left[e^{it \frac{e^{-\lambda_1\left(\zeta_n-y\right)}X_i}{\sqrt{\Lambda_n''\left(\eta_{k,n}\right)}}}e^{\eta_{k,n} e^{-\lambda_1\left(\zeta_n-y\right)}X_i}\right]\\
& = e^{-it \frac{Q_{k,n}}{\sqrt{\Lambda_n''\left(\eta_{k,n}\right)}}}e^{-\Lambda_n\left(\eta_{k,n}\right)} \frac{\int_0^{\zeta_n-y}e^{\lambda_0 s}\frac{e^{-\lambda_1\left(\zeta_n-y-s\right)}e^{f\left(t\right)e^{-\lambda_1\left(\zeta_n-y\right)}}}{e^{-\lambda_1\left(\zeta_n-y-s\right)}e^{f\left(t\right)e^{-\lambda_1\left(\zeta_n-y\right)}}-e^{f\left(t\right)e^{-\lambda_1\left(\zeta_n-y\right)}}+1} ds}{\int_0^{\zeta_n-y}e^{\lambda_0 s}ds}\\
& = e^{-it \frac{Q_{k,n}}{\sqrt{\Lambda_n''\left(\eta_{k,n}\right)}}} \frac{\int_0^{\zeta_n-y}e^{\lambda_0 s}\frac{e^{-\lambda_1\left(\zeta_n-y-s\right)}e^{f\left(t\right)e^{-\lambda_1\left(\zeta_n-y\right)}}}{e^{-\lambda_1\left(\zeta_n-y-s\right)}e^{f\left(t\right)e^{-\lambda_1\left(\zeta_n-y\right)}}-e^{f\left(t\right)e^{-\lambda_1\left(\zeta_n-y\right)}}+1} ds}{\int_0^{\zeta_n-y}e^{\lambda_0 s}\frac{e^{-\lambda_1\left(\zeta_n-y-s\right)}e^{e^{-\lambda_1\left(\zeta_n-y\right)}\eta_{k,n}}}{e^{-\lambda_1\left(\zeta_n-y-s\right)}e^{e^{-\lambda_1\left(\zeta_n-y\right)}\eta_{k,n}}-e^{e^{-\lambda_1\left(\zeta_n-y\right)}\eta_{k,n}}+1} ds},
\end{align*}
where $f\left(t\right)=\frac{i t}{\sqrt{\Lambda_n''\left(\eta_{k,n}\right)}}+\eta_{k,n}$. Let
\begin{align*}
& A_{n,s}=e^{-\lambda_1\left(\zeta_n-y-s\right)},\\
& B_{k,n}=e^{e^{-\lambda_1\left(\zeta_n-y\right)}\eta_{k,n}},\\
& C_{k,n,t}=\cos\left(\frac{t}{\sqrt{\Lambda_n''\left(\eta_{k,n}\right)}}e^{-\lambda_1\left(\zeta_n-y\right)}\right)-1.
\end{align*}
For simplicity, we omit all subscripts and refer to them as $A$, $B$ and $C$. For some positive $T$,
\begin{align*}
\left| \varphi_{k,n}\left(t\right) \right| & \le \frac{\int_0^{\zeta_n-y}e^{\lambda_0 s}\frac{AB}{\sqrt{\left(AB-B+1\right)^2+2\left(AB-B\right)C}} ds}{\int_0^{\zeta_n-y}e^{\lambda_0 s}\frac{AB}{AB-B+1} ds}\\
& =\frac{\int_0^{T}e^{\lambda_0 s}\frac{AB}{\sqrt{\left(AB-B+1\right)^2+2\left(AB-B\right)C}} ds+\int_T^{\zeta_n-y}e^{\lambda_0 s}\frac{AB}{\sqrt{\left(AB-B+1\right)^2+2\left(AB-B\right)C}} ds}{\int_0^{T}e^{\lambda_0 s}\frac{AB}{AB-B+1} ds+\int_T^{\zeta_n-y}e^{\lambda_0 s}\frac{AB}{AB-B+1} ds}\\
& < \frac{\int_0^{T}e^{\lambda_0 s}\frac{AB}{\sqrt{\left(AB-B+1\right)^2+2\left(AB-B\right)C}} ds+\int_T^{\zeta_n-y}e^{\lambda_0 s}\frac{AB}{AB-B+1} ds}{\int_0^{T}e^{\lambda_0 s}\frac{AB}{AB-B+1} ds+\int_T^{\zeta_n-y}e^{\lambda_0 s}\frac{AB}{AB-B+1} ds}.
\end{align*}
We notice that for $T_{k}/ \sqrt{k}\le \left|t\right|\le \tau_{k,n}/2$, the final term achieves its maximum value when $t=T_{k}/ \sqrt{k}$ ($\left|C\right|$ is minimized at $t=T_{k}/ \sqrt{k}$). By Lemma \ref{uniform_bounds_parameters_1}, let $\eta_1$ and $\eta_2$ be the lower bound and upper bound for $\eta_{k,n}$ when $n$ is sufficiently large and $\epsilon_1$ is sufficiently small. By checking the sign of derivative with respect to $\eta_{k,n}$ and $n$ (with the help of Leibniz integral rule), we can obtain that
\begin{align}\label{AB_limit}
\int_0^{\zeta_n-y}e^{\lambda_0 s}\frac{AB}{AB-B+1} ds \le \int_0^{\infty} e^{\lambda_0 s}\frac{e^{\lambda_1 s}}{e^{\lambda_1 s}-\eta_2}ds.
\end{align}
Moreover, for fixed $t$,
\begin{align*}
\lim\limits_{n\rightarrow \infty}\frac{1-\cos\left(\frac{t}{\sqrt{\Lambda_n''\left(\eta\right)}}e^{-\lambda_1\left(\zeta_n-y\right)}\right)}{e^{-2\lambda_1\left(\zeta_n-y\right)}}& =\lim\limits_{n\rightarrow \infty}\frac{t^2}{2\Lambda_n''\left(\eta\right)}=\frac{t^2}{2\Lambda''\left(\eta\right)},
\end{align*}
where
\begin{align*}
\Lambda''\left(\eta\right)=\frac{2\int_0^{\infty}e^{\lambda_0 s}\frac{e^{-2\lambda_1 s}}{\left(1-\eta e^{-\lambda_1 s}\right)^3}ds\int_0^{\infty}e^{\lambda_0 s}\frac{1}{1-\eta e^{-\lambda_1 s}}ds-\left(\int_0^{\infty}e^{\lambda_0 s}\frac{e^{-\lambda_1 s}}{\left(1-\eta e^{-\lambda_1 s}\right)^2}ds\right)^2}{\left(\int_0^{\infty}e^{\lambda_0 s}\frac{1}{1-\eta e^{-\lambda_1 s}}ds\right)^2}.
\end{align*}
Hence, 
\begin{align*}
\left| \varphi_{k,n}\left(t\right) \right| & <1+\frac{\int_0^{T}e^{\lambda_0 s}\frac{AB}{\sqrt{\left(AB-B+1\right)^2+2\left(AB-B\right)C}} ds-\int_0^{T}e^{\lambda_0 s}\frac{AB}{AB-B+1} ds}{\int_0^{T}e^{\lambda_0 s}\frac{AB}{AB-B+1} ds+\int_T^{\zeta_n-y}e^{\lambda_0 s}\frac{AB}{AB-B+1} ds}\\
& = 1+\frac{\int_0^{T}e^{\lambda_0 s}\frac{AB}{\sqrt{\left(AB-B+1\right)^2+2\left(AB-B\right)C}} ds-\int_0^{T}e^{\lambda_0 s}\frac{AB}{AB-B+1} ds}{\int_0^{\zeta_n-y}e^{\lambda_0 s}\frac{AB}{AB-B+1} ds}.
\end{align*}
We notice that
\begin{align*}
& \limsup\limits_{n\rightarrow\infty}\sup_{\eta_{k,n}\in \left[\eta_1, \eta_2\right]}\int_0^{T}e^{\lambda_0 s}\frac{AB}{\sqrt{\left(AB-B+1\right)^2+2\left(AB-B\right)C}} ds-\int_0^{T}e^{\lambda_0 s}\frac{AB}{AB-B+1} ds\\
\ge & \sup_{\eta_{k,n}\in \left[\eta_1, \eta_2\right]}\lim\limits_{n\rightarrow \infty}\int_0^{T}e^{\lambda_0 s}\frac{AB}{\sqrt{\left(AB-B+1\right)^2+2\left(AB-B\right)C}} ds-\int_0^{T}e^{\lambda_0 s}\frac{AB}{AB-B+1} ds\\
= & \sup_{\eta\in \left[\eta_1, \eta_2\right]}\int_0^{T}e^{\lambda_0 s}\frac{e^{\lambda_1 s}}{\sqrt{\left(e^{\lambda_1 s}-\eta\right)^2 +\frac{\left(T_{k}/ \sqrt{k}\right)^2}{\Lambda''\left(\eta\right)}}} ds-\int_0^{T}e^{\lambda_0 s}\frac{e^{\lambda_1 s}}{e^{\lambda_1 s}-\eta} ds,
\end{align*}
where the last line is less than zero and only depends on $\beta_3$ by recalling that $T_{k}=\frac{\sqrt{k}}{24\beta_3}$. From Lemma \ref{3rd_moment_bound}, we know that $\beta_3$ is bounded for sufficiently large $n$ and sufficiently small $\epsilon_1$. Therefore, for any $\epsilon>0$, by \eqref{AB_limit}, we have
\begin{align*}
& 1+\frac{\int_0^{T}e^{\lambda_0 s}\frac{AB}{\sqrt{\left(AB-B+1\right)^2+2\left(AB-B\right)C}} ds-\int_0^{T}e^{\lambda_0 s}\frac{AB}{AB-B+1} ds}{\int_0^{\zeta_n-y}e^{\lambda_0 s}\frac{AB}{AB-B+1} ds}\\
\le & 1+\sup_{\eta\in \left[\eta_1, \eta_2\right]}\frac{\int_0^{T}e^{\lambda_0 s}\frac{e^{\lambda_1 s}}{\sqrt{\left(e^{\lambda_1 s}-\eta\right)^2 +\frac{\left(T_{k}/ \sqrt{k}\right)^2}{\Lambda''\left(\eta\right)}}} ds-\int_0^{T}e^{\lambda_0 s}\frac{e^{\lambda_1 s}}{e^{\lambda_1 s}-\eta} ds}{\int_0^{\infty} e^{\lambda_0 s}\frac{e^{\lambda_1 s}}{e^{\lambda_1 s}-\eta_1}ds}+\epsilon,
\end{align*}
for sufficiently large $n$ and sufficiently small $\epsilon_1$. The desired result then follows.
\qed
\end{proof}

\subsection*{Acknowledgements}
The authors would like to thank Xuanming Zhang for helpful comments on the draft. The work of KL was partially supported by NSF award CMMI 2228034.

\newpage

\bibliographystyle{plain}
\bibliography{ref}

\newpage

\vspace{0.5in}

\noindent Kevin Leder\\
Department of Industrial and Systems Engineering\\
University of Minnesota\\
Minneapolis, MN 55455, USA \\
kevin.leder@isye.umn.edu
\\

\noindent Zicheng Wang \\
School of Data Science\\
The Chinese University of Hong Kong, Shenzhen\\
Shenzhen, Guangdong 518172, China \\
wangzicheng@cuhk.edu.cn


\end{document}


\maketitle

\begin{lemma}
$\EE\left[\left(Z\left(t\right)\right)^4\right]\sim 24e^{4\lambda_1 t}$.
\end{lemma}

We follow the argument in https://math.stackexchange.com/questions/2795940/expected-value-of-a-simple-birth-process. Let $P_{n+1}\left(t\right)=e^{-\lambda_1 t}\left(1-e^{-\lambda_1 t}\right)^n$. Then
\begin{align*}
\EE\left[Z\left(t\right)\right]&=\sum\limits_{n=1}^{\infty}nP_n\left(t\right)\\
&=e^{\lambda_1 t}.
\end{align*}
\begin{align*}
\EE\left[\left(Z\left(t\right)\right)^2\right]&=\sum\limits_{n=1}^{\infty}n^2P_n\left(t\right)\\
&=e^{2\lambda_1 t}\left(2-e^{-\lambda_1 t}\right).
\end{align*}
\begin{align*}
\EE\left[\left(Z\left(t\right)\right)^3\right]&=\sum\limits_{n=1}^{\infty}n^3P_n\left(t\right)\\
&=e^{3\lambda_1 t}\left(6-6e^{-\lambda_1 t}+e^{-2\lambda_1 t}\right).
\end{align*}
\begin{align*}
\EE\left[\left(Z\left(t\right)\right)^4\right]&=\sum\limits_{n=1}^{\infty}n^4P_n\left(t\right)\\
&=e^{4\lambda_1 t}\left(24-36e^{-\lambda_1 t}+14e^{-2\lambda_1 t}-e^{-3\lambda_1 t}\right).
\end{align*} \qed

\begin{lemma}
\begin{align*}
\lim\limits_{n\rightarrow \infty}\frac{1}{n^{1-\alpha}}\log\PP\left(Z_1^n\left(\zeta_n-y\right)>\left(1+\epsilon\right)n\right)<\lim\limits_{n\rightarrow \infty}\frac{1}{n^{1-\alpha}}\log\PP\left(Z_1^n\left(\zeta_n-y\right)>n\right).
\end{align*}
\end{lemma}

We know that
\begin{align*}
& \frac{1}{n^{1-\alpha}}\log \EE\left[e^{\frac{\lambda_1\theta}{r_1}e^{-\lambda_1\left(\zeta_n-y\right)} Z_1^n\left(\zeta_n-y\right)}\right]\\
=& \frac{1}{n^{1-\alpha}}\log \EE\exp\left(\mu n^{1-\alpha}\int_{0}^{\zeta_n-y}e^{-rs}\left(\phi_{\zeta_n-y}\left(\frac{\lambda_1\theta}{r_1}e^{-\lambda_1\left(\zeta_n-y\right)}\right)-1\right)ds\right)\\
\rightarrow & \frac{\mu \lambda_1}{r_1}\int_{0}^{\infty}\frac{\theta e^{-rs}}{e^{\lambda_1 s}-\theta}ds.
\end{align*}
By Gartner-Ellis theorem, we have
\begin{align*}
& \lim\limits_{n\rightarrow \infty}\frac{1}{n^{1-\alpha}}\log\PP\left(Z_1^n\left(\zeta_n-y\right)>\left(1+\epsilon\right)n\right)\\
= & \lim\limits_{n\rightarrow \infty}\frac{1}{n^{1-\alpha}}\log\PP\left(n^{\alpha-1}\frac{\lambda_1\theta}{r_1}e^{-\lambda_1\left(\zeta_n-y\right)} Z_1^n\left(\zeta_n-y\right)>n^{\alpha-1}\left(1+\epsilon\right)n \frac{\lambda_1\theta}{r_1}e^{-\lambda_1\left(\zeta_n-y\right)}\right)\\
=& -\sup\limits_{\theta\in \left(0,1\right)}\left[\left(1+\epsilon\right)\frac{\theta \lambda_1 \mu e^{\lambda_1 y}}{r_1\left(\lambda_1-\lambda_0\right)}-\frac{\mu \lambda_1}{r_1}\int_{0}^{\infty}\frac{\theta e^{-rs}}{e^{\lambda_1 s}-\theta}ds\right]. 
\end{align*}\qed

\begin{lemma}
\begin{align*}
\lim\limits_{n\rightarrow \infty}\left|n^{1-\alpha}\EE\left[R_{n,y}\middle| \gamma_n < \zeta_n-y\right]-n^{1-\alpha}\EE\left[\tilde{R}_{n,y}\middle| \gamma_n < \zeta_n-y\right]\right|=0.
\end{align*}
\end{lemma}
\begin{align*}
& n^{1-\alpha}\EE\left[R_{n,y}\middle| \gamma_n < \zeta_n-y\right]-n^{1-\alpha}\EE\left[\tilde{R}_{n,y}\middle| \gamma_n < \zeta_n-y\right]   \\
= & n^{1-\alpha}\EE_{A_{n,y}}\left[\left\lvert\sum_{i=1}^{I_n\left(\zeta_n-y\right)}\left(\frac{X_{i,n}}{Z_1^n\left(\zeta_n-y\right)}\right)^2-\sum_{i=1}^{I_n\left(\zeta_n-y\right)}\left(\frac{X_{i,n}}{n}\right)^2\right\rvert\right]\\
= & n^{1-\alpha}\EE_{A_{n,y}}\left[\left\lvert\frac{\sum_{i=1}^{I_n\left(\zeta_n-y\right)}\left(X_{i,n}\right)^2}{n^2}\left(\frac{n^2}{Z_1^n\left(\zeta_n-y\right)}-1\right)\right\rvert\right].
\end{align*}
We can first condition on the event that $\{Z_1^n\left(\zeta_n-y\right)\le \left(1+\epsilon\right)n\}$ and then obtain an arbitrarily small number depending on $\epsilon$. We then condition on the event that $\{Z_1^n\left(\zeta_n-y\right)> \left(1+\epsilon\right)n\}$ and calculate the second moment. From the previous lemma, we know that
\begin{align*}
&  \frac{1}{n^{1-\alpha}}\log\PP\left(\frac{Z_1^n\left(\zeta_n-y\right)}{n}>1+\epsilon \right)\\
=& \frac{1}{n^{1-\alpha}}\log\PP\left(Z_1^n\left(\zeta_n-y\right)>\left(1+\epsilon\right)n \right)\\
\le & \frac{1}{n^{1-\alpha}}\log\frac{\EE\left[e^{\frac{\lambda_1\theta}{r_1}e^{-\lambda_1\left(\zeta_n-y\right)} Z_1^n\left(\zeta_n-y\right)}\right]}{\EE\left[e^{\frac{\lambda_1\theta}{r_1}e^{-\lambda_1\left(\zeta_n-y\right)} \left(1+\epsilon\right)n}\right]}\\
\approx & \frac{\mu \lambda_1}{r_1}\int_{0}^{\infty}\frac{\theta e^{-rs}}{e^{\lambda_1 s}-\theta}ds-\left(1+\epsilon\right)\frac{\theta \lambda_1 \mu e^{\lambda_1 y}}{r_1\left(\lambda_1-\lambda_0\right)}
\end{align*}
Therefore, for $n$ sufficiently large, and $\delta$ sufficiently small,
\begin{align*}
\PP\left(\frac{Z_1^n\left(\zeta_n-y\right)}{n}>1+\epsilon \right)\le \exp\left(n^{1-\alpha}\left(\frac{\mu \lambda_1}{r_1}\int_{0}^{\infty}\frac{\theta e^{-rs}}{e^{\lambda_1 s}-\theta}ds-\left(1+\epsilon\right)\frac{\theta \lambda_1 \mu e^{\lambda_1 y}}{r_1\left(\lambda_1-\lambda_0\right)}+\delta\right)\right)    
\end{align*}
We can obtain that
\begin{align*}
& \EE\left[\left(\frac{Z_1^n\left(\zeta_n-y\right)}{n}\right)^2\middle| \frac{Z_1^n\left(\zeta_n-y\right)}{n}>1+\epsilon  \right]\\
\le & \int\limits_{0}^{\left(1+\epsilon\right)^2}x dx+\int\limits_{\left(1+\epsilon\right)^2}^{\infty} x \frac{\PP\left(\frac{Z_1^n\left(\zeta_n-y\right)}{n}>\sqrt{x} \right)}{\PP\left(\frac{Z_1^n\left(\zeta_n-y\right)}{n}>1+\epsilon \right)}dx.
\end{align*}
We can show that
\begin{align*}
& \int\limits_{\left(1+\epsilon\right)^2}^{\infty} x \PP\left(\frac{Z_1^n\left(\zeta_n-y\right)}{n}>\sqrt{x} \right)dx\\
\le & \int\limits_{\left(1+\epsilon\right)^2}^{\infty} x \exp\left(n^{1-\alpha}\left(\frac{\mu \lambda_1}{r_1}\int_{0}^{\infty}\frac{\theta e^{-rs}}{e^{\lambda_1 s}-\theta}ds-\sqrt{x}\frac{\theta \lambda_1 \mu e^{\lambda_1 y}}{r_1\left(\lambda_1-\lambda_0\right)}+\delta\right)\right) dx\\
= & \int\limits_{1+\epsilon}^{\infty} 2y^3 \exp\left(n^{1-\alpha}\left(\frac{\mu \lambda_1}{r_1}\int_{0}^{\infty}\frac{\theta e^{-rs}}{e^{\lambda_1 s}-\theta}ds-y\frac{\theta \lambda_1 \mu e^{\lambda_1 y}}{r_1\left(\lambda_1-\lambda_0\right)}+\delta\right)\right) dy.
\end{align*}
We can take out a term
\begin{align*}
 \exp\left(n^{1-\alpha}\left(\frac{\mu \lambda_1}{r_1}\int_{0}^{\infty}\frac{\theta e^{-rs}}{e^{\lambda_1 s}-\theta}ds-\frac{\theta \lambda_1 \mu e^{\lambda_1 y}}{r_1\left(\lambda_1-\lambda_0\right)}+\delta\right)\right)
\end{align*}
which is very small compared to $n^{\alpha-1}$. We then calculate the rest 
\begin{align*}
& \int\limits_{\epsilon}^{\infty} 2y^3 \exp\left(-n^{1-\alpha}\left(y-1\right)A\right) dy\\
= & \int\limits_{1+\epsilon}^{\infty} 2\left(1+x\right)^3 \exp\left(-n^{1-\alpha}xA\right) dx,
\end{align*}
which is finite.

\begin{lemma}
$t_2-t_1< -\frac{1}{\lambda_1}\log\left(\frac{1}{2-\theta^*}\right)$ implies that $\kappa^*<\delta^*$
\end{lemma}
It suffices to show that
\begin{align*}
\frac{\int_{t_1}^{t_2}\frac{e^{\lambda_1 s}}{e^{\lambda_1 s}-\theta^*}e^{\lambda_0 s}ds}{\int_{t_1}^{t_2}e^{\lambda_0 s}ds}<\frac{\int_{t_1}^{t_2}\frac{e^{\lambda_1s}}{\left(e^{\lambda_1 s}-\theta^*\right)^2}e^{\lambda_0 s}ds}{\int_{t_1}^{t_2}e^{-\left(\lambda_1-\lambda_0\right)s}ds}. 
\end{align*}
We know that
\begin{align*}
\frac{\int_{t_1}^{t_2}\frac{e^{\lambda_1 s}}{e^{\lambda_1 s}-\theta^*}e^{\lambda_0 s}ds}{\int_{t_1}^{t_2}e^{\lambda_0 s}ds}<\frac{1}{1-e^{-\lambda_1 t_1}\theta^*},
\end{align*}
and
\begin{align*}
\frac{\int_{t_1}^{t_2}\frac{e^{\lambda_1s}}{\left(e^{\lambda_1 s}-\theta^*\right)^2}e^{\lambda_0 s}ds}{\int_{t_1}^{t_2}e^{-\left(\lambda_1-\lambda_0\right)s}ds}>\frac{1}{\left(1-e^{-\lambda_1 t_2}\theta^*\right)^2}.
\end{align*}
We have
\begin{align*}
& \left(1-e^{-\lambda_1 t_2}\theta^*\right)^2<1-e^{-\lambda_1 t_1}\theta^*\\
\Leftarrow & 2e^{-\lambda_1\left(t_2-t_1\right)}-\theta^*e^{-\lambda_1\left(t_2-t_1\right)-\lambda_1 t_2}>1\\
\Leftarrow & 2e^{-\lambda_1\left(t_2-t_1\right)}-\theta^*e^{-\lambda_1\left(t_2-t_1\right)}>1\\
\Leftarrow & t_2-t_1< -\frac{1}{\lambda_1}\log\left(\frac{1}{2-\theta^*}\right)
\end{align*}

\textbf{Some useful results}

\begin{align*}
\EE[Z_0^n(s)^2]=\frac{r_0+d_0}{r}ne^{-rs}(1-e^{-rs})+n^2e^{-2rs}
\end{align*}

\begin{align*}
e^{\lambda_1 \zeta_n} \sim \frac{\lambda_1 + r}{\mu} n^{\alpha}
\end{align*}

\begin{align*}
f_1(\theta)\rightarrow \int_0^{\infty}e^{\lambda_0 t}\frac{e^{\lambda_1 t}}{e^{\lambda_1 t}-\theta}
\end{align*}

\begin{align*}
f_2(\theta)\rightarrow \int_0^{\infty}e^{\lambda_0 t}\frac{e^{\lambda_1 t}}{(e^{\lambda_1 t}-\theta)^2}
\end{align*}

\begin{align*}
f_3(\theta)\rightarrow \int_0^{\infty}e^{\lambda_0 t}\frac{2e^{\lambda_1 t}}{(e^{\lambda_1 t}-\theta)^3}
\end{align*}

\begin{align*}
f_4(\theta)\rightarrow \int_0^{\infty}e^{\lambda_0 t}\frac{6e^{\lambda_1 t}}{(e^{\lambda_1 t}-\theta)^4}
\end{align*}


